\documentclass[11pt]{article}

\usepackage[scale=.8]{geometry}
\usepackage[notref,notcite]{}
\usepackage{graphicx}
\usepackage{amssymb}
\usepackage{amsmath}
\usepackage{comment}
\usepackage{hyperref}
\usepackage{enumerate}
\usepackage{amsthm}
\usepackage{amsmath}
\usepackage{float}
\usepackage{mathtools}
\usepackage{lineno}
\usepackage{bbm}
\usepackage{textcomp}
\usepackage{esint}
\usepackage{cleveref}
\usepackage{extarrows}
\usepackage{xcolor}
\definecolor{lanse}{RGB}{0,112,192}
\definecolor{zise}{RGB}{112,48,160} 
\definecolor{hongse}{RGB}{200,0,0}
\definecolor{gai}{RGB}{224,0,224}
\usepackage{lineno}

\makeatletter
\renewenvironment{proof}[1][\proofname]{%
	\par\pushQED{\qed}\normalfont%
	\topsep6\p@\@plus6\p@\relax
	\trivlist\item[\hskip\labelsep\bfseries#1\@addpunct{.}]%
	\ignorespaces
}{%
	\popQED\endtrivlist\@endpefalse
}
\makeatother

\numberwithin{equation}{section}

\newtheorem{theorem}{Theorem}
\newtheorem{proposition}[theorem]{Proposition}
\numberwithin{theorem}{section}
\newtheorem{corollary}[theorem]{Corollary}
\newtheorem{lemma}[theorem]{Lemma}

\makeatletter
\def\keywords{\xdef\@thefnmark{}\@footnotetext}
\makeatother

\renewcommand{\P}{\mathbb{P}}
\newcommand{\E}{\mathbb{E}}
\newcommand{\R}{\mathbb{R}}

\newcommand{\cF}{\mathcal F}

\newcommand{\eps}{\varepsilon}

\newcommand{\nn}{\nonumber}

\newcommand{\Extra}[1]{{\color{blue}#1}}
\renewcommand{\Extra}[1]{}
\newcommand{\la}{\langle}
\newcommand{\ra}{\rangle}

\newcommand{\hk}{\textcolor{black}}

\allowdisplaybreaks[4]

\begin{document}

\title{Occupation times for superprocesses in random environments}

\author{
Ziling Cheng
\footnote{Department of Mathematics and SUSTech International Center for Mathematics, Southern University of Science and Technology, Shenzhen, China}
\footnote{E-mail:\ chengzl@sustech.edu.cn}, 
Jieliang Hong
\footnote{Department of Mathematics, Southern University of Science and Technology, Shenzhen, China}
\footnote{E-mail:\ hongjl@sustech.edu.cn}
\footnote{Jieliang Hong is supported by the National Natural Science Foundation of China (Grant No. 12571150), and the Shenzhen Science and Technology Program (Grant No. RCBS20231211090635057).}
\ and Dan Yao
\footnotemark[3]
\footnote{corresponding author,\ E-mail:\ yaod@sustech.edu.cn}
}

\date{}

\maketitle

	
\noindent{\bf Abstract}\quad
	Let $X=(X_t, t\geq 0)$ be a superprocess in a random environment governed by a Gaussian noise \hk{$W=\{W(t, x),t\geq 0,x\in\mathbb{R}^d\}$} white in time and colored in space with correlation kernel $g$. We consider the occupation time process of the model starting from a finite measure. It is shown that the occupation time process of $X$ is absolutely continuous with respect to \hk{Lebesgue measure} in $d\leq 3$, whereas it is singular with respect to Lebesgue measure in $d\geq 4$. Regarding the absolutely continuous case in $d\leq 3$, we further prove that the associated \hk{density function is jointly H\"older continuous} based on the Tanaka formula and moment formulas, and derive the H\"older exponents with respect to the spatial variable $x$ and the time variable $t$. 
\vspace{0.3cm}



\section{Introduction}

	The Dawson-Watanabe superprocess (also known as the super-Brownian motion), which arises as a high-density limit of the critical branching particle system, has been studied by many authors since the pioneering work of Dawson \cite{D75} and Watanabe \cite{W68}. The distribution of a classical superprocess is mainly determined by two factors: the branching mechanism and the spatial motion. As variants of classical superprocesses, superprocesses in random environments have been interestingly studied by incorporating the influences of these two factors. The model in which the random environment affects the spatial motion was introduced and studied by Wang \cite{W97, W98}. Subsequently, Dawson et al. \cite{DVW00}, Li et al. \cite{LWXZ12}, and Hu et al. \cite{HLN13} further investigated the existence and smoothness properties of the density processes for such a model. On the other hand, the model in which the random environment influences the branching mechanism was studied by Mytnik \cite{M96}. \hk{Later, Sturm \cite{S03} considered a related branching mechanism whose variance approaches $0$. In addition, Hu, Nualart, and Xia \cite{HNX19} considered the model where the random environment affects both the branching mechanism and spatial motion and studied its density process.} 
	
	As an essential tool for studying the superprocess, the occupation time process has
also been studied by several authors. Iscoe \cite{I861} derived the Laplace functional of the
occupation time process and the associated stochastic partial differential equation for its
cumulant semigroup, providing a theoretical foundation for studies on the distributional
properties of occupation time processes. Thereafter, several scholars investigated
the distributional properties of occupation time processes for superprocesses, such as Blount and Bose \cite{BB00}, Dawson and Fleischmann \cite{DF94}, Iscoe \cite{I861} and \cite{I862},
Sugitani \cite{S89}, ect. For the superprocess in random environments proposed by Mytnik \cite{M96}, \hk{when this model starts from the Lebesgue measure, and each particle moves as a $d$-dimensional Brownian motion, Mytnik and Xiong \cite{MX07} showed the process will be extinct in a finite time when $d=1$ and $2$.} Chen, Ren, and Zhao \cite{CRZ15} showed that under some proper conditions, the process will converge weakly to a \hk{non-trivial} random measure when $d\geq3$. Recently, \hk{Fan, Hong and Xiong \cite{FHX24} derived the law of large numbers and the central limit theorem for the occupation time process of the model.} Meanwhile, there are still many interesting questions of the occupation time process for the model remaining unresolved.

	The fixed-time density process of the superprocess has been investigated over the past decades. Konno and Shiga \cite{KS88} (independently Reimers \cite{R89}) proved the existence of the density \hk{process} of the super-Brownian motion on $\mathbb{R}$ and provided the stochastic partial differential equation satisfied by the density \hk{process}. A similar study of a super-Brownian motion with a single point catalyst and its occupation time process can be seen in \cite{DF94}. Moreover, Sugitani \cite{S89} investigated the density process of the occupation time process and got some smoothness properties. The density process of the occupation time process for other superprocesses without \hk{random environments} can be found in \cite{BB00}, \cite{DF94}, \cite{K93}, \cite{RW08}, etc. There is relatively little research on the density processes of models in random environments. Kwon, Cho, and Kang \cite{KCK02} established the absolute continuity for the superprocess in \hk{random environments} given by Mytnik \cite{M96}, where the path of each particle is a one-dimensional Feller process.  The research on the density of the occupation time process for this model is currently lacking. 

	In this paper, we focus on the occupation time process of the superprocesses in random environments introduced by Mytnik \cite{M96}. We shall establish the existence and smoothness properties of the density processes for the occupation times, where the model begins with a finite measure. To formally state our model, we introduce the necessary notations. Let $B(\mathbb{R}^d)$ be the space of Borel functions on $\mathbb{R}^d$. Let $C_b^2(\mathbb{R}^d)$ be the space of bounded continuous functions on $\mathbb{R}^d$ with bounded continuous derivatives up to order $2$. Let $C_c(\mathbb{R}^d)$ be the space of continuous functions with compact support. We use the superscripts ``$+$'' to denote the subsets of nonnegative elements, e.g., $B^+(\mathbb{R}^d)$, $C_b^{2,+}(\mathbb{R}^d)$. Let $M_F(\mathbb{R}^d)$ denote the space of finite measures on $\mathbb{R}^d$ equipped with the topology of weak convergence. For any function $f$ on $\mathbb{R}^d$ and any measure $\mu\in M_{F}(\mathbb{R}^d)$, define
		\begin{align*}
			&\la \mu, f \ra = \mu(f) := \int_{\mathbb{R}^d} f(x) \mu(dx),\\
			&\mu f(x):=\int_{\mathbb{R}^d} f(x-y) \mu(dy).
		\end{align*}  
	
	Let $X=(X_t,t\geq 0)$ be a superprocess in a random environment with covariance function $g(x,y)$ defined on some complete filtered probability space $(\Omega,\mathcal{F},\mathcal{F}_t,\mathbb{P})$ such that $X$ satisfies the following martingale problem:
	\begin{align}\label{eMP}
		(MP)_{X_0}:\quad & \text{For any } \phi\in C_b^2(\mathbb{R}^d), M_t(\phi)=X_t(\phi)-X_0(\phi)-\int_0^t X_s(\frac{\Delta}{2}\phi)ds \nonumber\\
		& \text{is a continuous $(\mathcal{F}_t)$-martingale with} \\
		& \langle M(\phi) \rangle_t = \int_0^t X_s(\phi^2) ds + \int_0^t ds \int_{\mathbb{R}^d} \int_{\mathbb{R}^d} g(u,v) \phi(u) \phi(v) X_s(du) X_s(dv).\nonumber
	\end{align}
	The uniqueness of the solution of $(MP)_{X_0}$ was established by Mytnik \cite{M96}.  Throughout this paper, we assume that 
	$$\|g\|_{\infty}=\sup_{x,y\in\mathbb{R}^d}|g(x,y)|<\infty.$$ 
	For $X_0=\mu\in M_F(\mathbb{R}^d)$, we denote by $\mathbb{P}_{\mu}$ the law of $X$, and the corresponding expectation is $\mathbb{E}_{\mu}$. The occupation time process $Y=(Y_t, t\geq 0)$ is defined by
	\begin{align*}
		Y_t(\cdot):=\int_0^t X_s(\cdot) ds.	
	\end{align*}
	Since $X$ is a measure-valued process, so is $Y$. 
	The absolute continuity of $X$ with respect to Lebesgue measure in $d=1$ was proved by Kwon, Cho, and Kang \cite{KCK02}. As an integral of the original process $X$, the occupation time $Y$ will exhibit better smoothness. It is easy to check the singularities of $X$ and $Y$ with respect to Lebesgue measure in $d\ge 2$ and $d\ge 4$, respectively. We shall give the proofs in Appendix \ref{appA}.
	
	\begin{theorem}\label{th.s.}
		Let $\mu\in M_F(\mathbb{R}^d)$ with $d\geq 2$. With $\mathbb{P}_{\mu}$-probability one, $X_t$ is singular with respect to Lebesgue measure for almost every $t>0$.
	\end{theorem}
	
	\begin{theorem}\label{th.s..}
		Let $\mu\in M_F(\mathbb{R}^d)$ with $d\geq4$. With $\mathbb{P}_{\mu}$-probability one, $Y_t$ is singular with respect to Lebesgue measure for almost every $t\geq 0$.
	\end{theorem}
	
	Therefore, we are interested in proving the existence and the smoothness properties of the density process of $Y$ for $d\leq 3$.
	Before stating our results, we introduce the following Green function and recall the definitions of the \hk{local} H\"older continuity. For any $x\in \R^d$, define
	\begin{align}\label{de1.2}
		g_d(x)=
		\begin{cases}
			1, &\text{ if } d=1;\\	
			\log^+(1/|x|), &\text{ if } d=2;\\
			|x|^{-1}, &\text{ if } d=3.
		\end{cases}
	\end{align}

	\noindent $\mathbf{Definition.}$ A function $f(x)$ on $\R^d$ is said to be \textit{locally $\gamma$-H\"older continuous}, if for any compact set $K$ in $\R^d$, there exists some $C_K>0$ such that
	\[|f(x)-f(y)|\leq C_K|x-y|^\gamma, \quad \forall\ x, y \in K.\]
	We refer to $\gamma>0$ as the H\"older exponent.
	
	\noindent $\mathbf{Definition.}$ A function $f(t,x)$ defined on $[0,\infty)\times \R^d$ is said to be \textit{locally jointly H\"older continuous} if there \hk{are} some $\alpha_1, \alpha_2>0$ such that for any compact set $K$ in $[0,\infty)\times \R^d$, there exists some $C_K>0$ such that
		\begin{align*}
			|f(t, x)-f(s, y)|\leq C_K\Big[|t-s|^{\alpha_1}+|x-y|^{\alpha_2}\Big], \quad \forall\  (t, x), (s, y) \in K.
		\end{align*}
	\hk{We refer to $\alpha_1>0$ and $\alpha_2>0$ as the H\"older exponents in $t$ and in $x$,  respectively.}

	
	\begin{theorem}\label{t1}
		Let $d\leq 3$ and $\mu\in M_F(\R^d)$. When $d=2$ or $3$, we assume that 
		\begin{align}\label{de.1.5}
			\mu g_d(x)\ \text{is continuous with respect to $x\in\mathbb{R}^d$}.
		\end{align}
		
		\begin{itemize}
			\item[\rm{(i)}] \hk{There exists a family of nonnegative random variables $\{Y(t,x): t\geq 0, x\in \R^d\}$ such that for every $\phi \in C_c(\R^d)$ and $t\geq 0$,} with $\P_\mu$-probability one,
				\begin{align}\label{de2}
					\la Y_t, \phi \ra=\int_{\mathbb{R}^d}Y(t,x) \phi(x) dx.
				\end{align}
			
			\item[\rm{(ii)}] Furthermore, suppose that when $d=2$ or $3$, there is some $0< \gamma<1 \wedge (2-\frac{d}{2})$ such that
			\begin{align}\label{de1.3}
			\mu  g_d(x)  \text{ is locally $\gamma$-H\"older continuous.}
		\end{align}
			Then for all $d\leq 3$, with $\P_\mu$-probability one, 
			\begin{itemize}
				\item[\rm{(a)}] $Y(t,x)$ is locally jointly H\"older continuous in $t> 0$ and $x\in \R^d  $;
				
				\item[\rm{(b)}] 
			the H\"older exponent of $Y(t,x)$ 
			is arbitrarily close to
		\begin{equation*}
			\left\{
			\begin{aligned}
				&1\
			\text{in}\ x\ \text{and}\ 1/2\ \text{in}\ t,\quad \text{if}\ d=1; \\
				&\gamma\
			\text{in}\ x\ \text{and}\ 1/4\ \text{in}\ t,\quad \text{if}\ d=2;\\
				&\gamma\
			\text{in}\ x\ \text{and}\ 1/12\ \text{in}\ t,\quad \text{if}\ d=3.
			\end{aligned}
			\right.
		\end{equation*}
		\end{itemize}	
		\end{itemize}
	\end{theorem}

	To prove the above theorem, we need to introduce the conditional Laplace transform of the occupation time process $Y_t$. Let \hk{$W=\{W(t, x), t\geq0, x\in\mathbb{R}^d\}$} be a Gaussian noise white in time and colored in space with covariance function $g(x, y)$ such that 
	\begin{align*}
		\mathbb{E}_{\mu}[W(t, x) W(s, y)] = g(x, y) (s\wedge t).	
	\end{align*}
	Fan, Hong, and Xiong \cite{FHX24} get the following conditional Laplace transform, which is an easy consequence of Mytnik and Xiong \cite[Theorem 2,15]{MX07}. For any $\mu\in M_{F}(\mathbb{R}^d)$ and $\phi, f\in C_b^{2,+}(\mathbb{R}^d)$,
		\begin{align}\label{clat}
			\E_\mu^W \big[e^{- X_t(f)-Y_t(\phi) }\big]= e^{-\la \mu, U^{f, \phi}(t,\cdot)\ra},
		\end{align}
		where $\E_\mu^{W}$ denotes the conditional expectation of $\E_\mu$ given $W$ and $U^{f, \phi}\geq 0$ is the solution to the following SPDE: 
		\begin{align*}
			\hk{U^{f, \phi}(t, x)= f(x)+ t \phi(x)+ \int_0^t \frac{\Delta}{2} U^{f, \phi}(s, x) ds - \int_0^t \frac{1}{2} \big( U^{f, \phi}(s, x) \big)^2 ds + \int_0^t U^{f, \phi}(s, x) W (ds, x).}
		\end{align*}

	Using \eqref{clat}, the existence of the density process $Y(t, x)$ follows similarly to that of Sugitani \cite{S89} by calculating the second moment. The main difficulty is that an additional term, containing the noise $W$ from the random environment, is present. For the same reason, we are unable to use the methods used in Sugitani \cite{S89} to establish the recursive relationships regarding time differences and spatial differences about $Y(t, x)$. We overcome this difficulty by using the Tanaka formula of $Y(t, x)$ and the Green function representation of the original superprocess \mbox{$X=(X_t, t\geq0)$} given below.  Set 
	\begin{align}\label{e10.21}
		p_t^x(y)=p_t(x,y)=  \big(\frac{1}{2\pi t}\big)^{d/2}  e^{-\frac{|y-x|^2}{2t}}
	\end{align}
	to be the transition density of the $d$-dimensional standard Brownian motion. 
	Furthermore, set 
	\begin{align}\label{q}
		P_t f(x)=\int_{\mathbb{R}^d} p_t(x,y) f(y) dy,\quad x\in\mathbb{R}^d.
	\end{align}
	For any $\alpha\geq 0$, let
	\begin{align}\label{eb.6.1}
		g_\alpha^x(y)=\int_0^\infty e^{-\alpha t} p_t^x(y) dt,\quad x,y\in\mathbb{R}^d.
	\end{align}

	Let the symbol $\xrightarrow[]{L^1}$ stand for $L^1$-convergence. 
	
	\begin{theorem}\label{prop7.1}
		{\rm (Tanaka formula)}
		Let $d\leq 3$ and $\mu\in M_F(\mathbb{R}^d)$. We assume that \eqref{de.1.5} holds for $d=2,3$. Then there exists a measurable function $L_t^x(\omega):  \Omega\times [0,\infty)\times \R^d \to [0,\infty)$ such that for any $t\geq0$ and $x\in\mathbb{R}^d$,
		\begin{align*}
			\int_0^t X_s (p_\eps^x)  ds \xrightarrow[]{L^1} L_t^x,\quad \text{ as }\ \eps\downarrow 0,
		\end{align*}
		and for all $\alpha\geq0\ (\alpha>0 \text{ if } d\leq 2)$,
		\begin{align*}
			L_t^x=X_0(g_{\alpha}^x)-X_t(g_{\alpha}^x)+\alpha \int_0^t X_s(g_{\alpha}^x) ds+M_t(g_{\alpha}^x)\quad \mathbb{P}_{\mu}\text{-a.s.},
		\end{align*}
		where $(M_t(g_{\alpha}^x))_{t\ge 0}$ is a continuous $(\cF_t)$-martingale 
		with 
		\begin{align*}
		\langle M(g_{\alpha}^x) \rangle_t = \int_0^t X_s((g_{\alpha}^x)^2) ds + \int_0^t ds \int_{\mathbb{R}^d} \int_{\mathbb{R}^d} g(u,v)g_{\alpha}^x(u) g_{\alpha}^x(v) X_s(du) X_s(dv).
	\end{align*}
	\end{theorem}

We may refer to $L_t^x$ as the {\it local time} of the superprocess. In fact, the local time is also a version of the density function of the occupation measure process $Y_t$.
 To see this,  we will show later in Corollary \ref{c4.2} (see \eqref{e4.4.0}) that for any $t\geq0$ and $a\in\mathbb{R}^d$,
\begin{align*}
	 Y_t(p_h^a)  \xrightarrow[]{L^2} Y(t,a),\quad \text{ as } h\downarrow 0.
\end{align*}
On the other hand, Theorem \ref{prop7.1}  implies that as $h\downarrow 0$,
\begin{align*}
	Y_t(p_h^a)=\int_0^t X_s(p_h^a) ds \xrightarrow[]{L^1} L_t^a. 
\end{align*}
Hence, we have $Y(t,a)=L_t^a,\ \mathbb{P}_{\mu}
$-a.s. It follows that
\begin{align} \label{ea2}
	Y(t,a) = X_0(g_{\alpha}^a)-X_t(g_{\alpha}^a)+\alpha \int_0^t X_s(g_{\alpha}^a) ds+M_t(g_{\alpha}^a),\quad \mathbb{P}_{\mu}\text{-a.s.}
\end{align}

Given the above, the proof of the continuity of $Y(t,x)$ can be done by studying the four terms on the right-hand side of \eqref{ea2}. We will present the moment formulas for $X_t$ and $Y_t$ in Section \ref{sec2}, and then apply Kolmogorov's continuity criterion by calculating high-order moments to obtain a continuous version of each term in \eqref{ea2}. When computing moments, the Burkholder-Davis-Gundy inequality and a generalized Gronwall inequality (see Lemma \ref{l4.3}) are used to deal with the cross terms about the noise $W$. \\

Nevertheless, there will be some problems when calculating the moments for $|X_t(g_{\alpha}^a)-X_s(g_{\alpha}^a)|$ when $t$ and $s$ are close to $0$. 
To deal with this issue, we will use the Green function representation of $X_t$ to study the continuity of $t\mapsto X_t(g_{\alpha}^a)$ inspired by Perkins \cite{Per91}. \\

Let $M$ be the martingale measure associated with the martingale in \eqref{eMP}. We first give an extensive definition of the martingale measure.  Let $\mathcal{P}$ denote the $\sigma$-field of $(\mathcal{F}_t)$-predictable sets in $[0,\infty)\times\Omega$. Define
	\begin{align*}
		\mathcal{L}^2&(M):=\bigg\{\psi:  [0,\infty)\times  \Omega\times\mathbb{R}^d\rightarrow \mathbb{R}:\ \psi\ \text{is}\ \mathcal{P}\times B(\mathbb{R}^d)\text{-measurable},\nn\\
		& \mathbb{E}_{\mu}\Big[\int_0^t X_s(\psi(s,\cdot)^2) ds +\! \int_0^t \!ds\!\int_{\mathbb{R}^d}\! \int_{\mathbb{R}^d} g(u,v)\psi(s,u) \psi(s,v) X_s(du) X_s(dv) \Big]<\infty,\ \forall\ t>0\bigg\}.
	\end{align*}

	\begin{lemma}
	Let $d\leq 3$ and $\mu\in M_F(\mathbb{R}^d)$. For any $\psi\in\mathcal{L}^2(M)$, the martingale
	\begin{align*}
		M_t(\psi) := \int_0^t \int_{\mathbb{R}^d} \psi(s,x) dM(s,x)
	\end{align*}
	is well defined, and moreover,
	 \begin{align*}
		\langle M(\psi) \rangle_t = \int_0^t X_s(\psi(s,\cdot)^2) ds + \int_0^t ds \int_{\mathbb{R}^d} \int_{\mathbb{R}^d} g(u,v)\psi(s,u) \psi(s,v) X_s(du) X_s(dv).
	\end{align*}
	\end{lemma}
	The proof of the above lemma follows similarly to the proof of Proposition II.5.4 in Perkins \cite{P02} and is therefore omitted.\\

	 We present the following Green function representation.
	 
	\begin{theorem}\label{cor3.4*}
	{\rm (Green function representation)}
	Let $d\leq 3$ and $\mu\in M_F(\mathbb{R}^d)$. Suppose that \eqref{de.1.5} holds for $d=2,3$. For every $t\geq0$, $a\in\mathbb{R}^d$ and $\alpha\geq0\ (\alpha>0 \text{ if } d\leq 2)$, we have
	\begin{align}\label{gfr}
		X_t(g^a_\alpha) = X_0(P_t g^a_\alpha) + \int_0^t\int_{\mathbb{R}^d} P_{t-s}g^a_\alpha(x) dM(s, x),\quad \mathbb{P}_{\mu}\text{-a.s.},	
	\end{align}
	where $M$ is the martingale measure associated with $X$.
\end{theorem}
	
	The proofs of Theorems \ref{prop7.1} and \ref{cor3.4*} shall be given in Section \ref{sec5} and Appendix \ref{appg}, respectively.\\

	\noindent \textbf{Convention on Constants.} Constants whose value is unimportant and may change from place to place. All these constants may depend on the dimension $d$, the covariance function $g$, and the initial measure $\mu$. All these parameters, $d,  g, \mu$, will be fixed before picking the constants $C$. Other important dependence on time index (e.g., $t$, $T$, etc), moment index (e.g., $n$, $N$), test function $\phi$, H\"older exponent $\gamma$,  space index $a$, and so on will be mentioned at the subscript of constants $C$.\\
	
	
	\noindent \textbf{Organization of the paper.}
	In Section \ref{sec2}, we state all the moment formulas of $X$ and $Y$. Furthermore, we utilize the Burkholder-Davis-Gundy inequality and a generalized Gronwall inequality to provide some moment bounds for $X$ and $Y$. Section \ref{sec3} is then devoted to proving the existence of the density of $Y$ in $d\le 3$ by using the second moment formula. In Section \ref{sec5} we derive a Tanaka formula, which provides a representation of the density.  Combining this representation with moment bounds, we proceed to estimate the spatial and time differences of the density in Section \ref{sec6}. Based on these estimates, we construct the local joint H\"older continuity of the density and then finish the proof of Theorem \ref{t1} at the end of Section \ref{sec6}. In particular, an essential tool used in estimating the time difference, called the Green function representation, is proved in Appendix \ref{appg}. As a complement to the existence of the density, the singularities of $X$ in $d>1$ and $Y$ in $d>3$  are proved in Appendix \ref{appA}.\\
	
	
	Throughout the remainder of this paper, we fix the initial measure $\mu\in M_{F}(\mathbb{R}^d)$. For simplicity, we write $\mathbb{P}=\mathbb{P}_{\mu}$ and $\mathbb{E}=\mathbb{E}_{\mu}$ when there is no ambiguity. We write \hk{$\int$ for $\int_{\R^d}$} to ease notation.

	\section{Moment formulas and bounds}\label{sec2}
	
	In this section, we will provide some moment formulas for $Y$ and $X$.  Before doing so, we give some bounds on the moments of $Y_t(\cdot)$ and $X_t(\cdot)$ in these two cases. We start with a proper estimate.
	
	\begin{lemma}\label{lem5.1'}
		Let $d\ge 2$, for any $0\le \gamma<2$, there exists a constant $C_\gamma>0$ such that for any $t\ge 0$ and $x,z\in \mathbb{R}^d$,
		\begin{align*}
\int_0^t d s \int p_{t-s}(x, y) \frac{1}{|z-y|^\gamma} d y \leq C_{\gamma}\cdot t^{1-\frac{\gamma}{2}}.
\end{align*}
	\end{lemma}
	
	\begin{proof}
		Note that
		\begin{align*}
\int p_{t-s}(x, y) \frac{1}{|z-y|^\gamma} d y = \int p_{t-s}(z-x, y) \frac{1}{|y|^\gamma} d y,
\end{align*}
and $1/r^\gamma$ is a decreasing non-negative function on $(0,\infty)$. It follows from Lemma 3.6 of \cite{S89} that for any $x,z\in\mathbb{R}^d$,
\begin{align}\label{8.1}
 \int p_{t-s}(z-x, y) \frac{1}{|y|^\gamma} d y
 &\le \int p_{t-s}(0, y) \frac{1}{|y|^\gamma} d y \nn\\
 &=\Big(\frac{1}{2 \pi(t-s)}\Big)^{d / 2} \int e^{-\frac{|y|^2}{2(t-s)}} \frac{1}{|y|^\gamma} d y\nn\\
 &=C\Big(\frac{1}{2 \pi(t-s)}\Big)^{d / 2} \int_0^{\infty} e^{-\frac{r^2}{2(t-s)}} r^{d-\gamma-1} d r\nn\\
 &=C\Big(\frac{1}{2 \pi}\Big)^{d / 2} (t-s)^{-\gamma/2}\int_0^{\infty} e^{-\frac{r^2}{2}} r^{d-\gamma-1} d r\nn\\
 &\leq C_{\gamma} (t-s)^{-\gamma/2},
\end{align}
where $\int_0^{\infty} e^{-\frac{r^2}{2}} r^{d-\gamma-1} d r<\infty$ as we note $d-\gamma-1>-1$ for $d\ge 2>\gamma$. Hence
\begin{align*}
\int_0^t d s \int p_{t-s}(x, y) \frac{1}{|z-y|^\gamma} d y 
\le C_{\gamma}\int_0^t(t-s)^{-\gamma/2}ds
\le C_{\gamma}  t^{1-\frac{\gamma}{2}}.
\end{align*}
The proof is now complete.  
	\end{proof}
	
	\subsection{Moment formulas of $Y_t$}
	
	\hk{For any $\lambda\geq0$ and $\phi \in C_b^{2,+}(\R^d)$}, based on the conditional Laplace transform \eqref{clat}, by setting $U^{\lambda\phi}(t,x):=U^{0,\lambda\phi}(t, x)$ we obtain the conditional Laplace transform of the occupation time process $Y=(Y_t, t\geq 0)$ as follows. 
	\begin{align}\label{e10.1}
		 \E_\mu^W \big[e^{-\lambda Y_t(\phi)}\big] = e^{-\la \mu, U^{\lambda\phi}(t,\cdot)\ra},
	\end{align}
	where $\E_\mu^{W}$ is the conditional expectation of $\E_\mu$ given ${W}$ and $U^{\lambda\phi}\geq 0$ is the solution to the following SPDE: 
	\begin{align}\label{e10.3}
		U^{\lambda\phi}(t,x)=t \lambda \phi(x)+\int_0^t \frac{\Delta}{2} U^{\lambda\phi}(s,x) ds-\int_0^t \frac{1}{2} (U^{\lambda\phi}(s,x))^2 ds+\int_0^t   U^{\lambda\phi}(s,x) W(ds,x).
	\end{align}
	In the above, $W_t$ is a Guassian field satifying $\la W(\cdot, x), W(\cdot, y) \ra_t =g(x,y) t$.
	Set 
	\begin{align}\label{q1}
		Q_t f(x)=\int_0^t P_s f(x) ds,
	\end{align}
	where $P_t$ is given by \eqref{q}.
	We may rewrite \eqref{e10.3} as
	\begin{align}\label{e9.1}
		U^{\lambda\phi}(t,x)=  \lambda  Q_t \phi(x)  &-\frac{1}{2}\int_0^t  ds \int p_{t-s}(x ,z) (U^{\lambda\phi}(s,z))^2  dz\nn \\
		&+\int_0^t  \int p_{t-s}(x,z)  U^{\lambda\phi}(s,z) W(ds,z) dz.
	\end{align}
	For each $n\geq 1$, define
	\begin{align}\label{e2.1}
		V^{\phi}_n(t, x):=(-1)^{n-1}\frac{\partial^n}{\partial \lambda^n} U^{\lambda \phi}(t,x)\Big|_{\lambda=0}.
	\end{align}
	By differentiating \eqref{e9.1} with respect to $\lambda$, \hk{it is easy to check that $V^{\phi}_1(t,x)$} satisfies 
	\begin{align}\label{e10.5}
		V^{\phi}_1(t,x)= Q_{t} \phi(x) +\int_0^t  \int p_{t-s}(x,z)  V^{\phi}_1(s, z) W(ds,z) dz.
	\end{align}
	For $n\geq 2$,  by iteratively differentiating \eqref{e9.1} with respect to $\lambda$, we get that \hk{$V^{\phi}_n(t,x)$} satisfies  
	\begin{align}\label{e6.11}
		V^{\phi}_n(t,x)=  &\sum_{k=1}^{n-1} \binom{n-1}{k}   \int_0^t ds \int p_{t-s} (x,z)V^{\phi}_{n-k}(s, z)V^{\phi}_{k}(s, z) dz\nn\\
		&+\int_0^t  \int p_{t-s}(x,z)  V^{\phi}_{n}(s,z) W(ds,z) dz.
	\end{align}
	
	The following two results give some useful moment formulas of $Y$. We omit the proofs of them since the arguments
are similar to Lemmas 3.1-3.2 and Corollary 3.3 of \cite{FHX24}.
	
	\begin{proposition}\label{l2.1}
		Let $d\ge 1$. For any \hk{$t\ge 0$} and $\phi\in C_b^{2,+}(\R^d)$, we have
		\begin{align}\label{e10.22}
			\E_\mu^W [Y_t(\phi)^n]= L^{\phi}_{n}(t), \quad \forall\ n\geq 1,
		\end{align}
		where
		\begin{align}\label{e10.11}
			L^{\phi}_0(t)=&1, \quad L^{\phi}_1(t)=   \la \mu, V^{\phi}_1(t,\cdot)\ra, \nn\\
			L^{\phi}_n(t)=&\sum_{k=0}^{n-1} \binom{n-1}{k}   \la \mu, V^{\phi}_{n-k}(t,\cdot) \ra L^{\phi}_k(t) , \quad \forall\ n\geq 2.
		\end{align}
	\end{proposition}
	
	\begin{corollary}\label{c2.1}
		Let $d\ge 1$. For every \hk{$t\ge 0$} and $\phi, \psi \in B^+(\R^d)$, we have
		\begin{align} \label{e3}
			\E_\mu [Y_t(\phi)]= \la \mu, Q_t\phi\ra,
		\end{align}
		and
		\begin{align} \label{e10.7}
			\E_\mu [Y_t(\phi)Y_t(\psi)]= \int &\int V^{\phi, \psi}_t(x,y)\mu(dx) \mu(dy)\nn\\
			&+\int_0^t ds \int \mu(dx)\int p_{t-s}(x,y) V^{\phi, \psi}_s(y,y) dy,
		\end{align}
		where
		\begin{align} \label{e10.23}
			V^{\phi, \psi}_t(x,y):=\Pi_{(x,y)}\Big\{\int_0^t  \Big[\phi(B_s) Q_{t-s}\psi(\tilde{B}_s)+\psi(\tilde{B}_s) Q_{t-s}\phi({B}_s)\Big] e^{\int_0^s g(B_u, \tilde{B}_u) du} ds\Big\}
		\end{align}
		and $B_t, \tilde{B}_t$ are independent \hk{Brownian motions} starting respectively  from $x,y\in \R^d$ under $\Pi_{(x,y)}$.  
	\end{corollary}
	
	\subsection{Moment bounds of $Y_t$}
	
		 Based on the estimates of $V^{\phi}_n(t,x)$ in Appendix \ref{appb1}, we first give the moment bounds of $Y_t$ with respect to a general text function $\phi\in C_b^{2,+}(\mathbb{R}^d)$ as follows.	
	\begin{lemma}\label{lem5.5}
		Let $d\ge 1$. For any $n\ge 1$, $T>0$ and $\phi\in C_b^{2,+}(\R^d)$, there exists some constant $C_{n,T,\phi}>0$ such that for all $0\le t\le T$,
		  \begin{align*}
\mathbb{E}_\mu[Y_t(\phi)^n] \leq C_{n, T,\phi}.
\end{align*}
	\end{lemma}
	\begin{proof}
		Recall from Proposition \ref{l2.1} that
		\begin{align}\label{'y1}
\mathbb{E}_\mu[Y_t(\phi)^n]=\mathbb{E}[L^{\phi}_n(t)] ,
\end{align}
where $L^{\phi}_n(t)$ is given by (\ref{e10.11}). Fix $T>0$ and let $0\le t\le T$, for any $n\ge 1$, we shall prove that
\begin{align}\label{'y2}
\mathbb{E}_\mu[L^{\phi}_n(t)^N] \leq C_{n,N, T,\phi},\quad \forall\ N\ge 1.
\end{align}
The conclusion follows immediately from (\ref{'y1}) and the above with $N=1$. 
When $n=1$, for any $N\ge 1$ we have
\begin{align}\label{'5.16}
\mathbb{E}[L^{\phi}_1(t)^N] =\mathbb{E}[\langle \mu,V^{\phi}_1(t,\cdot)\rangle^N]&=\int \mu(d x_1) \cdots \int \mathbb{E}\bigg[\prod_{i=1}^N V^{\phi}_1(t, x_i)\bigg] \mu(d x_n)\nn\\
	&\le C_{N,T}\int \mu(d x_1) \ldots \int \prod_{i=1}^N Q_t \phi(x_i) \mu(d x_n)\nn\\
	&=C_{N, T}\bigg[\int Q_t \phi(x) \mu(d x)\bigg]^N\nn\\
	&\le C_{N, T}[\mu(1)G(\phi,T)]^N\le C_{N, T}\big[\mu(1)(1+T\|\phi\|_{\infty})\big]^N,
\end{align}
where $G(\phi,t)$ is defined by (\ref{g.1}) and the first inequality follows by Lemma \ref{lem5.1}. Hence (\ref{'y2}) holds for the case $n=1$. Assume that (\ref{'y2}) holds for all $1\le k\le n-1$ with some $n\ge 2$. Then, for the case $n$, by Cauchy-Schwarz's inequality, we get
\begin{align}\label{'y3}
\mathbb{E}[L^{\phi}_n(t)^N] &\le C_{n,N}\sum_{k=0}^{n-1}\mathbb{E}\big[\langle \mu,V^{\phi}_{n-k}(t,\cdot)\rangle^N L^{\phi}_{k}(t)^N\big]\nn\\
&\le C_{n,N}\sum_{k=0}^{n-1}\Big(\mathbb{E}\big[\langle \mu,V^{\phi}_{n-k}(t,\cdot)\rangle^{2N}\big]\Big)^{1/2} \Big(\mathbb{E}\big[L^{\phi}_{k}(t)^{2N}\big]\Big)^{1/2}\nn\\
&\le C_{n,N}\sum_{k=0}^{n-1}C_{k,N,T,\phi}\Big(\mathbb{E}\big[\langle \mu,V^{\phi}_{n-k}(t,\cdot)\rangle^{2N}\big]\Big)^{1/2} ,
\end{align} 
where we have used the induction hypothesis on the last line. When $k=n-1$, for any $0\le t\le T$ we get from (\ref{'5.16}) that
\begin{align*}
\mathbb{E}[\langle\mu, V^{\phi}_{n-k}(t,\cdot)\rangle^{2 N}]=\mathbb{E}[\langle\mu, V^{\phi}_1(t,\cdot)\rangle^{2 N}]=\mathbb{E}[L^{\phi}_1(t)^{2 N}] \leq C_{N,T,\phi}.
\end{align*}
When $0\le k\le n-2$, by H\"older's inequality and Lemma \ref{lem5.3} we have
\begin{align*}
\mathbb{E}[\langle\mu, V^{\phi}_{n-k}(t,\cdot)\rangle^{2N}]
&\le \mu(1)^{2N-1}\int \mathbb{E}[V^{\phi}_{n-k}(t, x)^{2 N}] \mu(d x)\nn\\ 
&\le \mu(1)^{2N-1}\int C_{k,n,N,T}G(\phi,T)^{N2^{n-k}} \mu(d x)\nn\\
&= \mu(1)^{2N} C_{k,n,N,T}G(\phi,T)^{N2^{n-k}}.
\end{align*} 
Therefore (\ref{'y3}) becomes
\begin{align*}
\mathbb{E}[L^{\phi}_n(t)^N] 
&\le C_{n,N,T,\phi}<\infty.
\end{align*} 
The proof is complete by induction.
	\end{proof}
	
	Now we consider a special set of text functions $\phi_{a,\gamma}$ (on $\mathbb{R}^d$) defined by 
\begin{align}\label{phig}
	\phi_{a,\gamma}(x)=\frac{1}{|a-x|^{\gamma}},\quad 0<\gamma<d,\ a\in \mathbb{R}^d.
\end{align}
Since $\phi_{a,\gamma}$ is a non-negative measurable function, there exists a sequence $\{\phi_{a,\gamma,r}\}\subset C_b^{2,+}(\mathbb{R}^d)$ such that
\begin{align}\label{phir}
	\phi_{a,\gamma,r}\uparrow \phi_{a,\gamma},\quad \text{as}\ r\rightarrow\infty.
\end{align}
Based on Lemmas \ref{'lem5.14} and \ref{'lem5.17}, we immediately obtain the following results by monotone convergence.
	
	\begin{corollary}\label{'cor5.15}
		Let $d=2$. For every $n \geq 1$ and $T>0$, there exists some constant $C_{n, T, \gamma}>0$ such that for any $0 \leq t \leq T$ and $a \in \mathbb{R}^2$,
		\begin{align*}
			\mathbb{E}_{\mu}\big[Y_{t}(\phi_{a,\gamma})^{n}\big] \leq C_{n, T, \gamma}.
		\end{align*}
	\end{corollary}
	
	\begin{corollary}\label{'cor5.18}
		Let $d=3$. Suppose \eqref{de.1.5} holds for $d=3$. For every $n \geq 1$ and $T>0$, there exists some constant $C_{n, T, \gamma}>0$ such that for any $0 \leq t \leq T$ and $a\in\mathbb{R}^3$ with $|a|\leq T$,
		\begin{align*}
			\mathbb{E}_{\mu}\big[Y_{t}(\phi_{a,\gamma})^{n}\big] \leq C_{n, T, \gamma}.
		\end{align*}
	\end{corollary}

	\subsection{Moment formulas of $X_t$}
	
	For any $\phi \in C_b^{2,+}(\R^d)$, by \eqref{clat} and letting $\widetilde{U}^{\phi}(t, x):=U^{\phi,0}(t, x)$ yields the conditional Laplace transform of the original superprocess $X=(X_t, t\geq0)$, i.e.,
	\begin{align*} 
		\E_\mu^W \big[e^{- X_t(\phi) }\big]= e^{-\la \mu, \widetilde{U}^{\phi}(t,\cdot)\ra},
	\end{align*}
	where  $\widetilde{U}^{\phi}\geq 0$ is the pathwise unique solution to the following SPDE: 
	\begin{align*}
		\hk{\widetilde{U}^{\phi}(t,x)=f(x)+  \int_0^t \frac{\Delta}{2} \widetilde{U}^{\phi}(s,x)ds-\frac{1}{2}\int_0^t (\widetilde{U}^{\phi}(s,x))^2ds+\int_0^t \widetilde{U}^{\phi}(s,x) W(ds,x).}
	\end{align*}
	Based on \hk{this}, by repeating the calculations of \hk{Proposition \ref{l2.1} and} Corollary \ref{c2.1}, we get the following \hk{results}, which \hk{give} the $n$-th moment formulas of $X_t$.
	
	\begin{proposition}\label{l2.2}
		Let $d\ge 1$. For any \hk{$t\ge 0$} and $\phi\in C_b^{2,+}(\R^d)$, we have
		\begin{align*}
			\E_\mu^W [X_t(\phi)^n]= \widetilde{L}^{\phi}_{n}(t), \quad \forall\ n\geq 1,
		\end{align*}
		where
		\begin{align*}
			\widetilde{L}^{\phi}_0(t)=&1, \quad \widetilde{L}^{\phi}_1(t)=   \la \mu, \widetilde{V}^{\phi}_1(t,\cdot)\ra, \\
			\widetilde{L}^{\phi}_n(t)=&\sum_{k=0}^{n-1} \binom{n-1}{k}   \la \mu, \widetilde{V}^{\phi}_{n-k}(t,\cdot) \ra \widetilde{L}^{\phi}_k(t) , \quad \forall\ n\geq 2,\nn
		\end{align*}
		with 
		\begin{align}\label{ve10.5'}
		\widetilde{V}^{\phi}_1(t,x)&= P_{t} \phi(x) +\int_0^t  \int p_{t-s}(x,z)  \widetilde{V}^{\phi}_1(s, z) W(ds,z) dz,\nn\\
		\widetilde{V}^{\phi}_n(t,x)&= \sum_{k=1}^{n-1} \binom{n-1}{k}   \int_0^t ds \int p_{t-s} (x,z)\widetilde{V}^{\phi}_{n-k}(s, z)\widetilde{V}^{\phi}_{k}(s, z) dz\nn\\
		&\qquad+\int_0^t  \int p_{t-s}(x,z)  \widetilde{V}^{\phi}_{n}(s,z) W(ds,z) dz,\quad \forall\  n\geq 2.
	\end{align}
	\end{proposition}
	
	\begin{corollary}\label{cor3.2}
		Let $d\ge 1$. For any $t\geq0$ and $\phi, \psi \in B^+(\mathbb{R}^d)$, we have
		\begin{align*}
			\mathbb{E}_\mu[X_t(\phi)]=\langle\mu, P_t \phi\rangle,
		\end{align*}
		and 
		\begin{align*}
			\mathbb{E}_\mu[X_t(\phi) X_t(\psi)]=\int & \int \widetilde{V}_t^{\phi, \psi}(x, y) \mu(d x) \mu(d y) \\
			& +\int_0^t d s \int \mu(d x) \int p_{t-s}(x, y) \widetilde{V}_s^{\phi, \psi}(y, y) d y,
		\end{align*}
		where
		\begin{align*}
			\widetilde{V}_t^{\phi, \psi}(x, y):=\Pi_{(x, y)}\bigg\{\exp \bigg(\int_0^t g(B_s, \tilde{B}_s) d s\bigg) \phi(B_t) \psi(\tilde{B}_t)\bigg\}
		\end{align*}
		and $B_t, \tilde{B}_t$ are independent \hk{Brownian motions} starting respectively from $x, y \in \mathbb{R}^d$ under $\Pi_{(x, y)}$.
	\end{corollary}
	
		\subsection{Moment bounds of $X_t$}
	
\begin{lemma}\label{lemx5.1}
		Let $d\ge 1$. Given $T\geq 0$, there is some constant $C_{N,T}>0$ such that for all $0\le t\le T$,
		\begin{align*}
			\mathbb{E}_{\mu} \Big[ \big(\sup_{0\le s\le t}X_s(1)\big)^{2N} \Big] \le  C_{N,T}.
		\end{align*}
	\end{lemma}

	\begin{proof}
		By applying the martingale problem \eqref{eMP} with $\phi\equiv1$ we get that
		\begin{align*}
			X_t(1) =   X_0(1) + M_t(1),
		\end{align*}
		where 
		\begin{align*}
			\langle M(1) \rangle_t 
			&= \int_0^t X_s(1) ds + \int_0^t ds \int\!\!\int g(u,v)  X_s(du) X_s(dv).
		\end{align*}
		Then by Lemma \ref{lem5.5}, for any $0\le t\le T$ we easily obtain
		\begin{align*}
			\mathbb{E}_{\mu} \Big[ \big(\sup_{0\le s\le t}X_s(1)\big)^{2N} \Big]
			&\leq 2^{2 N}\mathbb{E}_\mu\big[X_0(1)^{2 N}\big]+2^{2 N} \mathbb{E}_\mu\Big[\big(\sup _{0\le s\le t} M_s(1)\big)^{2 N}\Big] \nn\\
			&\le (2 \mu(1))^{2 N}+2^{2 N}C_N \mathbb{E}_\mu\big[\langle M(1)\rangle_t^N\big]\nn\\
			&\leq C_N+C_N \mathbb{E}_\mu\Big[\Big(\int_0^t X_s(1) ds + \int_0^t ds \int\!\!\int g(u,v)  X_s(du) X_s(dv)
		\Big)^{N}\Big]\nn\\
		&\leq C_N+C_N \Big\{\mathbb{E}_\mu\big[Y_t(1)^{N}\big]+\|g\|_{\infty}^{N}\mathbb{E}_\mu\Big[\Big(\int_0^t X_s(1)^2 ds
		\Big)^{N}\Big]\Big\}\nn\\
		&\leq C_N+C_N \Big\{\mathbb{E}_\mu\big[Y_t(1)^{N}\big]+\|g\|_{\infty}^{N}t^{N-1}\mathbb{E}_\mu\Big[\int_0^t X_s(1)^{2N} ds
		\Big]\Big\}\nn\\
		&\leq C_N+C_N \Big\{\mathbb{E}_\mu\big[Y_t(1)^{N}\big]+t^{N-1}\int_0^t\mathbb{E}_\mu\Big[\big(\sup_{0\le r\le s}X_s(1)\big)^{2N}
		\Big]ds\Big\}\nn\\
		&\le C_{N,T}+C_{N,T} \int_0^t\mathbb{E}_\mu\Big[\big(\sup_{0\le r\le s}X_s(1)\big)^{2N}
		\Big]ds
		\end{align*}
		Therefore, by Gronwall's inequality, we have
		\begin{align*}
			\mathbb{E}_{\mu} \Big[ \big(\sup_{0\le s\le t}X_s(1)\big)^{2N} \Big]
			&\leq C_{N,T}\Big[1+C_{N,T}\int_{0}^t e^{C_{N,T}(t-s)}ds \Big]\nn\\
			&\le C_{N,T}e^{TC_{N,T}}=C_{N,T},
		\end{align*}
	as required.
 	\end{proof}

	Now we discuss a special set of text functions $\phi_{a,\gamma}$ (on $\mathbb{R}^d$) with the same form as (\ref{phig}), i.e.,
\begin{align}\label{xphig}
	\phi_{a,\gamma}(x)=\frac{1}{|a-x|^{\gamma}},\quad a\in \mathbb{R}^d.
\end{align}
We shall consider $0<\gamma<1$ when $d=2$ and $0<\gamma<5/2$ when $d=3$.
Similarly, there exists a sequence $\{\phi_{a,\gamma,r}\}\subset C_b^{2,+}(\mathbb{R}^d)$ such that
\begin{align}\label{xphir}
	\phi_{a,\gamma,r}\uparrow \phi_{a,\gamma},\quad \text{as}\ r\rightarrow\infty.
\end{align}
Based on Lemmas \ref{'lemb5.7} and \ref{'lemb5.10}, we immediately obtain the following results by monotone convergence.
	
	\begin{corollary}\label{'corb5.8}
		Let $d=2$ and $0<\gamma<1$. For every $n \geq 1$ and $0<\eps_0<T$, there exists some constant $C_{n, \eps_0, T, \gamma}>0$ such that for any $\eps_0 \leq t \leq T$ and $a\in\mathbb{R}^2$,
		\begin{align*}
			\mathbb{E}_{\mu}\big[X_{t}(\phi_{a,\gamma})^{n}\big] \leq C_{n, \eps_0, T, \gamma}.
		\end{align*}
	\end{corollary}
	
	\begin{corollary}\label{'corb5.11}
		Given $d=3$ and $0<\gamma<5/2$. For every $n \geq 1$ and $0<\eps_0<T$, there exists some constant $C_{n, \eps_0, T, \gamma}>0$ such that for any $\eps_0 \leq t \leq T$ and $a \in \mathbb{R}^3$,
		\begin{align*}
			\mathbb{E}_{\mu}\big[X_{t}(\phi_{a,\gamma})^{n}\big] \leq C_{n, \eps_0, T, \gamma}.
		\end{align*}
	\end{corollary}
	
	\section{Existence of the density}
	\label{sec3}
	
	In this section, we shall prove the existence of the density in $d\le 3$. Specifically, for any $h>0$ , $t\geq 0$ and $x\in \R^d$, we define
	\begin{align}\label{e2.13}
		Y_h(t,x):=\int p_h^x(y) Y_t(dy).
	\end{align}
	\hk{By using} the moment formulas of $Y_t$, we shall show that the limit of $Y_{h}(t,x)$ (as $h\downarrow 0$) exists. An application of Corollary \ref{c2.1} first gives the following results.
	
	\begin{proposition}\label{l4.1}
		Let $d\leq3$. 
		Then for any fixed \hk{$t\ge 0$} we have
		\begin{align}\label{e4.1} 
			\limsup_{h, h^{\prime} \downarrow 0}\int \E_\mu\big[Y_h(t, x) Y_{h^{\prime}}(t, x)\big] dx<\infty.
		\end{align}
	\end{proposition}
	
	\begin{proof}
		We use \eqref{e10.7} with $\phi(u)=p_h^x(u)$ and $\psi(u)=p_{h^{\prime}}^x(u)$ to get that for any \hk{$t\ge 0$}, $x\in\mathbb{R}^d$ and $h,h^{\prime}>0$,
		\begin{align}\label{e.4.14}
			\E_\mu[Y_h(t, x) Y_{h^{\prime}}(t, x)] & =\int\!\!\int V_t^{p_h^x, p_{h^{\prime}}^x} (u, v) \mu(d u) \mu(d v)+\int_0^t d s \int \mu(d u) \int p_{t-s}(u, v) V_{s}^{p_h^x,p_{h^{\prime}}^x} (v, v) d v \nonumber\\
			& :=I_1+I_2,
		\end{align}
		where 
		\begin{align*}
			 V_t^{p_h^x, p_{h^{\prime}}^x}(u, v)=&\Pi_{(u, v)}\bigg\{\int_0^t \Big[p_h^x(B_s) Q_{t-s} p_{h^{\prime}}^x(\tilde{B}_s) + p_{h^{\prime}}^x(\tilde{B}_s) Q_{t-s} p_{h}^x(B_s) \Big] e^{\int_0^s g(B_r, \tilde{B}_r) d r} d s\bigg\} \nonumber\\
			\leq& e^{t\|g\|_{\infty}} \Pi_{(u, v)} \bigg\{\int_0^t \Big[p_h^x(B_s) Q_{t-s} p_{h^{\prime}}^x(\tilde{B}_s) + p_{h^{\prime}}^x(\tilde{B}_s) Q_{t-s} p_{h}^x(B_s)\Big] d s\bigg\} . 
		\end{align*}
		Note that $Q_{t-s} p_{h^{\prime}}^x(u)=\int_0^{t-s} P_r p_{h^{\prime}}^x(u) d r=\int_0^{t-s} p_{r+h^{\prime}}^x(u) d r$, it then follows that
		\begin{align*}
			\Pi_{(u, v)} \bigg[\int_0^t p_h^x(B_s) Q_{t-s} p_{h^{\prime}}^x(\tilde{B}_s) ds\bigg] 
			&= \int_0^t d s \int_0^{t-s} \Pi_{(u, v)} \big[p_h^x(B_s) p_{r+h^{\prime}}^x(\tilde{B}_s) \big] d r \nonumber\\
			&= \int_0^t d s \int_0^{t-s} p_{s+h}^x(u) p_{s+r+h^{\prime}}^x(v) d r \nonumber\\
			&= \int_0^t p_{s+h}^x(u) d s \int_s^{t}  p_{r+h^{\prime}}^x(v) dr.
		\end{align*}
		Similarly, we have
		\begin{align*}
			\Pi_{(u, v)} \bigg[\int_0^t p_{h^{\prime}}^x(\tilde{B}_s) Q_{t-s} p_h^x(B_s) ds\bigg] = \int_0^t p_{s+h^{\prime}}^x(v) d s \int_s^{t}  p_{r+h}^x(u) dr.
		\end{align*}	
		Therefore, 
		\begin{align}\label{e4.3}
			V_t^{p_h^x, p_{h^{\prime}}^x}(u, v) \leq e^{t\|g\|_{\infty}} \bigg[ \int_0^t p_{s+h}^x(u) d s \int_s^{t}  p_{r+h^{\prime}}^x(v) dr + \int_0^t p_{s+h^{\prime}}^x(v) d s \int_s^{t}  p_{r+h}^x(u) dr \bigg].
		\end{align}
		Hence,
		\begin{align*}
			I_1 &\leq e^{t\|g\|_{\infty}} \int_{0}^t \mu(p_{s+h}^x) ds \int_s^t   \mu(p_{r+h^{\prime}}^x) dr + e^{t\|g\|_{\infty}} \int_{0}^t \mu(p_{s+h^{\prime}}^x) ds \int_s^t  \mu(p_{r+h}^x)  dr \nonumber\\
			&\leq \mu(1) e^{t\|g\|_{\infty}} \int_{0}^t  \mu(p_{s+h}^x+p_{s+h^{\prime}}^x)  ds \int_{s}^t r^{-d/2} dr,
		\end{align*}
		where the last inequality follows by $\mu(p_{r+h^{\prime}}^x) \leq \mu(1) r^{-d/2}$. We next use \eqref{e4.3} to bound $I_2$ by
		\begin{align}\label{e4.4}
			I_2 \leq &\int_{0}^t e^{s\|g\|_{\infty}}ds \int \mu(du) \int p_{t-s}(u,v)   dv\int_0^s p_{s_1+h}^x(v) d s_1 \int_{s_1}^{s}  p_{r+h^{\prime}}^x(v) dr \nonumber \\
			 &+ \int_{0}^t e^{s\|g\|_{\infty}}ds \int \mu(du) \int p_{t-s}(u,v)  dv\int_0^s p_{s_1+h^{\prime}}^x(v) d s \int_{s_1}^{s}  p_{r+h}^x(v) dr.
		\end{align}
		Since $p_{r+h^{\prime}}^x(v)\leq r^{-d/2}$, the first term on the right-hand side of \eqref{e4.4} is bounded by
		\begin{align}\label{be4.4}
			&e^{t\|g\|_{\infty}} \int_{0}^t ds \int \mu(du) \int p_{t-s}(u,v)  dv\int_0^s p_{s_1+h}^x(v) d s_1 \int_{s_1}^{s}  p_{r+h^{\prime}}^x(v) dr \nonumber\\
			&\leq e^{t\|g\|_{\infty}} \int_{0}^t ds \int_0^s \bigg( \int \mu(du) \int p_{t-s}(u,v)  p_{s_1+h}^x(v)   dv \bigg) d s_1 \int_{s_1}^{s}   r^{-d/2} dr  \nonumber\\
			&=e^{t\|g\|_{\infty}} \int_{0}^t ds \int_0^s \mu(p_{t-s+s_1+h}^x) d s_1 \int_{s_1}^{s}   r^{-d/2} dr,
		\end{align}
		By a similar estimation, the second term on the right-hand side of \eqref{e4.4} is bounded by
		\begin{align}\label{e4.5}
			e^{t\|g\|_{\infty}} \int_{0}^t ds \int_0^s \mu(p_{t-s+s_1+h^{\prime}}^x) d s_1 \int_{s_1}^{s}   r^{-d/2} dr.
		\end{align}
		Hence, in view of \eqref{e4.4}-\eqref{e4.5} we have
		\begin{align*}
			I_2 \leq e^{t\|g\|_{\infty}} \int_{0}^t ds \int_0^s \mu(p_{t-s+s_1+h}^x + p_{t-s+s_1+h^{\prime}}^x) d s_1 \int_{s_1}^{s}   r^{-d/2} dr.
		\end{align*}
		Based on the estimations of $I_1$ and $I_2$, and the fact
		\begin{align*}
			\int p_{s+h}^x (u)dx= \int p_{s+h^{\prime}}^x (u)dx = \int p_{t-s+s_1+h}^x (u)dx=\int p_{t-s+s_1+h^{\prime}}^x (u)dx =1,
		\end{align*}
		by \hk{\eqref{e.4.14}} and a simple calculation, we see that
		\begin{align}\label{e4.6}
			&\int \E_\mu\big[Y_h(t, x) Y_{h^{\prime}}(t, x)\big] dx\nn\\ \leq&\  \mu(1) e^{t\|g\|_{\infty}} \int dx \int_{0}^t  \mu(p_{s+h}^x+p_{s+h^{\prime}}^x)  ds \int_{s}^t r^{-d/2} dr \nn\\
			&+ e^{t\|g\|_{\infty}} \int dx \int_{0}^t ds \int_0^s \mu(p_{t-s+s_1+h}^x + p_{t-s+s_1+h^{\prime}}^x) d s_1 \int_{s_1}^{s}   r^{-d/2} dr \nonumber\\
			=& 2\mu(1)^2 e^{t\|g\|_{\infty}} \int_{0}^t ds \int_s^t r^{-d/2} dr + 2\mu(1) e^{t\|g\|_{\infty}} \int_{0}^t ds \int_0^s d s_1 \int_{s_1}^{s}   r^{-d/2} dr<\infty, 
		\end{align}
		as required.
		\end{proof}
		
		\begin{proposition}\label{l4.2}
		Let $d\leq3$. Suppose that \eqref{de.1.5} holds for $d=2,3$, then for any fixed \hk{$t\ge 0$} we have
		\begin{align}\label{e4.2}
			\lim _{h, h^{\prime} \downarrow 0} \int \E_\mu\big[ \big(Y_h(t, x)-Y_{h^{\prime}}(t, x)\big)^2 \big] dx=0 .
		\end{align}
	\end{proposition}
	\begin{proof}
		It suffices to show that for any fixed \hk{$t\ge 0$},
		\begin{align}\label{e4.7}
			\lim_{h, h^{\prime} \downarrow 0}\int \E_\mu\big[Y_h(t, x) Y_{h^{\prime}}(t, x)\big] dx = N_t <\infty
		\end{align} 
		and $N_t$ is independent of the ways of $h$ and $h^{\prime}$ \hk{approaching zero}. We first note that the finiteness of $N_t$ follows from Proposition \ref{l4.1}. Then it remains to prove the \hk{existence and uniqueness of the} limit in \eqref{e4.7}. We claim that for any $x\in\mathbb{R}^d$, 
		\begin{align}\label{e.4.15}
				\lim_{h, h^{\prime} \downarrow 0} V_t^{p_h^x, p_{h^{\prime}}^x}(u, v) = V_t^{x}(u, v)
		\end{align}
		holds for $\mu\times\mu$-almost every  $(u,v)\in\mathbb{R}^{2d}$, and for any $x\in\mathbb{R}^d$ and $s\in[0,t]$,
		\begin{align}\label{bue.4.15}
			\lim_{h, h^{\prime} \downarrow 0} V_s^{p_h^x, p_{h^{\prime}}^x}(v, v) = V_s^{x}(v, v)	
		\end{align}
		holds for $ds\mu(du)p_{t-s}(u,v)dv$-almost every $v\in\mathbb{R}^d$, where
		\begin{align*}
			V_t^{x}(u, v)=&\int_0^t p_{s}^x(u) ds \int_{s}^t p_{r}^x(v) dr+ \int_0^t p_{s}^x(v) ds \int_{s}^t p_{r}^x(u) dr \nn\\
			&+\int_0^t ds 
			\int_0^s \Pi_{(u, v)}\Big\{ p_{s-r}^x(B_r) Q_{t-s} p_{s-r}^x(\tilde{B}_r)\cdot  g(B_r, \tilde{B}_r) e^{\int_0^r g(B_a, \tilde{B}_a) da}  \Big\} dr\nonumber\\
			&+\int_0^t ds 
			\int_0^s \Pi_{(u, v)}\Big\{ p_{s-r}^x(\tilde{B}_r) Q_{t-s} p_{s-r}^x(B_r)\cdot  g(B_r, \tilde{B}_r) e^{\int_0^r g(B_a, \tilde{B}_a) da}  \Big\} dr.
		\end{align*}
It follows from \eqref{e.4.14}, \eqref{e.4.15} and \eqref{bue.4.15} that
		\begin{align}\label{e4.28}
			\lim_{h, h^{\prime} \downarrow 0} \E_\mu[Y_h(t, x) Y_{h^{\prime}}(t, x)] = &\int\!\!\int V_t^x(u,v) \mu(du) \mu(dv) \nonumber\\
			& \hk{+\int_0^t ds \int \mu(du) \int p_{t-s}(u,v) V_s^{x}(v,v) dv.}
		\end{align}
		We defer the proofs of \eqref{e.4.15} and \eqref{e4.28} to Appendix \ref{appc} since they are essentially applications of the dominated convergence theorem. Moreover, \eqref{bue.4.15} is proved using a similar argument to that employed in the proof of \eqref{e.4.15}.	Note that \hk{\eqref{e4.8} holds for $d\leq 3$. Thus we have
		\begin{align*}
			\mathbb{E}_{\mu}[Y_h(t, x) Y_{h^{\prime}}(t, x)] \leq &\ 2e^{t\|g\|_{\infty}}\Big(\int\mu(du)\int_0^{t+1} p_s^x(u)ds \Big)^2 \\
			&+ 2e^{t\|g\|_{\infty}}\int_0^t ds \int \mu(du) \int p_{t-s}(u,v)\Big(\int_0^{s+1} p_{s_1}^x(v)ds_1 \Big)^2 dv.
		\end{align*}
		By applying Cauchy-Schwarz's inequality to the first term, we see
		\begin{align*}
			\Big(\int\mu(du)\int_0^{t+1} p_s^x(u)ds \Big)^2 \leq \mu(1) \int\Big(\int_0^{t+1} p_s^x(u)ds \Big)^2\mu(du).
		\end{align*}
		A simple calculation shows that}
		\begin{align*}
			\bigg( \int_0^t p_s^x(u) ds \bigg)^2 &=\int_{0}^{t} ds \int_{0}^{t} (2\pi s)^{-d/2} (2\pi r)^{-d/2} e^{-\frac{(x-u)^2}{2s}-\frac{(x-u)^2}{2r}} dr \nonumber\\
			&= \int_{0}^{t} ds \int_{0}^{t} (2\pi (r+s))^{-d/2} p_{rs/(r+s)}^x(u) dr.
		\end{align*}	
		\hk{Combining the above yields} that for any $0<h,h^{\prime}<1$ and $x\in\mathbb{R}^d$,
		\begin{align*}
			\E_\mu[Y_h(t, x) &Y_{h^{\prime}}(t, x)] \leq 2\mu(1) e^{t \|g\|_{\infty} } \int_0^{t+1} d s \int \mu(du) \int_0^{t+1}(2 \pi(r+s))^{-d / 2} p_{r s /(r+s)}^x(u) d r \nonumber\\
			&+2 e^{t\|g\|_{\infty}} \int_0^t d s \int \mu(d u) \int_0^{s+1} d s_1 \int_0^{s+1} \big( 2 \pi(r+s_1) \big)^{-d / 2} p_{t-s+\frac{rs_1}{r+s_1}}^{x}(u) d r .
		\end{align*}
		\hk{It is easy to check that}
		\begin{align*}
			 \int dx \int_0^{t+1} d s \int \mu(du) &\int_0^{t+1}(2 \pi(r+s))^{-d / 2} p_{r s /(r+s)}^x(u) d r \nonumber\\
			 &= \mu(1) \int_0^{t+1} d s  \int_0^{t+1}(2 \pi(r+s))^{-d / 2} dr <\infty
		\end{align*}
		and 
		\begin{align*}
			\int dx \int_0^t d s \int \mu(d u) \int_0^{s+1} d s_1 &\int_0^{s+1} \big( 2 \pi(r+s_1) \big)^{-d / 2} p_{t-s+\frac{rs_1}{r+s_1}}^{x}(u) d r \nonumber\\
			&=\mu(1) \int_0^t d s  \int_0^{s+1} d s_1 \int_0^{s+1} \big( 2 \pi(r+s_1) \big)^{-d / 2} dr<\infty.
		\end{align*}
		Applying dominated convergence with the above arguments yields that
		\begin{align*}
			\lim_{h, h^{\prime} \downarrow 0} \int \E_\mu[Y_h(t, x) Y_{h^{\prime}}(t, x)] dx
			=  &\hk{\int \lim_{h, h^{\prime} \downarrow 0}\E_\mu[Y_h(t, x) Y_{h^{\prime}}(t, x)] dx} \nonumber \\
			=  &\int dx \int \mu(du) \int V_t^x(u,v) \mu(dv) \nonumber\\
			&\hk{+ \int dx \int_0^t ds \int \mu(du) \int p_{t-s}(u,v) V_s^{x}(v,v) dv=N_t}.
		\end{align*}
		The proof is now complete.
	\end{proof}

	Based on the proofs of Propositions \ref{l4.1} and \ref{l4.2}, it is easy to get the following result.
	\begin{corollary}\label{bc4.2}
		If $d\leq3$ and \eqref{de.1.5} holds for $d=2,3$, then for any $t\geq0$ and $x\in\mathbb{R}^d$ we have
		\begin{align}\label{e4.23}
			\lim _{h, h^{\prime} \downarrow 0}  \E_\mu\big[ \big(Y_h(t, x)-Y_{h^{\prime}}(t, x)\big)^2 \big] =0 .
		\end{align}
	\end{corollary}
	
	\begin{proof}
		This is similar to the proof of \eqref{e4.2}, so we only need to show that for any $t\geq 0$ and $x\in\mathbb{R}^d$, 
		\begin{align*}
			\lim _{h, h^{\prime} \downarrow 0}  \E_\mu\big[ Y_h(t, x)Y_{h^{\prime}}(t, x) \big] = N_{t,x}<\infty.
		\end{align*} 
		The existence and uniqueness of the above limit follow from \eqref{e4.28}. To prove the finiteness of $N_{t,x}$, it is sufficient to show that 
		\begin{align}\label{e4.29}
			\limsup _{h, h^{\prime} \downarrow 0} \E_\mu[Y_h(t, x) Y_{h^{\prime}}(t, x)] <\infty.
		\end{align}  
		Recall from \eqref{e.4.14} to see that
		\begin{align*}
			\E_\mu[Y_h(t, x) Y_{h^{\prime}}(t, x)]=\int\!\!\int V_t^{p_h^x, p_{h^{\prime}}^x} (u, v) \mu(d u) \mu(d v)+\int_0^t d s \int \mu(d u) \int p_{t-s}(u, v) V_{s}^{p_h^x,p_{h^{\prime}}^x} (v, v) d v.
		\end{align*}
		Hence, by using \eqref{e4.25}-\eqref{e4.27} it is easy to obtain that \eqref{e4.29} holds for $d=1$. By \eqref{e4.8}, \eqref{e4.2.2}, \eqref{e4.2.3} and \eqref{e4.11}, we see \eqref{e4.29} holds for $d=2$.  We use \eqref{e4.10}, \eqref{e4.3.2}-\eqref{e4.3.4} to prove \eqref{e4.29} holds for $d=3$. Thus, we complete the proof.
	\end{proof}
	
	By a combination of Propositions \ref{l4.1}-\ref{l4.2} and Corollary \ref{bc4.2}, we get the absolute continuity of $Y_t(dx)$ as follows. 
	
	\begin{corollary}\label{c4.2}
		Suppose \eqref{de.1.5} holds for $d=2,3$. Then with $\P_{\mu}$-probability one, $Y_t(dx)$ is absolutely  continuous with respect to Lebesgue measure $dx$ for all $t\geq0$.
	\end{corollary}
		
	\begin{proof}
		By Propositions \ref{l4.1}-\ref{l4.2} and Corollary \ref{bc4.2}, it is easily seen that there exists a measurable function $Y(t,x, \omega):[0, \infty) \times \mathbb{R}^d \times \Omega \rightarrow[0, \infty)$ such that  
		\begin{align}
			 & \E_\mu \big[Y(t,x)^2\big] <\infty, \quad \lim _{h \downarrow 0}  \E_\mu\Big[\big(Y_h(t,x)-Y(t,x)\big)^2\Big]  =0,\quad t\geq0,\ x\in\mathbb{R}^d,\label{e4.4.0}\\
			 &\int \E_\mu\big[ Y(t,x)^2 \big] d x <\infty, 		\quad \lim _{h \downarrow 0}  \int \E_\mu\Big[\big(Y_h(t,x)-Y(t,x)\big)^2\Big] d x =0,\quad t>0\label{e4.4.1},
		\end{align}
		where $Y_h(t,x)$ is defined by (\ref{e2.13}). For any $\phi\in C_c(\mathbb{R}^d)$, note that $|\phi(x)p_h^x(y)|\leq \|\phi\|_{\infty}p_h^x(y),\ p_h^x(y)\rightarrow \delta_x(y)$ as $h\downarrow0$ and
		\begin{align*}
			\lim_{h\downarrow0}\int dx \int p_h^x(y) Y_t(dy) = Y_t(1)= \int dx \int \lim_{h\downarrow0} p_h^x(y) Y_t(dy)  <\infty,\quad \P_{\mu}\text{-a.s.},
		\end{align*}
		where the finiteness comes from \eqref{e3}. Then by the generalized dominated convergence, we get
		\begin{align*}
			Y_t(\phi)  = \lim_{h\downarrow0} \int \phi(x) dx \int p_h^x(y) Y_t(dy) = \lim_{h\downarrow0} \int Y_h(t,x) \phi(x) d x,\quad \P_{\mu}\text{-a.s.}.
		\end{align*}
		Therefore, by Fatou's lemma and \eqref{e4.4.1} we have
		\begin{align*}
			\E_\mu\bigg[\Big|Y_t(\phi)-\int Y(t,x) \phi(x) d x\Big|^2\bigg]&\leq\liminf _{h \downarrow 0} \E_\mu\bigg[\Big|\int Y_h(t,x) \phi(x) d x-\int Y(t,x) \phi(x) d x\Big|^2\bigg] \nonumber\\
			&\leq \lim _{h \downarrow 0} \int \E_\mu\bigg[\Big( \hk{\text{$Y_h(t,x)  - Y(t,x)$}}\Big)^2\bigg] d x \int \phi(x)^2 dx, \nonumber\\
			&=0,
		\end{align*}
		where the second inequality follows from Cauchy-Schwarz's inequality. Summing up, we obtain that with $\P_{\mu}$-probability one, $Y_t(dx)$ is absolutely continuous with $dx$, the density of which coincides with $Y(t,x)$ for all $t\geq0$.
	\end{proof}
	
	The assertion (i) in Theorem \ref{t1} now follows by Corollary \ref{c4.2}.
		
	\section{Tanaka formula}\label{sec5}
	
	Let $d\leq 3$. Throughout this section, we assume $\alpha\geq0$ and require $\alpha>0$ if $d\leq2$. Recall that
	\begin{align*}
		g_\alpha^x(y)=\int_0^\infty e^{-\alpha t} p_t^x(y) dt,\quad x,y\in\mathbb{R}^d,
	\end{align*}
	where $p_t^x(y)=p_t(x,y)=\frac{1}{(2\pi t)^{d/2}} e^{-|y-x|^2/(2t)}$ is the transition density of the $d$-dimensional Brownian motion. It is easily seen that if $d=3$,
	\begin{align*} 
		g_0^x(y)= \frac{1}{2\pi |y-x|}, \quad y\neq x.
	\end{align*}
	For $d\leq 2$, we set
	\begin{align*} 
		g_0^x(y):= 
		\begin{cases}
			1+\log^+(1/|y-x|), & \text{ if } d=2;\\
			1,&\text{ if } d=1.
		\end{cases}
	\end{align*}
	Note that there exists a constant $C_\alpha>0$ such that
	\begin{align}\label{e.7.0}
		g_\alpha^x(y) \leq C_\alpha  g_0^x(y).
	\end{align}
	Indeed, this inequality follows by simple calculations when $d=1,3$ and by \cite[Appendix C(i)]{H18} when $d=2$. Furthermore, 		
		\begin{align}\label{e7.24}
			\int g_0^x(y) p_t^z(y) dy \leq C g_0^x(z),
		\end{align}
	for some constant $C>0$. In fact, \eqref{e7.24} obviously holds for $d=1$. We use \cite[Lemma 7.3]{H18} and \cite[Lemma 3.1]{H18} to prove \eqref{e7.24} holds for $d=2$ and $d=3$, respectively. 
	
	We now consider the martingale problem for the special test function $g_{\alpha, \eps}^x$ (on $\mathbb{R}^d$) defined by 
	\begin{align*} 
		g_{\alpha, \eps}^x(y):=\int_0^\infty e^{-\alpha t} p_{t+\eps}^x(y) dt,\quad 0< \eps< 1\ \text{and}\ x\in\mathbb{R}^d.
	\end{align*}
	In view of \eqref{e.7.0}, we have
	\begin{align*} 
		g_{\alpha, \eps}^x(y)\leq e^{\alpha}g_{\alpha}^x(y) \leq C_\alpha g_0^x(y).
	\end{align*}
	Thus,		
	\begin{align}\label{bu.e.100}
			\big|g_{\alpha,\eps}^x(y)- g_{\alpha}^x(y) \big| \leq g_{\alpha,\eps}^x(y)+ g_{\alpha}^x(y) \leq C_\alpha g_0^x(y).	
	\end{align}
	It is easy to check that as $\eps\downarrow 0$, 
	\begin{align}\label{bu.e1.0}
		\eta(\eps):=\sup_{x,y\in\mathbb{R}}	|g_{\alpha, \eps}^x(y)-g_{\alpha}^x(y)|   \to 0,\quad \text{when $d=1$,}
	\end{align}
	and for any $\delta>0$,
	\begin{align}\label{e1.0}
		\eta(\eps, \delta):=\sup_{|y-x|>\delta} |g_{\alpha, \eps}^x(y)-g_{\alpha}^x(y)|   \to 0,\quad\text{when $d=2,3$}.
	\end{align}
	Note that $g_{\alpha, \eps}^x \in C_b^2(\R^d)$ and
	\begin{align*} 
		\frac{\Delta}{2} g_{\alpha, \eps}^x(y)= \alpha g_{\alpha, \eps}^x(y)-p_\eps^x(y).
	\end{align*}
	By applying the martingale problem \eqref{eMP} with $\phi=g_{\alpha, \eps}^x$ we see that
	\begin{align} \label{e1}
		\int_0^t X_s(p_\eps^x) ds=X_0(g_{\alpha, \eps}^x)-X_t(g_{\alpha, \eps}^x)+\alpha \int_0^t X_s(g_{\alpha, \eps}^x) ds+M_t(g_{\alpha, \eps}^x),
	\end{align}
	where $(M_t(g_{\alpha, \eps}^x))_{t\ge 0}$ is a continuous $(\cF_t)$-martingale with 
	\begin{align*} 
		\langle M(g_{\alpha, \eps}^x)\rangle_t=\int_0^t X_s((g_{\alpha, \eps}^x)^2)ds+\int_0^t ds \int \int g(u,v) g_{\alpha, \eps}^x(u)g_{\alpha, \eps}^x(v) X_s(du) X_s(dv).
	\end{align*}
	To prove Tanaka formula (Theorem \ref{prop7.1}), it remains to prove that the right-hand side of \eqref{e1} converges in $L^1$.  We shall give these in the following three lemmas.

	\begin{lemma}\label{lem7.2}
		Let $d\leq 3$. Suppose that \eqref{de.1.5} holds for $d=2,3$. For any $t\geq0$ and $x\in\mathbb{R}^d$, then we have
		\begin{align}\label{e.7.2}
			\mathbb{E}_{\mu}\big[ | X_0(g_{\alpha,\eps}^x)-X_t(g_{\alpha,\eps}^x) - X_0(g_{\alpha}^x)+X_t(g_{\alpha}^x) | \big]\rightarrow 0,\quad\text{ as }\ \eps\downarrow0.
		\end{align}
	\end{lemma}

	\begin{proof}
		When $t\geq0$, for any fixed $x\in\mathbb{R}^d$ and $0<\eps<1$, a simple calculation shows 
		\begin{align}\label{e.7.8}
			&\ \mathbb{E}_{\mu}\big[ | X_0(g_{\alpha,\eps}^x)-X_t(g_{\alpha,\eps}^x) - X_0(g_{\alpha}^x)+X_t(g_{\alpha}^x) | \big] \nn\\
			\leq &\ \mathbb{E}_{\mu}\big[ | X_0(g_{\alpha,\eps}^x)- X_0(g_{\alpha}^x) | \big] + \mathbb{E}_{\mu}\big[ | X_t(g_{\alpha,\eps}^x)  -X_t(g_{\alpha}^x) | \big] \nn\\
			\leq &\ \langle \mu, |g_{\alpha,\eps}^x- g_{\alpha}^x | \rangle + \mathbb{E}_{\mu}\big[  X_t(|g_{\alpha,\eps}^x -g_{\alpha}^x|)  \big].
		\end{align} 
		When $d=1$, note that $\mathbb{E}_{\mu}[X_t(1)]=\mu(1)<\infty$ by Corollary \ref{cor3.2}, thus \eqref{e.7.2} follows by \eqref{bu.e1.0} and \eqref{e.7.8}.
		It remains to prove \eqref{e.7.2} holds when $d=2,3$. Indeed, it follows from \eqref{bu.e.100} and \eqref{e1.0} that for any $\delta>0$,
		\begin{align}\label{e.7.3}
			\langle \mu, |g_{\alpha,\eps}^x- g_{\alpha}^x | \rangle 
			&= \int_{|x-y|\leq \delta} \big|g_{\alpha,\eps}^x(y)- g_{\alpha}^x(y) \big| \mu(dy) + \int_{|x-y|> \delta} \big|g_{\alpha,\eps}^x(y)- g_{\alpha}^x(y) \big| \mu(dy)\nn\\
			&\leq C_{\alpha}\int_{|x-y|\leq \delta} g_0^x(y) \mu(dy) + \mu(1) \eta(\eps,\delta).
		\end{align}
		Suppose that	 the assumption \eqref{de.1.5} holds for $d=2,3$, then $\mu$ on one-point sets is $0$ and $\mu(g_0^x)<\infty$. Hence, by a truncation argument to $g_0^x$ and monotone convergence, we get 
		\begin{align}\label{xiu.e.1}
			\int_{|x-y|\leq \delta} g_0^x(y) \mu(dy) \rightarrow 0,\quad \text{ as }\ \delta\downarrow0.
		\end{align}
		By combining this with \eqref{e1.0}, we let $\eps\downarrow0$ first and then $\delta\downarrow0$ in both sides of \eqref{e.7.3} to obtain 
		\begin{align}\label{e.7.9}
			\lim_{\eps\downarrow0}\langle \mu, |g_{\alpha,\eps}^x- g_{\alpha}^x | \rangle=0.
		\end{align}
		Turning to the second term on the right-hand side of \eqref{e.7.8}, by arguments similar to those used in \eqref{e.7.3} we get that
		for any $\delta>0$,		
		\begin{align}\label{xiu.e.3}
			\mathbb{E}_{\mu} \big[  X_t(|g_{\alpha,\eps}^x  -g_{\alpha}^x|) \big]  
			&\leq C_{\alpha} \mathbb{E}_{\mu} \Big[ \int_{|x-y|\leq \delta} g_{0}^x(y) X_t(dy)  \Big] + \eta(\eps,\delta) \mathbb{E}_{\mu} [X_t(1)] \nn\\
			&= C_{\alpha} \int \mu(dz) \int_{|x-y|\leq\delta} g_0^x(y) p_t^z(y) dy + \mu(1)\eta(\eps,\delta),
		\end{align}
		where the last equality comes from Corollary \ref{cor3.2}. In view of \eqref{e7.24}, we see that
		\begin{align*}
			\int \mu(dz) \int g_0^x(y) p_t^z(y) dy\leq C \int g_0^x(z) \mu(dz)<\infty,
		\end{align*}
		as \eqref{de.1.5} holds for $d=2,3$. Also note that $\int \mu(dz) \int \mathbf{1}_{\{|x-y|\leq \delta\}}(y) dy\to 0$ as $\delta \downarrow 0$, so
		\begin{align}\label{xiu.e.4}
			\lim_{\delta\downarrow0}\int \mu(dz) \int_{|x-y|\leq\delta} g_0^x(y) p_t^z(y) dy=0
		\end{align}
		follows similarly by \eqref{xiu.e.1}. This combined with \eqref{e1.0} gives that the right-hand side of \eqref{xiu.e.3} tends to $0$ as $\eps\downarrow0$ first and then $\delta\downarrow0$, which means 
		\begin{align}\label{e.7.10}
			\mathbb{E}_{\mu} \big[ | X_t(g_{\alpha,\eps}^x)  -X_t(g_{\alpha}^x) | \big]\rightarrow0,\quad\text{as $\eps\downarrow0$.}
		\end{align}
		Summarizing, the desired result follows from \eqref{e.7.8}, \eqref{e.7.9} and \eqref{e.7.10}.
	\end{proof}

	\begin{lemma}\label{lem7.3}
		Suppose that $d\leq 3$ and \eqref{de.1.5} holds for $d=2,3$. For any $t\geq0,x\in\mathbb{R}^d$ and $\alpha\geq0\ (\alpha>0 \text{ if } d\leq 2)$, then we have
		\begin{align*}
			\mathbb{E}_{\mu} \bigg[ \Big| \int_0^t X_s(g_{\alpha, \eps}^x) ds - \int_0^t X_s(g_{\alpha}^x) ds \Big| \bigg]\rightarrow 0,\quad\text{ as }\ \eps\downarrow0.
		\end{align*}
	\end{lemma}

	\begin{proof}
		When $t\geq0$, for any fixed $x\in\mathbb{R}^d$ and $0<\eps<1$, we have
		\begin{align*}
			\mathbb{E}_{\mu} \bigg[ \Big| \int_0^t X_s(g_{\alpha, \eps}^x) ds-  \int_0^t X_s(g_{\alpha}^x) ds \Big| \bigg] = \mathbb{E}_{\mu} \big[ | Y_t(g_{\alpha, \eps}^x - g_{\alpha}^x) | \big] \leq \mathbb{E}_{\mu} \big[  Y_t(|g_{\alpha, \eps}^x - g_{\alpha}^x|)  \big].
		\end{align*}
		It then suffices to prove 
		\begin{align}\label{xiu.e.2}
			\mathbb{E}_{\mu} \big[  Y_t(|g_{\alpha, \eps}^x - g_{\alpha}^x|)  \big]\to0,\quad \text{as }\eps\downarrow0.
		\end{align}
		Since $\mathbb{E}_{\mu}[Y_t(1)] = t \mu(1)$ by \eqref{e3}, we use \eqref{bu.e1.0} to get \eqref{xiu.e.2} holds if $d=1$. Turning to $d=2,3$, by replacing $X_t$ with $Y_t$ in \eqref{xiu.e.3} and again using the moment for $Y_t$ to obtain
		\begin{align*}
			\mathbb{E}_{\mu} \big[  Y_t(|g_{\alpha,\eps}^x  -g_{\alpha}^x|) \big]  
			\leq C_{\alpha} \int \mu(dz) \int_0^t ds\int_{|x-y|\leq\delta} g_0^x(y) p_s^z(y) dy + t\mu(1)\eta(\eps,\delta).
		\end{align*}
		Suppose that \eqref{de.1.5} holds for $d=2,3$. Note that
		\begin{align*}
			\int \mu(dz) \int_0^t ds \int \mathbf{1}_{\{|x-y|\leq\delta\}}(y) dy \rightarrow0,\quad\text{as $\delta\downarrow0$},
		\end{align*}
		and by \eqref{e7.24},
		\begin{align*}
			\int \mu(dz) \int_0^t ds \int g_0^x(y) p_s^z(y) dy \leq C t \int g_0^x(z) \mu(dz)  <\infty.
		\end{align*}
		So by repeating arguments for deriving \eqref{e.7.10}, we have \eqref{xiu.e.2} holds if $d=2,3$. The proof is now finished.
	\end{proof}

	Let $X$ be a superprocess satisfying the martingale problem \eqref{eMP} on $(\Omega,\mathcal{F},\mathcal{F}_t,\mathbb{P})$ and $M$ be the martingale measure associated with $X$. Recall that for any $\psi\in\mathcal{L}^2(M)$, 	\begin{align*}
		M_t(\psi) = \int_0^t \int \psi(s,x) dM(s,x)
	\end{align*}
	is a continuous $(\mathcal{F}_t)$-martingale and 
	\begin{align*}
		\langle M(\psi) \rangle_t = \int_0^t X_s(\psi(s,\cdot)^2) ds + \int_0^t ds \int \int g(u,v)\psi(s,u) \psi(s,v) X_s(du) X_s(dv).
	\end{align*}
	
	\begin{lemma}\label{lem7.4}
		Suppose that the assumptions in Lemma \ref{lem7.3} hold. For any $t\geq0$ and $x\in\mathbb{R}^d$, let $M_t(g_{\alpha}^x)=\int_0^t\int g_{\alpha}^x(y) dM(r, y)$. Then
		\begin{itemize}
			\item [{\rm(i)}] $(M_t(g_{\alpha}^x), t\ge 0)$ is a continuous $L^2$-bounded $(\cF_t)$-martingale and		
			\begin{align*}
				\langle M(g_{\alpha}^x)\rangle_t=\int_0^t X_s((g_{\alpha}^x)^2)ds+\int_0^t ds \int\int g(u,v) g_{\alpha}^x(u)g_{\alpha}^x(v) X_s(du) X_s(dv);
			\end{align*}
			
			\item [{\rm(ii)}] $\lim_{\eps\downarrow0} \mathbb{E}_{\mu} [ | M_t(g_{\alpha, \eps}^x) - M_t(g_{\alpha}^x) | ] =0$.
		\end{itemize}
	\end{lemma}
	
	\begin{proof}
		(i) Suppose that \eqref{de.1.5} holds for $d=2,3$, it suffices to prove $g_{\alpha}^x\in\mathcal{L}^2(M)$ for any $x\in\mathbb{R}^d$.
	To do this, note that for any $t\geq 0$, 
	\begin{align}\label{e.7.21}
		&\ \mathbb{E}_{\mu} \Big[ \int_0^t  X_s((g_{\alpha}^x)^2)  ds +\int_0^t ds  \int\!\!\int g(u,v) g_{\alpha}^x(u)g_{\alpha}^x(v) X_s(du) X_s(dv) \Big] \nn\\
		\leq &\ \mathbb{E}_{\mu} [Y_t((g_{\alpha}^x)^2) ]+ \|g\|_{\infty} \int_0^t \mathbb{E}_{\mu} \big[(X_s(g_{\alpha}^x))^2 \big]  ds.
	\end{align}
	By \eqref{e.7.0} and \eqref{e3}, we have
	\begin{align}\label{ebc.1}
		\mathbb{E}_{\mu} [Y_t((g_{\alpha}^x)^2) ] &\leq C_\alpha \mathbb{E}_{\mu} [Y_t((g_{0}^x)^2) ] = C_\alpha \langle\mu,Q_t((g_{0}^x)^2)\rangle \nn\\
		&=C_\alpha \int\mu(dz) \int_0^t dr \int (g_{0}^x(y))^2 p_r^z(y) dy=C_\alpha I_3(t),
	\end{align}
	where 
	\begin{align}\label{xiu.e.8}
		I_3(t):=\int\mu(dz)\int_0^t dr \int (g_{0}^x(y))^2 p_r^z(y) dy.
	\end{align}
	We next prove $I_3(t)$ is finite.
	
	(a) When $d=1$, then $I_3(t)= t\mu(1)$.
	
	(b) When $d=2$, recall from \eqref{'bu5.63} to see that 
		\begin{align*}
			\Big( \log^+\frac{1}{|x-y|} \Big)^2 \leq \frac{1}{|x-y|}\quad \text{and} \quad \log^+\frac{1}{|x-y|} \leq \frac{1}{|x-y|}.
		\end{align*}
		Thus we have
		\begin{align*}
			I_3(t) 
			&=  \int \mu(dz) \int_0^t dr \int  \Big[ 1+\log^+\frac{1}{|x-y|} \Big]^2 p_r^z(y) dy \nn\\
			&\leq 3 \int \mu(dz) \int_0^t dr  \int  \frac{1}{|x-y|} p_r^z(y)dy + t \mu(1).
		\end{align*}
		It follows from \cite[Lemma 3.2]{H18} that
		\begin{align*}
			\int_0^t dr \int  \frac{1}{|x-y|} p_r^z(y)dy \leq 2\sqrt{2t}.
		\end{align*}
		Hence, 
		\begin{align*}
			I_3(t) \leq 6\sqrt{2t} \mu(1) + t \mu(1).
		\end{align*}
		
	(c) When $d=3$, we see
		\begin{align*}
			I_3(t) &\leq \int \mu(dz) \int_0^t d r  \int \frac{1}{|x-y|^2} p_r^z(y)dy 
			\leq 2 \int \Big[ \log^+\frac{1}{|x-z|} +1+\sqrt{3t} \Big] \mu(d z)\nn\\
			&\leq 2 \Big[\int \frac{1}{|x-z|} \mu(d z)+\mu(1)(1+\sqrt{3 t})\Big],
		\end{align*}
		where the second inequality comes from \cite[Lemma 5.1]{H18}. Then $I_3(t)<\infty$ if \eqref{de.1.5} holds for $d=3$. 
	This together with \eqref{ebc.1} implies $\mathbb{E}_{\mu} [Y_t((g_{\alpha}^x)^2) ]<\infty.$ 
	
	Turning to the second term on the right-hand side of \eqref{e.7.21}, recall from \eqref{e.7.0} and Corollary \ref{cor3.2} to see that
	\begin{align*}
		\mathbb{E}_{\mu} \big[(X_s(g_{\alpha}^x))^2 \big] &\leq C_\alpha \mathbb{E}_{\mu} \big[(X_s(g_{0}^x))^2 \big] \\
		&= C_\alpha \bigg\{\int\!\!\int \widetilde{V}_s^{g_0^x, g_0^x}(u, v) \mu(d u) \mu(d v) + \int_0^s d r \int \mu(d u) \int p_{s-r}(u, v) \widetilde{V}_r^{g_0^x, g_0^x}(v, v) d v \bigg\},\nn
	\end{align*}
	and 
	\begin{align*}
		\widetilde{V}_s^{g_0^x, g_0^x}(u, v) 
		&\leq e^{s \|g\|_{\infty}} \Pi_{(u, v)} \big[ g_0^x (B_s)  g_0^x(\tilde{B}_s) \big] = e^{s \|g\|_{\infty}} \Pi_{(u, v)} \big[ g_0^x (B_s) \big] \Pi_{(u, v)} \big[ g_0^x(\tilde{B}_s) \big]\nn\\
		&= e^{s \|g\|_{\infty}} \Big(\int g_0^x(y) p_{s}^u(y)dy\Big) \Big(\int g_0^x(y) p_{s}^v(y)dy\Big) \nn\\
		&\leq Ce^{s \|g\|_{\infty}}   g_{0}^x(u)  g_{0}^x(v),
	\end{align*}
	where we use \eqref{e7.24} in the last line. 
	It then follows that
	\begin{align*}
		\mathbb{E}_{\mu} \big[(X_s(g_{\alpha}^x))^2 \big] \leq C_\alpha  e^{s \|g\|_{\infty}} \bigg\{  \Big[\int g_0^x(u) \mu(du) \Big]^2 + \int_0^s dr \int \mu(du) \int p_{r}(u,v) \big(g_0^x(v)\big)^2 dv \bigg\}.
	\end{align*}
	Thus, 
	\begin{align}\label{e.7.28}
		\int_0^t \mathbb{E}_{\mu} \big[(X_s(g_{\alpha}^x))^2 \big] ds 
		&\leq C_\alpha  t e^{t \|g\|_{\infty}} \bigg\{  \Big[\int g_0^x(y) \mu(dy) \Big]^2+ \int_0^t dr \int \mu(du) \int p_{r}(u,v) \big(g_0^x(v)\big)^2 dv \bigg\} \nn\\
		&= C_\alpha  t e^{t \|g\|_{\infty}} \bigg\{  \Big[\int g_0^x(y) \mu(dy) \Big]^2 +I_3(t)\bigg\} <\infty,
	\end{align}
	if \eqref{de.1.5} holds for $d=2,3$. Summing up, we get $g_{\alpha}^x\in\mathcal{L}^2(M)$. 
	
	(ii) Since $g_{\alpha,\eps}^x\in C_b^2(\mathbb{R}^d)$, we get $g_{\alpha,\eps}^x - g_{\alpha}^x\in\mathcal{L}^2(M)$. Then $(M_t(g_{\alpha,\eps}^x - g_{\alpha}^x))_{t\geq 0}$ is a $(\mathcal{F}_t)$-martingale with quadratic variation 
	\begin{align*}
		\langle M(g_{\alpha,\eps}^x-g_{\alpha}^x)\rangle_t=&\int_0^t X_s((g_{\alpha,\eps}^x-g_{\alpha}^x)^2)ds\nn\\
		&+\int_0^t ds \int\!\!\int g(u,v) \big( g_{\alpha,\eps}^x(u)-g_{\alpha}^x(u) \big) \big( g_{\alpha,\eps}^x(v)-g_{\alpha}^x(v) \big) X_s(du) X_s(dv).
	\end{align*}
	It follows that 
	\begin{align}\label{e.7.26}
		&\ \mathbb{E}_{\mu} \Big[ \big| M_t(g_{\alpha, \eps}^x) - M_t(g_{\alpha}^x) \big| \Big]  = \mathbb{E}_{\mu} \Big[ \big| M_t(g_{\alpha, \eps}^x-g_{\alpha}^x) \big| \Big] \nn\\ 
		\leq &\ \Big\{ \mathbb{E}_{\mu} \big[ \big( M_t(g_{\alpha, \eps}^x-g_{\alpha}^x) \big)^2 \big] \Big\}^{1/2} = \Big\{ \mathbb{E}_{\mu} \big[ \langle M(g_{\alpha,\eps}^x-g_{\alpha}^x)\rangle_t \big] \Big\}^{1/2} \nn\\
		\leq &\ \bigg\{ \mathbb{E}_{\mu}\Big[ \int_0^t X_s \big( (g_{\alpha, \eps}^x-g_{\alpha}^x)^2 \big) ds \Big]  \bigg\}^{1/2} \nn\\
		&+ \bigg\{  \mathbb{E}_{\mu}\Big[ \int_0^t ds \int\!\!\int g(u,v) \big( g_{\alpha,\eps}^x(u)-g_{\alpha}^x(u) \big) \big( g_{\alpha,\eps}^x(v)-g_{\alpha}^x(v) \big) X_s(du) X_s(dv) \Big] \bigg\}^{1/2} \nn\\
		\leq &\ \bigg\{ \mathbb{E}_{\mu} \big[ Y_t \big( (g_{\alpha, \eps}^x-g_{\alpha}^x)^2 \big) \big] \bigg\}^{1/2} + \bigg\{ \|g\|_{\infty} \mathbb{E}_{\mu} \Big[  \int_0^t  \Big(\int |g_{\alpha,\eps}^x(u) - g_{\alpha}^x(u)| X_s(du)\Big)^2  ds \Big] \bigg\}^{1/2}. 
	\end{align}
	Note that if \eqref{de.1.5} holds for $d=2,3$,
	\begin{align*}
		\int\mu(dz) \int_0^t ds \int \big( g_{0}^x(y) \big)^2 p_s^z(y) dy=I_3(t)<\infty,
	\end{align*}
	by repeating arguments for proving \eqref{xiu.e.2} with $|g_{\alpha,\eps}^x-g_{\alpha}^x|$ replaced by $(g_{\alpha,\eps}^x-g_{\alpha}^x)^2$ we get
	\begin{align}\label{xiu.e.5}
		\mathbb{E}_{\mu} \big[ Y_t \big( (g_{\alpha, \eps}^x-g_{\alpha}^x)^2 \big) \big] \to 0,\quad\text{as $\eps\downarrow0$}.	
	\end{align}
	Turning to the second term on the right-hand side of \eqref{e.7.26}, when $d=1$, from \eqref{bu.e1.0} we have
	\begin{align}\label{xiu.e.6}
		\mathbb{E}_{\mu} \Big[  \int_0^t  \Big(\int |g_{\alpha,\eps}^x(u) - g_{\alpha}^x(u)| X_s(du)\Big)^2  ds \Big] 
		&\leq \eta(\varepsilon)^2\int_0^t \mathbb{E}_{\mu}[  X_s(1)^2 ]ds \nn\\
		&\leq C_t \eta(\varepsilon)^2 \rightarrow 0,\quad\text{as $\varepsilon\downarrow 0$},
	\end{align}
	where the last inequality comes from Lemma \ref{lemx5.1}. When $d=2,3$, by \eqref{bu.e.100} and \eqref{e1.0} we obtain that for any $\delta>0$,
	\begin{align}\label{e.7.27}
		 &\ \mathbb{E}_{\mu} \Big[  \int_0^t  \Big(\int |g_{\alpha,\eps}^x(u) - g_{\alpha}^x(u)| X_s(du)\Big)^2  ds \Big] \nn\\
		\leq&\ C_{\alpha} \mathbb{E}_{\mu} \Big[ \int_0^t \Big( \int_{|x-u|\leq\delta} g_{0}^x(u) X_s(du) \Big)^2  ds \Big] +  \eta(\eps,\delta)^2 \int_0^t
		\mathbb{E}_{\mu}[ X_s(1)^2 ]ds.
	\end{align}
	Furthermore, the second term on the right-hand side in the above tends to $0$ as $\varepsilon\downarrow0$. On the other hand,  
	let $\phi_{\delta}^x(u):=\mathbf{1}_{\{ |x-u|\leq\delta \}}(u) g_0^x(u)$. By Corollary \ref{cor3.2} we see
	\begin{align*}
		\mathbb{E}_{\mu} \Big[  \Big( \int_{|x-u|\leq\delta} g_{0}^x(u) X_s(du) \Big)^2  \Big] = \int\!\!\int& \widetilde{V}_s^{\phi_{\delta}^x, \phi_{\delta}^x}(u, v) \mu(d u) \mu(d v) \nn\\
		&+ \int_0^s d r \int \mu(d u) \int p_{s-r}(u, v) \widetilde{V}_r^{\phi_{\delta}^x, \phi_{\delta}^x}(v, v) d v,
	\end{align*}
	and
	\begin{align*}
		\widetilde{V}_s^{\phi_{\delta}^x, \phi_{\delta}^x}(u, v) 
		&\leq e^{s \|g\|_{\infty}} \Pi_{(u, v)} \big[\phi_{\delta}^x (B_s) \big] \Pi_{(u, v)} \big[ \phi_{\delta}^x(\tilde{B}_s) \big] \\
		&=e^{s \|g\|_{\infty}} \Big( \int_{|x-y|\leq\delta} g_0^x(y) p_s^u(y) dy\Big) \Big( \int_{|x-y|\leq\delta} g_0^x(y) p_s^v(y) dy\Big) := \widetilde{V}_s^{x,\delta}(u,v)\nn.
	\end{align*}
	Therefore, 
	\begin{align*}
		\mathbb{E}_{\mu} \Big[ \int_0^t \Big( \int_{|x-u|\leq\delta} g_{0}^x(u) X_s(du) \Big)^2 & ds \Big] \leq \int_0^t ds \int\!\!\int \widetilde{V}_s^{x, \delta}(u, v) \mu(d u) \mu(d v) \\
		&+ \int_0^t ds \int_0^s d r \int \mu(d u) \int p_{s-r}(u, v) \widetilde{V}_r^{x, \delta}(v, v) d v.\nn
	\end{align*}
	From \eqref{e7.24}, it is easy to check that
	\begin{align*}
		&\int_0^t  \Big( \int \mu(du) \int g_0^x(y) p_s^u(y) dy\Big)^2 ds \leq C t \big[\mu(g_0^x)\big]^2 <\infty,\nn\\
		&\int_0^t ds \int_0^s d r \int \mu(d u) \int p_{s-r}(u, v)  \Big( \int g_0^x(y) p_s^v(y) dy \Big)^2 dv \leq C \int_0^t I_3(s) ds<\infty,
	\end{align*}
	where $I_3(s)$ is given by \eqref{xiu.e.8}. By arguments similar to the proof of \eqref{xiu.e.4}, we get
	\begin{align*}
		 \lim_{\delta\downarrow0}\mathbb{E}_{\mu} \Big[ \int_0^t \Big( \int_{|x-u|\leq\delta} g_{0}^x(u) X_s(du) \Big)^2  &ds \Big] 
		\leq \lim_{\delta\downarrow0}\int_0^t ds \int\!\!\int \widetilde{V}_s^{x, \delta}(u, v) \mu(d u) \mu(d v) \nn\\
		&+ \lim_{\delta\downarrow0}\int_0^t ds \int_0^s d r \int \mu(d u) \int p_{s-r}(u, v) \widetilde{V}_r^{x, \delta}(v, v) d v=0.
	\end{align*}
	Hence, by letting $\eps\downarrow0$ first and then $\delta\downarrow0$ in both sides of \eqref{e.7.27} we obtain that if $d=2,3$,
	\begin{align*}
		\mathbb{E}_{\mu} \Big[  \int_0^t  \Big(\int |g_{\alpha,\eps}^x(u) - g_{\alpha}^x(u)| X_s(du)\Big)^2  ds \Big]\rightarrow0,\quad\text{as $\varepsilon\downarrow0$}.	
	\end{align*}
	This together with \eqref{e.7.26}, \eqref{xiu.e.5} and \eqref{xiu.e.6} gives the assertion (ii).
	\end{proof}

	

\section{Local joint H\"older continuity}\label{sec6}
We have proved the existence of the density process $Y(t,x)$ for $d\leq 3$ in Section \ref{sec3}. In this section, we shall investigate the local joint H\"older continuity of $Y(t,x)$. To this end, it is sufficient to give the estimates on the spatial and time differences.

Before deriving the estimates, we first give another characterization of the density $Y(t,x)$. Suppose that \eqref{de.1.5} holds for $d=2,3$. For any $t\geq0$ and $a\in\mathbb{R}^d$, recall from Corollary \ref{c4.2} that
\begin{align*}
	Y_{h}(t,a)=Y_t(p_h^a)=\int p_h^a(y) Y_t(dy) \xrightarrow[]{\mathbb{P}_{\mu}} Y(t,a),\quad \text{ as } h\downarrow 0,
\end{align*}
where the symbol $\xrightarrow[]{\mathbb{P}_{\mu}}$ denotes the convergence in probability. By Tanaka formula (see Theorem \ref{prop7.1}), we obtain that as $h\downarrow 0$,
\begin{align*}
	Y_t(p_h^a)=\int_0^t X_s(p_h^a) ds \xrightarrow[]{L^1} L_t^a. 
\end{align*}
Hence, we have $Y(t,a)=L_t^a,\ \mathbb{P}_{\mu}
$-a.s., i.e.,
\begin{align}\label{eb.7.5}
	Y(t,a) = X_0(g_{\alpha}^a)-X_t(g_{\alpha}^a)+\alpha \int_0^t X_s(g_{\alpha}^a) ds+M_t(g_{\alpha}^a),\quad \mathbb{P}_{\mu}\text{-a.s.},
\end{align}
where 
\begin{align*}		
g_\alpha^a(x):=\int_0^\infty e^{-\alpha t} p_t^a(x) dt,  \quad \alpha \geq 0\ (\alpha>0 \text{ if } d=1,2)
\end{align*}
is defined by \eqref{eb.6.1} and $(M_t(g_{\alpha}^a))_{t\ge 0}$ is a continuous $(\mathcal{F}_t)$-martingale given in Lemma \ref{lem7.4}. Throughout this section    we always assume $\alpha\geq 0\ (\alpha>0$ if $d=1,2)$.

\subsection{Estimates on the spatial difference}

In view of \eqref{eb.7.5}, to obtain the estimate on the spatial difference, we need the following H\"older continuity of $g_\alpha^a$. \begin{lemma}\label{lem.eb.7.1}
	Let $d\leq 3$ and $0<\gamma<1$. There is some constant $C_{\alpha, \gamma}>0$ such that the following holds for all $x,a,b \in \R^d$.
	\begin{itemize}
		\item [\rm{(i)}] When $d=1$, we have
	\begin{align*}
		|g_\alpha^a(x)-g_\alpha^{b}(x)| \leq C_{\alpha, \gamma} |a-b|^\gamma;
	\end{align*}
	\item[\rm{(ii)}] When $d=2, 3$, we have
	\begin{align*}
		|g_\alpha^a(x)-g_\alpha^{b}(x)| \leq 
		C_{\alpha, \gamma} |a-b|^\gamma \Big(\frac{1}{|a-x|^{\gamma+d-2}}+\frac{1}{|b-x|^{\gamma+d-2}}\Big).
	\end{align*}
	\end{itemize}
\end{lemma}

\begin{proof}
	Recall from (3.44) of Sugitani \cite{S89} that for any $0<\gamma<1$, there is a constant $C_\gamma>0$ such that 
	\begin{align*}
		|p_t^a(x)-p_t^b(x)| \leq C_\gamma t^{-\gamma/2} |a-b|^\gamma \big(p_{2t}^a(x)+p_{2t}^b(x)\big), \quad \forall\  t>0,\ x, a,b\in \R^d.
	\end{align*}
	Hence, it follows that
	\begin{align*}
		|g_\alpha^a(x)-g_\alpha^{b}(x)| &\leq  \int_0^\infty e^{-\alpha t} |p_t^a(x)-p_t^b(x)| dt\\
		&\leq C_\gamma  |a-b|^\gamma \int_0^\infty e^{-\alpha t} t^{-\gamma/2} \big(p_{2t}^a(x)+p_{2t}^b(x)\big) dt.
	\end{align*}

	(i) When $d=1$, we obtain
	\begin{align*}
		|g_\alpha^a(x)-g_\alpha^{b}(x)| &\leq   C_\gamma  |a-b|^\gamma \int_0^\infty e^{-\alpha t} t^{-(1+\gamma)/2}   dt\leq C_{\alpha,\gamma}  |a-b|^\gamma,
	\end{align*}
	for some constant $C_{\alpha,\gamma}>0$ and the last inequality follows by $(1+\gamma)/2<1$. 
	
	(ii) When $d=2$ or $3$, we obtain
	\begin{align*}
		|g_\alpha^a(x)-g_\alpha^{b}(x)| &\leq    C_{\gamma}  |a-b|^\gamma \int_0^\infty   t^{-(d+\gamma)/2} (e^{-\frac{|a-x|^2}{4t}}+e^{-\frac{|b-x|^2}{4t}}) dt\\
		&\leq C_{\gamma}  |a-b|^\gamma \Big(\frac{1}{|a-x|^{\gamma+d-2}}+\frac{1}{|b-x|^{\gamma+d-2}}\Big) \int_0^\infty    s^{(d+\gamma-4)/2} e^{-s} ds \\
		&\leq C_{\gamma}  |a-b|^\gamma \Big(\frac{1}{|a-x|^{\gamma+d-2}}+\frac{1}{|b-x|^{\gamma+d-2}}\Big),
	\end{align*}
	where the last inequality follows by $(d+\gamma-4)/2>-1$.
\end{proof}

Based on \eqref{eb.7.5} and Lemma \ref{lem.eb.7.1}, we derive the estimates on the spatial difference under the strengthened assumption \eqref{de1.3}.

\begin{theorem}\label{th.e.7.2}
	Let $d\leq 3$, $N\geq 1$ and $0<\eps_0<T$. 
	\begin{itemize}
		\item [\rm{(i)}] When $d=1$, for any $0<\gamma_1<1$, there exists a constant $C_{N,T,\alpha,\gamma_1}>0$ such that for any $0\leq s,t\leq T$ and $a,b\in\mathbb{R}$ with $|a|,|b|\leq T$,
			\begin{align*}
				\mathbb{E}_{\mu} \big[ | Y(t,a) -Y(t,b) |^{2N} \big] \leq C_{N,T,\alpha,\gamma_1} |a-b|^{2 N\gamma_1}.
			\end{align*}

		\item [\rm{(ii)}] When $d=2$ or $3$, we further assume that \eqref{de1.3} holds. Let $\gamma\in(0,1\wedge (2-\frac{d}{2}))$ be the H\"older exponent given in this assumption. Then there is a constant $C_{N,\eps_0,T,\alpha,\gamma}>0$ such that for all $\eps_0\leq s,t\leq T$ and $a,b\in\mathbb{R}^d$ with $|a|,|b|\leq T$,
	\begin{align}\label{eb.7.6}
		\mathbb{E}_{\mu} \big[ | Y(t,a) -Y(t,b) |^{2N} \big] \leq C_{N,\eps_0,T,\alpha,\gamma} |a-b|^{2 N\gamma}.
	\end{align}
	\end{itemize} 
\end{theorem}

\begin{proof}
	Fix $N\geq 1$ and $0<\varepsilon_0<T$. By \eqref{eb.7.5}, we see that for any $t\geq0$ and $a,b\in\mathbb{R}^d$,
	\begin{align*}
		\mathbb{E}_{\mu} \big[  | Y(t,a) - Y(t,b) |^{2N} \big] \leq 2^{4N} \Big\{ &\mathbb{E}_{\mu} \big[ | X_0(g_{\alpha}^a)- X_0(g_{\alpha}^b) |^{2N} \big] +  \mathbb{E}_{\mu} \big[ | X_t(g_{\alpha}^a)- X_t(g_{\alpha}^b)|^{2N} \big] \nonumber\\
		&+ \alpha^{2N} \mathbb{E}_{\mu} \bigg[ \Big| \int_{0}^{t} X_s(g_{\alpha}^a) ds - \int_{0}^{t} X_s(g_{\alpha}^b) ds \Big|^{2N} \bigg] \nonumber\\
		&+ \mathbb{E}_{\mu} \big[ |M_t(g_{\alpha}^a)-M_t(g_{\alpha}^b)|^{2N} \big] \Big\}, \nonumber\\
		\leq 2^{4N} \Big\{ & \big| \mu(g_{\alpha}^a-g_{\alpha}^b) \big|^{2N}  +  \mathbb{E}_{\mu} \big[ X_t\big( |g_{\alpha}^a-g_{\alpha}^b| \big)^{2N} \big] \nonumber\\
		&+ \alpha^{2N} \mathbb{E}_{\mu} \big[ Y_t\big( |g_{\alpha}^a-g_{\alpha}^b| \big)^{2N} \big] + \mathbb{E}_{\mu} \big[ M_t(g_{\alpha}^a-g_{\alpha}^b)^{2N} \big] \Big\}. 
	\end{align*}
We next evaluate the four terms on the right-hand side of the above inequality. In the remainder of this proof, let $\gamma$ be given in the assumption \eqref{de1.3}, and we suppose $\gamma_1\in(0,1)$.

\textbf{The first term:} 

(i) When $d=1$, we use Lemma \ref{lem.eb.7.1}(i) to get
	\begin{align*}
		\big| \mu(g_{\alpha}^a-g_{\alpha}^b) \big| \leq C_{\alpha,\gamma_1} \mu(1) |a-b|^{\gamma_1}.
	\end{align*}
	
	(ii) When $d=2$, let
	\begin{align*}
		f_{\alpha}(x):=\int_0^{\infty} e^{-\alpha s} p_s(0,x) ds - \frac{1}{\pi} \log^+\frac{1}{|x|}=\int_0^{\infty} e^{-\alpha s} p_s(0,x) ds - \frac{1}{\pi}g_2(x).
	\end{align*}
	Then we see
	\begin{align*}
		g_{\alpha}^a(x)=\frac{1}{\pi}g_2(a-x) + f_{\alpha}(a-x).
	\end{align*}
	Hence,
	\begin{align*}
		\mu(g_{\alpha}^a)= \frac{1}{\pi} \int g_2(a-x) \mu(dx) + \int f_{\alpha}(a-x) \mu(dx) = \frac{1}{\pi} \mu g_2(a) + \mu f_{\alpha}(a).
	\end{align*}
	The assumption \eqref{de1.3} implies $\mu g_2(a)$ is locally $\gamma$-H\"older continuous. Since 
	\begin{align}\label{hc1}
		\text{$\mu f_{\alpha}(a)$ is locally $\gamma$-H\"older continuous with respect to $a$}
	\end{align}
	 (see Appendix \ref{AA}), so is $\mu(g_{\alpha}^a)$. 
	
	(iii)\ When $d=3$, we define a different function 
	\begin{align*}
			\widetilde{f}_{\alpha}(x):=\int_0^{\infty} e^{-\alpha s} p_s(0,x) ds - \frac{1}{2\pi|x|} = \int_0^{\infty} e^{-\alpha s} p_s(0,x) ds - \frac{1}{2\pi}g_3(x).
	\end{align*}
	A simple calculation shows 
	\begin{align*}
		\mu(g_{\alpha}^a)=\frac{1}{2\pi}\mu g_3(a)+\mu\widetilde{f}_{\alpha}(a).	
	\end{align*}
By the assumption \eqref{de1.3}, we see $\mu g_3(a)$ is locally $\gamma$-H\"older continuous. Since
\begin{align}\label{hc2}
	\text{$\mu\widetilde{f}_{\alpha}(a)$ is locally $\gamma$-H\"older continuous with respect to $a$}
\end{align}
  (see Appendix \ref{AA}), so is $\mu(g_{\alpha}^a)$.

\textbf{The second term:} In view of Lemma \ref{lem.eb.7.1} we obtain the following estimates.  


(i) When $d=1$, since $\gamma_1\in(0,1)$, we have
\begin{align*}
	\mathbb{E}_{\mu} \big[ X_t\big( |g_{\alpha}^a-g_{\alpha}^b| \big)^{2N} \big]
	&\leq C_{N,\alpha,\gamma_1}  |a-b|^{2N\gamma_1} \mathbb{E}_{\mu} \big[ X_t(1)^{2N}\big] \nonumber\\
	&\leq C_{N,T,\alpha,\gamma_1} |a-b|^{2N\gamma_1},\quad \forall\ 0\leq t \leq T\text{ and }a,b\in\mathbb{R},
\end{align*}
where $C_{N,T,\alpha,\gamma_1}>0$ is a constant and the last inequality comes from Lemma \ref{lemx5.1}.

(ii) When $d=2,3$, recall that the text function 
\begin{align*}
		\phi_{a,\gamma}(x)=\frac{1}{|a-x|^{\gamma}},\quad x\in \mathbb{R}^d
\end{align*}
is defined by \eqref{xphig}. Since $\gamma \in(0, 1\wedge(2-\frac{d}{2}))\subset(0,1)$, we get
\begin{align*}
	\mathbb{E}_{\mu} \big[ X_t\big( |g_{\alpha}^a-g_{\alpha}^b| \big)^{2N} \big]
	&\leq C_{N,\alpha,\gamma} |a-b|^{2N\gamma } \mathbb{E}_{\mu}\bigg[ \Big(\int \Big(\frac{1}{|a-x|^{\gamma+d-2}}+\frac{1}{|b-x|^{\gamma+d-2}}\Big) X_t(dx) \Big)^{2N} \bigg] \nonumber\\
	&= C_{N,\alpha,\gamma} |a-b|^{2N\gamma } \mathbb{E}_{\mu}\Big[ \big(  X_t(\phi_{a,\gamma+d-2}) + X_t(\phi_{b,\gamma+d-2})  \big)^{2N} \Big] \nonumber\\
	&\leq C_{N,\alpha,\gamma} |a-b|^{2N\gamma } \Big\{ \mathbb{E}_{\mu}\big[ X_t(\phi_{a,\gamma+d-2})^{2N} \big] +  \mathbb{E}_{\mu}\big[ X_t(\phi_{b,\gamma+d-2})^{2N} \big] \Big\}.
\end{align*}
We assume $\eps_0\leq t\leq T$ and $a\in\mathbb{R}^d$. If $d=2$, since $0<\gamma<1\wedge(2-\frac{d}{2})=1$, by Corollary \ref{'corb5.8} we get
\begin{align*}
	\mathbb{E}_{\mu}\big[ X_t(\phi_{a,\gamma+d-2})^{2N} \big] =\mathbb{E}_{\mu}\big[ X_t(\phi_{a,\gamma})^{2N} \big] \leq C_{N,\eps_0,T,\gamma}.
\end{align*}
If $d=3$, since $0<\gamma<1\wedge(2-\frac{d}{2})=1/2$, by Corollary \ref{'corb5.11} we see 
\begin{align*}
	\mathbb{E}_{\mu}\big[ X_t(\phi_{a,\gamma+d-2})^{2N} \big] =\mathbb{E}_{\mu}\big[ X_t(\phi_{a,\gamma+1})^{2N} \big] \leq C_{N,\eps_0,T,\gamma}.
\end{align*} 
Therefore, when $d=2,3$, we have
\begin{align*}
	\mathbb{E}_{\mu} \big[ X_t\big( |g_{\alpha}^a-g_{\alpha}^b| \big)^{2N} \big] \leq C_{N,\eps_0,T,\alpha,\gamma} |a-b|^{2N\gamma }, \quad\forall\ \eps_0\leq t\leq T\text{ and } a,b\in\mathbb{R}^d.
\end{align*}

\textbf{The third term:} By repeating the arguments in the second term and replacing $X_t$ by $Y_t$, we get the following estimates.

(i) When $d=1$, by Lemma \ref{lem5.5} with $\phi\equiv1$ we obtain
\begin{align*}
	\mathbb{E}_{\mu} \big[ Y_t\big( |g_{\alpha}^a-g_{\alpha}^b| \big)^{2N} \big] &\leq C_{N,\alpha,\gamma_1} |a-b|^{2N\gamma_1} \mathbb{E}_{\mu} \big[ Y_t(1)^{2N} \big] \nonumber\\
	&\leq C_{N,T,\alpha,\gamma_1} |a-b|^{2N\gamma_1},\quad \forall\  0\leq t \leq T\text{ and }a,b\in\mathbb{R}, 
\end{align*}
for some $C_{N,T,\alpha,\gamma_1}>0$.

(ii) When $d=2,3$,
\begin{align*}
	\mathbb{E}_{\mu} \big[ Y_t\big( |g_{\alpha}^a-g_{\alpha}^b| \big)^{2N} \big]
	\leq C_{N,\alpha,\gamma} |a-b|^{2N\gamma } \Big\{ \mathbb{E}_{\mu}\big[ Y_t(\phi_{a,\gamma+d-2})^{2N} \big] +  \mathbb{E}_{\mu}\big[ Y_t(\phi_{b,\gamma+d-2})^{2N} \big] \Big\}.
\end{align*}
If $d=2$, then $0< \gamma < 1\wedge(2-\frac{d}{2})$ implies $\gamma+d-2=\gamma\in(0,1)$. Hence, there is a constant $C_{N,T,\alpha,\gamma}>0$ such that
\begin{align*}
	\mathbb{E}_{\mu} \big[ Y_t\big( |g_{\alpha}^a-g_{\alpha}^b| \big)^{2N} \big]
	&\leq C_{N,\alpha,\gamma} |a-b|^{2N\gamma } \Big\{ \mathbb{E}_{\mu}\big[ Y_t(\phi_{a,\gamma})^{2N} \big] +  \mathbb{E}_{\mu}\big[ Y_t(\phi_{b,\gamma})^{2N} \big] \Big\}\nonumber\\
	&\leq C_{N,T,\alpha,\gamma} |a-b|^{2N\gamma},\quad \forall\ 0\leq t\leq T\text{ and }a,b\in\mathbb{R}^2,
\end{align*}
where the last inequality follows from Corollary \ref{'cor5.15}. If $d=3$, since $0< \gamma < 1\wedge(2-\frac{d}{2})$ implies $\gamma+d-2=\gamma+1\in(1,3/2)$, by Corollary \ref{'cor5.18} we obtain  
\begin{align*}
	\mathbb{E}_{\mu} \big[ Y_t\big( |g_{\alpha}^a-g_{\alpha}^b| \big)^{2N} \big]
	&\leq C_{N,\alpha,\gamma} |a-b|^{2N\gamma } \Big\{ \mathbb{E}_{\mu}\big[ Y_t(\phi_{a,\gamma+1})^{2N} \big] +  \mathbb{E}_{\mu}\big[ Y_t(\phi_{b,\gamma+1})^{2N} \big] \Big\}\nonumber\\
	&\leq C_{N,T,\alpha,\gamma} |a-b|^{2N\gamma} 
	\end{align*}
holds for any $0\leq t\leq T$ and $a,b\in\mathbb{R}^3$ with $|a|,|b|\leq T$.

\textbf{The forth term:} Note that $(M_t(g_{\alpha}^a-g_{\alpha}^b), t\geq0)$ is a continuous martingale with quadratic variation
\begin{align*}
	\langle M(g_{\alpha}^a-g_{\alpha}^b)\rangle_t=\int_0^t &X_s\big((g_{\alpha}^a-g_{\alpha}^b)^2\big)ds\\
	&+\int_0^t ds \int \int g(x,y) \big( g_{\alpha}^a(x) - g_{\alpha}^b(x)\big) \big( g_{\alpha}^a(y) - g_{\alpha}^b(y)\big) X_s(dx) X_s(dy).\nonumber
\end{align*}
Applying Burkholder-Davis-Gundy's inequality yields
\begin{align}\label{eb.7.17}
	&\mathbb{E}_{\mu} \big[ M_t(g_{\alpha}^a-g_{\alpha}^b)^{2N} \big] \nonumber\\
	\leq &C_N \mathbb{E}_{\mu} \bigg[ \Big( \int_0^t X_s\big((g_{\alpha}^a-g_{\alpha}^b)^2\big)ds \nonumber\\  
	&\qquad\qquad + \int_0^t ds \int \int g(x,y) \big( g_{\alpha}^a(x) - g_{\alpha}^b(x)\big) \big( g_{\alpha}^a(y) - g_{\alpha}^b(y)\big) X_s(dx) X_s(dy) \Big)^N \bigg] \nonumber\\
	\leq\ &C_N \mathbb{E}_{\mu} \Big[ \Big( \int_0^t X_s\big((g_{\alpha}^a-g_{\alpha}^b)^2\big)ds \Big)^N \Big] + C_N \|g\|_{\infty}^N \mathbb{E}_{\mu} \Big\{ \Big[ \int_0^t \big(X_s(|g_{\alpha}^a-g_{\alpha}^b|) \big)^2 ds \Big]^N \Big\} \nonumber\\
	\leq\ & C_N \mathbb{E}_{\mu} \Big\{ \big[Y_t((g_{\alpha}^a-g_{\alpha}^b)^2)\big]^N \Big\} + C_N \mathbb{E}_{\mu} \Big\{ \Big[ \int_0^t \big(X_s(|g_{\alpha}^a-g_{\alpha}^b|) \big)^2 ds \Big]^N \Big\}.
\end{align}
Here and in what follows, the constant $C$ depends on $\|g\|_{\infty}$ and we omit this dependence. 
For the first term on the right-hand side of \eqref{eb.7.17}, we again use Lemma \ref{lem.eb.7.1} to get that
\begin{equation*}
\mathbb{E}_{\mu} \Big\{ \big[Y_t((g_{\alpha}^a-g_{\alpha}^b)^2)\big]^N \Big\}\leq 
	\begin{cases}
		 C_{N,\alpha,\gamma_1}  |a-b|^{2N\gamma_1}\mathbb{E}_{\mu} \big[Y_t(1)^{N}\big]\leq C_{N,T,\alpha,\gamma_1} |a-b|^{2N\gamma_1}, & \text{if } d=1; \\
		 C_{N,\alpha,\gamma} |a-b|^{2N\gamma} \mathbb{E}_{\mu} \Big\{ \big[ Y_t\big((\phi_{a,\gamma+d-2}+\phi_{b,\gamma+d-2})^2\big)\big]^N \Big\}, & \text{if } d=2,3,
	\end{cases}
\end{equation*}
holds for every $0\leq t \leq T$ and the inequality in $d=1$ follows by $\gamma_1\in(0,1)$ and Lemma \ref{lem5.5}. When $d=2,3$, by the definition of $\phi_{a,\gamma}$ we see
\begin{align*}
	\mathbb{E}_{\mu} \Big\{ \big[ Y_t\big((\phi_{a,\gamma+d-2}+\phi_{b,\gamma+d-2})^2\big)\big]^N \Big\} 
	&\leq 2^N \mathbb{E}_{\mu} \Big\{ \big[ Y_t(\phi_{a,2(\gamma+d-2)})+Y_t(\phi_{b,2(\gamma+d-2)})\big]^N \Big\} \nonumber\\
	&\leq 2^{2N} \Big\{ \mathbb{E}_{\mu}\big[ Y_t(\phi_{a,2(\gamma+d-2)})^N \big] +\mathbb{E}_{\mu}\big[ Y_t(\phi_{b,2(\gamma+d-2)})^N \big] \Big\}.
\end{align*}
If $d=2$, then $0\leq \gamma \leq 1\wedge(2-\frac{d}{2})$ implies $2\gamma \in(0,2)$. Hence, by Corollary \ref{'cor5.15} we have
\begin{align*}
	\mathbb{E}_{\mu}\big[ Y_t(\phi_{a,2(\gamma+d-2)})^N \big] = \mathbb{E}_{\mu}\big[ Y_t(\phi_{a,2\gamma})^N \big] \leq C_{N,T,\gamma},\quad \forall\ 0\leq t\leq T \text{ and } a\in\mathbb{R}^2.
\end{align*}
If $d=3$, since $0\leq \gamma \leq 1\wedge(2-\frac{d}{2})$ implies $2(\gamma +1)\in(2,3)\subset(0,3)$, by Corollary \ref{'cor5.18} we get
\begin{align*}
	\mathbb{E}_{\mu}\big[ Y_t(\phi_{a,2(\gamma+d-2)})^N \big] = \mathbb{E}_{\mu}\big[ Y_t(\phi_{a,2(\gamma+1)})^N \big] \leq C_{N,T,\gamma}
\end{align*}
holds for any $0\leq t\leq T$ and $a\in\mathbb{R}^3$ with $|a|\leq T$. Summing up, when $d=2$ or $3$,
\begin{align*}
	\mathbb{E}_{\mu} \Big\{ \big[Y_t((g_{\alpha}^a-g_{\alpha}^b)^2)\big]^N \Big\}\leq C_{N,T,\alpha,\gamma} |a-b|^{2N\gamma}
\end{align*}
holds for any $0\leq t\leq T$ and $a, b\in\mathbb{R}^d$ with $|a|,|b|\leq T$.

 We finally consider the second term on the right-hand side of \eqref{eb.7.17}. 

(i) When $d=1$, by applying H\"older's inequality with $p=N$ and $q=N/(N-1)$ we get
\begin{align}\label{eb.7.21}
	\Big[ \int_0^t \big(X_s(|g_{\alpha}^a-g_{\alpha}^b|) \big)^2 ds \Big]^N \leq t^{N-1} \int_0^t \big(X_s(|g_{\alpha}^a-g_{\alpha}^b|)\big)^{2N} ds.
\end{align}
It then follows from Lemmas \ref{lemx5.1} and \ref{lem.eb.7.1} that for any $ \gamma_1 \in(0,1)$, 
\begin{align*}
	\mathbb{E}_{\mu} \Big\{ \Big[ \int_0^t \big(X_s(g_{\alpha}^a-g_{\alpha}^b) \big)^2 ds \Big]^N \Big\} 
	&\leq t^{N-1} \int_0^t \mathbb{E}_{\mu}\Big[ \big(X_s(g_{\alpha}^a-g_{\alpha}^b)\big)^{2N} \Big]ds \nonumber\\
	&\leq C_{N,\alpha,\gamma_1} t^{N-1} |a-b|^{2N\gamma_1} \int_0^t \mathbb{E}_{\mu} \big[X_s(1)^{2N}\big] ds \nonumber\\
	&\leq C_{N,\alpha,\gamma_1} t^{N}  C_{N,t} |a-b|^{2N\gamma_1} \nonumber\\
	&\leq C_{N,T,\alpha,\gamma_1} |a-b|^{2N\gamma_1},\quad \forall\ 0\leq t\leq T \text{ and } a,b\in\mathbb{R}.
\end{align*}

(ii) When $d=2,3$, we apply Cauchy-Schwarz's inequality twice to get
\begin{align}\label{eb.7.23}
	&\mathbb{E}_{\mu} \Big\{ \Big[ \int_0^t \big(X_s(|g_{\alpha}^a-g_{\alpha}^b|) \big)^2 ds \Big]^N \Big\} \nn\\
	\leq &\ \mathbb{E}_{\mu}\Big[ \Big( \int_{0}^t X_s(1) X_s \big( (g_{\alpha}^a-g_{\alpha}^b)^2 \big) ds \Big)^N \Big] \nonumber\\
	\leq &\ \mathbb{E}_{\mu}\Big[ \Big( \sup_{0\leq s \leq t}X_s(1) \cdot \int_{0}^t  X_s \big( (g_{\alpha}^a-g_{\alpha}^b)^2 \big) ds \Big)^N \Big] \nonumber\\
	\leq &\ \Big\{ \mathbb{E}_{\mu}\Big[ \Big( \sup_{0\leq s \leq t}X_s(1) \Big)^{2N}\Big]\cdot\mathbb{E}_{\mu} \Big[ \Big(\int_0^t X_s \big( (g_{\alpha}^a-g_{\alpha}^b)^2 \big) ds\Big)^{2N}\Big] \Big\}^{1/2} \nonumber\\
	= &\ \Big\{ \mathbb{E}_{\mu}\Big[ \Big( \sup_{0\leq s \leq t}X_s(1) \Big)^{2N}\Big]\cdot\mathbb{E}_{\mu} \Big[ \Big(Y_t \big( (g_{\alpha}^a-g_{\alpha}^b)^2\big) \Big)^{2N}\Big] \Big\}^{1/2}.
\end{align}
By Lemma \ref{lemx5.1} we see that 
\begin{align}\label{bu.7.26}
	\mathbb{E}_{\mu}\Big[ \big( \sup_{0\leq s \leq t}X_s(1) \big)^{2N}\Big] \leq C_{N,T},\quad \forall\  0\leq t \leq T.
\end{align}
For the second term on the right-hand side of \eqref{eb.7.23}, by Lemma \ref{lem.eb.7.1} we get
\begin{align*}
	\mathbb{E}_{\mu} \Big[ \Big(Y_t \big( (g_{\alpha}^a-g_{\alpha}^b)^2 \big)\Big)^{2N}\Big] &\leq C_{N,\alpha,\gamma} |a-b|^{4N\gamma} \mathbb{E}_{\mu} \Big[ \Big(Y_t \big( (\phi_{a,\gamma+d-2}+\phi_{b,\gamma+d-2})^2 \big)\Big)^{2N}\Big] \nonumber\\
	&\leq C_{N,\alpha,\gamma} |a-b|^{4N\gamma} \mathbb{E}_{\mu} \Big[ \Big( Y_t ( \phi_{a,2(\gamma+d-2)}+\phi_{b,2(\gamma+d-2)})\Big)^{2N}\Big] \nonumber\\
	&\leq C_{N,\alpha,\gamma} |a-b|^{4N\gamma} \Big( \mathbb{E}_{\mu}\big[ Y_t ( \phi_{a,2(\gamma+d-2)})^{2N} \big] + \mathbb{E}_{\mu}\big[ Y_t ( \phi_{b,2(\gamma+d-2)})^{2N} \big] \Big)\\
	&\leq C_{N,T,\alpha,\gamma}  |a-b|^{4N\gamma}\end{align*}
holds for any $0\leq t\leq T$ and $a,b\in\mathbb{R}^d$ with $|a|,|b|\leq T$, and the last inequality respectively follows by Corollaries \ref{'cor5.15} and \ref{'cor5.18} when $d=2$ and $d=3$. This together with \eqref{eb.7.23} and \eqref{bu.7.26} gives that for any $0\leq t\leq T$ and $a, b\in\mathbb{R}^d$ with $|a|, |b|\leq T$,
\begin{align*}
	\mathbb{E}_{\mu} \Big\{ \Big[ \int_0^t \big(X_s(|g_{\alpha}^a-g_{\alpha}^b|) \big)^2 ds \Big]^N \Big\} \leq C_{N,T,\alpha,\gamma} |a-b|^{2N\gamma}.
\end{align*}
The desired result follows by summing up the above arguments.
\end{proof}

\subsection{Estimates on the time difference}
 In this section, we turn to give the estimate on the time difference. By the expression \eqref{eb.7.5} of $Y(t,x)$, it suffices to evaluate each term on the right-hand side of this expression. 
 When $d=2$ or $3$, we fix $\gamma_0\in(0,1\wedge(2-\frac{d}{2}))$ and the index $\gamma_0$ is independent of the H\"older exponent in the assumption \eqref{de1.3}. 
 We  still assume $\alpha\geq0$ ($\alpha>0$ if $d=1,2$) throughout this section.

We first estimate the time difference of the third term on the right-hand side of \eqref{eb.7.5}.
\begin{proposition}\label{prop.e.7.3}
	Let $d\leq 3$, $N\geq 1$ and $T>0$. 	
	\begin{itemize}
		\item [\rm{(i)}] When $d=1$, there is some constant $C_{N,T,\alpha}>0$ such that for any $0\leq s, t\leq T$ and $a\in\mathbb{R}$,
			\begin{align*}
				\mathbb{E}_{\mu}\Big[ \Big| \int_0^t X_{r}(g_{\alpha}^a)dr - \int_0^s X_{r}(g_{\alpha}^a)dr  \Big|^{2N} \Big] \leq C_{N,T,\alpha} |t-s|^{2N}.
			\end{align*}

		\item [\rm{(ii)}] When $d=2$ or $3$, we further assume \eqref{de.1.5} holds. Then for any $\gamma_0\in(0,1\wedge(2-\frac{d}{2}))$, there exists a constant $C_{N,T,\alpha,\gamma_0}>0$ such that 
			\begin{align*}
				\mathbb{E}_{\mu}\Big[ \Big| \int_0^t X_{r}(g_{\alpha}^a)dr - \int_0^s X_{r}(g_{\alpha}^a)dr  \Big|^{2N} \Big] \leq C_{N,T,\alpha,\gamma_0} |t-s|^{N\gamma_0}
			\end{align*}
			holds for all $0\leq s, t\leq T$ and $a\in\mathbb{R}^d$ with $|a|\leq T$.
	\end{itemize}
\end{proposition}

\begin{proof}
	Let $N\geq1$. For any $s, t\geq 0$ and $a\in\mathbb{R}^d$,
	\begin{align}\label{eb.7.37}
		\mathbb{E}_{\mu}\Big[ \Big| \int_0^t X_{r}(g_{\alpha}^a)dr - \int_0^s X_{r}(g_{\alpha}^a)dr  \Big|^{2N} \Big]
		&=\mathbb{E}_{\mu}\Big[ \Big|\int_s^t 	X_r(g_{\alpha}^a) dr\Big|^{2N}\Big] \nn\\
		&\leq C_{N,\alpha} \mathbb{E}_{\mu}\Big[ \Big|\int_s^t	 X_r(g_{0}^a) dr\Big|^{2N}\Big],
	\end{align}
	where the inequality follows by \eqref{e.7.0}. We next fix $T>0$.
	
	(i) When $d=1$, by the definition of $g_0^a$ we see $g_0^a(x)\equiv 1$. Hence, for any $0\leq s, t \leq T$,
	\begin{align*}
		\mathbb{E}_{\mu}\Big[ \Big|\int_s^t	 X_r(g_{0}^a) dr\Big|^{2N}\Big] 
		&= \mathbb{E}_{\mu}\Big[ \Big|\int_s^t	 X_r(1) dr\Big|^{2N}\Big]	\nn\\
		&\leq |t-s|^{2N} \mathbb{E}_{\mu}\big[ (\sup_{0\leq r \leq T}X_r(1))^{2N} \big] \nn\\
		&\leq C_{N,T} |t-s|^{2N},
	\end{align*}
	where the last inequality comes from Lemma \ref{lemx5.1}. Thus, Proposition \ref{prop.e.7.3}(i) follows by the preceding two inequalities.
		
	(ii) When $d=2,3$, for any $\gamma_0\in(0,1\wedge(2-\frac{d}{2}))$, we apply H\"older's inequality with $p=2/\gamma_0>1$ and $q=2/(2-\gamma_0)>1$ to obtain that for any $0\leq s, t\leq T$,
	\begin{align}\label{xiu.e.90}
		\Big|\int_s^t X_r(g_0^a) dr \Big|^{2N} 
		&\leq \Big| \int_s^t X_r(1) dr\Big|^{N\gamma_0} \times \Big| \int_s^t dr\int (g_0^a(x))^{2/(2-\gamma_0)} X_r(dx) \Big|^{N(2-\gamma_0)} \nn\\
		&\leq |t-s|^{N\gamma_0}\Big(\sup_{0\leq r\leq T} X_r(1) \Big)^{N\gamma_0}  \Big[Y_T\big( (g_0^a)^{2/(2-\gamma_0)} \big)\Big]^{N(2-\gamma_0)}.
	\end{align}
	It then follows from Cauchy-Schwarz's inequality that
	\begin{align*}
		\mathbb{E}_{\mu}&\Big[ \Big|\int_s^t	X_r(g_{0}^a) dr\Big|^{2N}\Big] \nn\\
		&\leq |t-s|^{N\gamma_0} \bigg\{\mathbb{E}_{\mu} \Big[ \big(\sup_{0\leq r\leq T} X_r(1) \big)^{2N\gamma_0}\Big] \times \mathbb{E}_{\mu}\Big[ \Big(Y_T\big((g_0^a)^{2/(2-\gamma_0)}\big)\Big)^{2N(2-\gamma_0)}\Big] \bigg\}^{1/2} \nn\\
		&\leq |t-s|^{N\gamma_0} \bigg\{\mathbb{E}_{\mu} \Big[ \big(\sup_{0\leq r\leq T} X_r(1) \big)^{2N}\Big] \bigg\}^{\gamma_0/2} \bigg\{\mathbb{E}_{\mu}\Big[ \Big(Y_T\big((g_0^a)^{2/(2-\gamma_0)}\big)\Big)^{4N}\Big] \bigg\}^{(2-\gamma_0)/4},
	\end{align*}
	where we use Jensen's inequality twice to get the last inequality as we note $\gamma_0<1$ and $(2-\gamma_0)/2<1$. We again apply Lemma \ref{lemx5.1} to see
	\begin{align}\label{eb.7.39}
		\mathbb{E}_{\mu}\Big[ \Big|\int_s^t	X_r(g_{0}^a) dr\Big|^{2N}\Big] 
		\leq C_{N,T,\gamma_0} |t-s|^{N\gamma_0} \bigg\{\mathbb{E}_{\mu}\Big[ \Big(Y_T\big((g_0^a)^{\widetilde{\gamma}}\big)\Big)^{4N}\Big] \bigg\}^{(2-\gamma_0)/4},
	\end{align}
	where  
	\begin{align}\label{eb.7.42}
			\widetilde{\gamma}:=2/(2-\gamma_0)\in(1,\frac{2}{2-1\wedge(2-\frac{d}{2})}).
	\end{align}
	If $d=2$, then $\widetilde{\gamma}\in(1,2)$. Recall that $\phi_{a,\widetilde{\gamma}}$ is given by \eqref{phig}, 	it then follows from \eqref{'bu5.63} to obtain
	\begin{align}\label{eb.7.43}
		\mathbb{E}_{\mu}\Big[ \Big(Y_T\big((g_0^a)^{\widetilde{\gamma}}\big)\Big)^{4N}\Big] 
		&= \mathbb{E}_{\mu}\Big[ \Big(\int \Big(1+\log^+\frac{1}{|a-y|}\Big)^{\widetilde{\gamma}} Y_T(dy)\Big)^{4N} \Big] \nn\\
		&\leq \mathbb{E}_{\mu}\Big[ \Big(\int \Big(1+\frac{1}{|a-y|}\Big)^{\widetilde{\gamma}} Y_T(dy)\Big)^{4N} \Big] \nn\\
		&\leq \mathbb{E}_{\mu}\Big[ \Big(2\int \Big(1+\frac{1}{|a-y|^{\widetilde{\gamma}}}\Big) Y_T(dy)\Big)^{4N} \Big] \nn \\
		&= 2^{4N} \mathbb{E}_{\mu} \Big[ \big(Y_T(1)+Y_T(\phi_{a,\widetilde{\gamma}})\big)^{4N}\Big] \nn\\
		&\leq 2^{8N} \Big( \mathbb{E}_{\mu}\big[Y_T(1)^{4N}\big]+\mathbb{E}_{\mu}\big[Y_T(\phi_{a,\widetilde{\gamma}})^{4N}\big]\Big) \nn\\
		&\leq C_{N,T,\gamma_0},
 	\end{align}
	where the last inequality comes from Lemma \ref{lem5.5} and Corollary \ref{'cor5.15}. If $d=3$, since \eqref{de.1.5} holds in this case, by applying Corollary \ref{'cor5.18} with $\widetilde{\gamma}\in(1,4/3)$ we get that for any $a\in\mathbb{R}^3$ with $|a|\leq T$,
	\begin{align*}
		\mathbb{E}_{\mu} \Big[ \Big(Y_T\big((g_0^a)^{\widetilde{\gamma}}\big)\Big)^{4N}\Big] = \mathbb{E}_{\mu} \Big[\Big(\int \Big(\frac{1}{2\pi |a-y|}\Big)^{\widetilde{\gamma}} Y_T(dy)\Big)^{4N} \Big] \leq \mathbb{E}_{\mu}\big[Y_T(\phi_{a,\widetilde{\gamma}})^{4N}\big] \leq C_{N,T,\gamma_0}.
	\end{align*}
	This combined with \eqref{eb.7.39} and \eqref{eb.7.43} implies that for any $a\in\mathbb{R}^d$ with $|a|\leq T$,
	\begin{align*}
		\mathbb{E}_{\mu}\Big[ \Big|\int_s^t	X_r(g_{0}^a) dr\Big|^{2N}\Big] \leq C_{N,T,\gamma_0} |t-s|^{N\gamma_0},	
	\end{align*}
	which together with \eqref{eb.7.37} gives Proposition \ref{prop.e.7.3}(ii).
\end{proof}

We next estimate the time difference of the fourth term on the right-hand side of \eqref{eb.7.5}. To do this, we start with the stochastic integral $M_t(g_{\alpha}^a)$ for any $a\in\mathbb{R}^d$. By Lemma \ref{lem7.4}, if \eqref{de.1.5} holds for $d=2, 3$, then $(M_t(g_{\alpha}^a), t\geq 0)$ is a continuous $(\mathcal{F}_t)$-martingale with
\begin{align*}
	\langle M(g_{\alpha}^a)\rangle_t=\int_0^t X_r((g_{\alpha}^a)^2)dr+\int_0^t dr \int\!\!\int g(u,v) g_{\alpha}^a(u)g_{\alpha}^a(v) X_r(du) X_r(dv).
\end{align*}
Obviously, $(M_t(g_{\alpha}^a)-M_s(g_{\alpha}^a), t\geq s)$ is a continuous $(\mathcal{F}_t)$-martingale and for any $t\geq s$, its quadratic variation
\begin{align*}
	\langle M_{\cdot}(g_{\alpha}^a)-M_{s}(g_{\alpha}^a)\rangle_t=\int_s^t X_r((g_{\alpha}^a)^2)dr+\int_s^t dr \int\!\!\int g(u, v) g_{\alpha}^a(u)g_{\alpha}^a(v) X_r(du) X_r(dv).
\end{align*}
Based on this, we get the following result, which describes the time difference about $M_{t}(g_{\alpha}^a)$.
\begin{proposition}\label{prop.e.7.4}
	Let $d\leq 3$, $N\geq 1$ and $T>0$.
	\begin{itemize}
		\item [\rm{(i)}] When $d=1$, for any $0\leq s, t\leq T$ and $a\in\mathbb{R}$, we have
		\begin{align*}
			\mathbb{E}_{\mu} \big[ \big|M_t(g_{\alpha}^a)-M_s(g_{\alpha}^a)\big|^{2N} \big] \leq C_{N,T,\alpha} |t-s|^{N}.
		\end{align*}
		for some constant $C_{N,T,\alpha}>0$.
		
		\item [\rm{(ii)}] When $d=2$ or $3$, we further assume \eqref{de.1.5} holds. Then for any $\gamma_0\in(0,1\wedge(2-\frac{d}{2}))$, there is a constant $C_{N,T,\alpha,\gamma_0}>0$ (depends on the dimension $d$) such that
		\begin{equation*}
			\mathbb{E}_{\mu} \big[ \big|M_t(g_{\alpha}^a)-M_s(g_{\alpha}^a)\big|^{2N} \big]\leq 	
			\begin{cases}
				C_{N,T,\alpha,\gamma_0}  |t-s|^{N\gamma_0/2}, & \text{if } d=2; \\
				C_{N,T,\alpha,\gamma_0} |t-s|^{2N\gamma_0/3}, & \text{if } d=3,
			\end{cases}
		\end{equation*}
	holds for any $0\leq s, t \leq T$ and $a\in\mathbb{R}^d$ with $|a|\leq T$.
	\end{itemize}
\end{proposition}

\begin{proof}
	Fix $N\geq1$ and $T>0$. Without loss of generality, we may assume $0\leq s\leq t\leq T$. By Burkholder-Davis-Gundy's inequality we obtain that for any $a\in\mathbb{R}^d$,
	\begin{align*}
		\mathbb{E}_{\mu} &\big[ \big|M_t(g_{\alpha}^a)-M_s(g_{\alpha}^a)\big|^{2N} \big] \nn \\
		&\leq C_N \mathbb{E}_{\mu} \big[ \langle M_{\cdot}(g_{\alpha}^a)-M_{s}(g_{\alpha}^a)\rangle_t^N \big] \nn\\
		&= C_N \mathbb{E}_{\mu}\Big[ \Big( \int_s^t X_r((g_{\alpha}^a)^2)dr+\int_s^t dr \int\!\!\int g(u, v) g_{\alpha}^a(u)g_{\alpha}^a(v) X_r(du) X_r(dv) \Big)^{N} \Big] \nn\\
		&\leq 2^{N}C_N \bigg\{ \mathbb{E}_{\mu}\Big[\Big( \int_s^t X_r((g_{\alpha}^a)^2)dr\Big)^{N}\Big] + \|g\|_{\infty}^N \mathbb{E}_{\mu}\Big[\Big(\int_s^t \big(X_r(g_{\alpha}^a)\big)^2 dr\Big)^{N}\Big] \bigg\}.
	\end{align*}
	For the second term, by applying Cauchy-Schwarz's inequality twice and using Lemma \ref{lemx5.1} and \eqref{e.7.0} we get
	\begin{align}\label{xiu.e.57}
		\mathbb{E}_{\mu}	\Big[\Big(\int_s^t \big(X_r(g_{\alpha}^a)\big)^2 dr\Big)^{N}\Big]
		&\leq \mathbb{E}_{\mu}\Big[\Big(\int_s^t X_r(1)X_r\big((g_{\alpha}^a)^2\big) dr\Big)^{N}\Big] \nn\\
		&\leq \mathbb{E}_{\mu}\Big[\Big(\sup_{0\leq r \leq T}X_r(1)\Big)^{N}\Big(\int_s^t X_r\big((g_{\alpha}^a)^2\big) dr\Big)^{N}\Big] \nn\\
		&\leq \bigg\{ \mathbb{E}_{\mu}\Big[\Big(\sup_{0\leq r \leq T}X_r(1)\Big)^{2N} \Big] \mathbb{E}_{\mu}\Big[\Big(\int_s^t X_r\big((g_{\alpha}^a)^2\big) dr\Big)^{2N}\Big] \bigg\}^{1/2} \nn\\
		&\leq C_{N,T,\alpha} \bigg\{ \mathbb{E}_{\mu}\Big[\Big(\int_s^t X_r\big((g_{0}^a)^2\big) dr\Big)^{2N}\Big]  \bigg\}^{1/2}.
	\end{align}
	Turning to the first term, by again applying Cauchy-Schwarz's inequality and \eqref{e.7.0} we see
	\begin{align*}
		\mathbb{E}_{\mu}\Big[\Big( \int_s^t X_r((g_{\alpha}^a)^2)dr\Big)^{N}\Big] \leq C_{N,\alpha}\bigg\{ \mathbb{E}_{\mu}\Big[\Big(\int_s^t X_r\big((g_{0}^a)^2\big) dr\Big)^{2N}\Big]  \bigg\}^{1/2}.
	\end{align*}
	Combining the three preceding inequalities gives that for any $0\leq s\leq t\leq T$,
	\begin{align}\label{eb.7.50}
		\mathbb{E}_{\mu} \big[ \big|M_t(g_{\alpha}^a)-M_s(g_{\alpha}^a)\big|^{2N} \big] 
		\leq C_{N,T,\alpha} \bigg\{ \mathbb{E}_{\mu}\Big[\Big|\int_s^t X_r\big((g_{0}^a)^2\big) dr\Big|^{2N}\Big]  \bigg\}^{1/2}.
	\end{align}
	By a similar argument, we get \eqref{eb.7.50} holds for $0\leq s, t\leq T$. So it suffices to consider the expectation on the right-hand side of \eqref{eb.7.50} for any $0\leq s, t\leq T$.
	
	(i) When $d=1$, we use Lemma \ref{lemx5.1} to get 
	\begin{align}\label{6.53}
		\mathbb{E}_{\mu}\Big[\Big|\int_s^t X_r\big((g_{0}^a)^2\big) dr\Big|^{2N}\Big] 
		&= \mathbb{E}_{\mu} \Big[ \Big|\int_s^t X_r(1) dr\Big|^{2N}\Big] \nn\\
		&\leq |t-s|^{2N} \mathbb{E}_{\mu}\Big[\big(\sup_{0\leq r \leq T}X_r(1)\big)^{2N}\Big] \nn\\
		&\leq C_{N,T} |t-s|^{2N}.
	\end{align}
	Hence, Proposition \ref{prop.e.7.4}(i) follows by the above inequality and \eqref{eb.7.50}.
	
	(ii) When $d=2$ or $3$, we fix $\gamma_0\in(0, 1\wedge(2-\frac{d}{2}))$. If $d=2$, by repeating arguments for deriving \eqref{eb.7.39} and replacing $g_0^{a}$ with $(g_0^a)^2$ we get that for any fixed $\gamma_0\in(0,1)$, there is a constant $C_{N,T,\gamma_0}>0$ such that
	\begin{align}\label{eb.7.53}
		\mathbb{E}_{\mu}\Big[\Big|\int_s^t X_r\big((g_{0}^a)^2\big) dr\Big|^{2N}\Big] \leq C_{N,T,\gamma_0} |t-s|^{N \gamma_0} \bigg\{ \mathbb{E}_{\mu}\Big[ \Big(Y_T\big((g_0^a)^{2\widetilde{\gamma}}\big)\Big)^{4N}\Big]\bigg\}^{(2-\gamma_0)/4},
	\end{align}
	where $\widetilde{\gamma}=2/(2-\gamma_0)\in(1,2)$ is given by \eqref{eb.7.42}. Note that $g_0^a(x)=1+\log^+\frac{1}{|a-x|}$ for $d=2$. It then follows from \eqref{'bu5.63} that
	\begin{align*}
		\big(g_0^a(x)\big)^{2\widetilde{\gamma}}=\Big[\Big(1+\log^+\frac{1}{|a-x|}\Big)^2\Big]^{\widetilde{\gamma}} \leq \Big(1+\frac{3}{|a-x|}\Big)^{\widetilde{\gamma}}.	
	\end{align*}
	Hence, similar to \eqref{eb.7.43}, we easily get
	\begin{align*}
		\mathbb{E}_{\mu}\Big[ \Big(Y_T\big((g_0^a)^{2\widetilde{\gamma}}\big)\Big)^{4N}\Big] 
		&\leq \mathbb{E}_{\mu}\Big[\Big(\int\Big(1+\frac{3}{|a-y|}\Big)^{\widetilde{\gamma}} Y_T(dy)\Big)^{4N}\Big] \nn\\
		&\leq C_{N}\Big( \mathbb{E}_{\mu}\big[Y_T(1)^{4N}\big] + \mathbb{E}_{\mu}\big[Y_T(\phi_{a,\widetilde{\gamma}})^{4N} \big] \Big)
		\leq C_{N,T,\gamma_0},
	\end{align*}
	where $\phi_{a,\widetilde{\gamma}}$ is defined by \eqref{phig}. This combined with \eqref{eb.7.53} gives
	\begin{align}\label{6.57}
		\mathbb{E}_{\mu}\Big[\Big|\int_s^t X_r\big((g_{0}^a)^2\big) dr\Big|^{2N}\Big] \leq C_{N,T,\gamma_0} |t-s|^{N \gamma_0},
	\end{align}
	which together with \eqref{eb.7.50} shows that Proposition \ref{prop.e.7.4}(ii) holds for $d=2$. 
	
	If $d=3$, then $\gamma_0\in(0,1/2)$ and $g_0^a(x)=1/(2\pi|a-x|)$. Similar to \eqref{xiu.e.90}, but here we apply H\"older's inequality with $p=3/(2\gamma_0)$ and $q=3/(3-2\gamma_0)$ to obtain that for any $0\leq s, t\leq T$,
	\begin{align*}
		 \Big|\int_s^t X_r\big((g_{0}^a)^2\big) dr\Big|^{2N}
		&\leq  |t-s|^{4N\gamma_0/3} \Big(\sup_{0\leq r \leq T}X_r(1)\Big)^{4N\gamma_0/3}\Big[ Y_T\big( (g_0^a)^{6/(3-2\gamma_0)}\big) \Big]^{2N(3-2\gamma_0)/3} \nn\\
		&\leq   |t-s|^{4N\gamma_0/3} \Big(\sup_{0\leq r \leq T}X_r(1)\Big)^{4N\gamma_0/3} \Big(\int\frac{1}{|a-y|^{6/(3-2\gamma_0)}}Y_T(dy)\Big)^{2N(3-2\gamma_0)/3}.	
	\end{align*}
	Let $\bar{\gamma}:=6/(3-2\gamma_0)$. Then $\bar{\gamma}\in(2,3)$. By Cauchy-Schwarz's inequality and Jensen's inequality, we get that if the assumption \eqref{de.1.5} holds for $d=3$, 
	\begin{align}\label{6.59}
		&\mathbb{E}_{\mu}\Big[\Big|\int_s^t X_r\big((g_{0}^a)^2\big) dr\Big|^{2N}\Big] \nn\\
		\leq &\ |t-s|^{4N\gamma_0/3} \mathbb{E}_{\mu}\Big[\Big(\sup_{0\leq r \leq T}X_r(1)\Big)^{4N\gamma_0/3} \Big(Y_T(\phi_{a,\bar{\gamma}})\Big)^{2N(3-2\gamma_0)/3}\Big]	\nn\\
		\leq &\ |t-s|^{4N\gamma_0/3} \bigg\{\mathbb{E}_{\mu}\Big[\Big(\sup_{0\leq r \leq T}X_r(1)\Big)^{8N\gamma_0/3} \Big] \mathbb{E}_{\mu}\Big[\Big(Y_T(\phi_{a,\bar{\gamma}})\Big)^{4N(3-2\gamma_0)/3}\Big]\bigg\}^{1/2} \nn\\
		\leq &\ |t-s|^{4N\gamma_0/3} \bigg\{\mathbb{E}_{\mu}\Big[\Big(\sup_{0\leq r\leq T}X_r(1)\Big)^{8N}\Big]\bigg\}^{\gamma_0/6} \bigg\{\mathbb{E}_{\mu}\Big[Y_T(\phi_{a,\bar{\gamma}})^{4N}\Big]\bigg\}^{(3-2\gamma_0)/6} \nn\\
		\leq &\ C_{N,T,\gamma_0} |t-s|^{4N\gamma_0/3},\quad \forall\ a\in\mathbb{R}^3 \text{ with }|a|\leq T,	\end{align}
where the last inequality is obtained by Corollary \ref{'cor5.18} and Lemma \ref{lemx5.1}. By taking the above inequality into \eqref{eb.7.50} we prove Proposition \ref{prop.e.7.4}(ii) holds for $d=3$.
\end{proof}

At the end of this section, we shall consider the estimate on the time difference about $X_t(g_{\alpha}^a)$. Recall the Green function representation \eqref{gfr}:
\begin{align*}
		X_t(g^a_\alpha) = X_0(P_t g^a_\alpha) + \int_0^t\int_{\mathbb{R}^d} P_{t-s}g^a_\alpha(x) dM(s, x),\quad \mathbb{P}_{\mu}\text{-a.s.},	
	\end{align*}
	which is an important tool given in Theorem \ref{cor3.4*}.	Based on this, in order to estimate the time difference of the second term on the right-hand side of (\ref{eb.7.5}), we only need to consider the two terms on the right-hand side of (\ref{gfr}).
	
\begin{lemma}\label{*lem6.1}
	Let $d\le 3$ and $p_t(0,x)$ be the transition density defined by \eqref{e10.21}. For any $0\le \tilde{\delta}\le 1$ we have
	\begin{align*}
		\big|p_t(0,x)-p_s(0,x)\big| \leq\Big[(t-s) s^{-\frac{d}{2}-1}\Big]^{\tilde{\delta}}\Big[p_t(0,x)^{1-\tilde{\delta}}+p_s(0,x)^{1-\tilde{\delta}}\Big], \quad \forall\ 0<s \leq t  \text { and } x \in \mathbb{R}^d .
	\end{align*}
\end{lemma} 

\begin{proof}
	The result follows by taking $p=1$ in Lemma III.4.5(a) of \cite{P02} when $d=1$. Similarly, when $d=2$ or $d=3$, by the mean value theorem, there exists a constant $u\in (s,t)$ such that
	\begin{align*}
		\big|p_t(0,x)-p_s(0,x)\big|=(t-s)\Big|\frac{\partial p_u(0,x)}{\partial u}\Big|=(t-s) \frac{p_u(0,x)}{2 u}\Big|\frac{|x|^2}{u}-d\Big|.
	\end{align*}
	Notice that
	\begin{align*}
		\frac{u^{d / 2}}{2} p_u(0,x)\Big|\frac{|x|^2}{u}-d\Big|= \Big(\frac{1}{2 \pi}\Big)^{d / 2} e^{-\frac{|x|^2}{2 u}}\Big|\frac{|x|^2}{2u}-\frac{d}{2}\Big|\le \Big(\frac{1}{2}\Big)^{d/2}\sup _{r \geq 0}\Big|e^{-r}\Big(r-\frac{d}{2}\Big)\Big|\le 1.
	\end{align*}
	Then we have $|p_t(0,x)-p_s(0,x)|\le (t-s)u^{-d / 2-1}$. For any $0\le \tilde{\delta}\le 1$, we further get
	\begin{align*}
\big|p_t(0,x)-p_s(0,x)\big| \leq\Big[(t-s) s^{-\frac{d}{2}-1}\Big]^{\tilde{\delta}}\big|p_t(0,x)-p_s(0,x)\big|^{1-\tilde{\delta}},
\end{align*}
which follows the desired result.
\end{proof}

By an application of Lemma \ref{*lem6.1}, we get the following time difference about the first term on the right-hand side of \eqref{gfr}.

\begin{lemma}\label{bu.lem6.1}
	Let $d\le 3$, $N\ge 1$ and $\varepsilon_0>0$. There exists a constant $C_{N,\varepsilon_0,\alpha}$ such that
	\begin{align*}
		\mathbb{E}_\mu\Big[\big|X_0(P_t g_\alpha^a)-X_0(P_s g_\alpha^a)\big|^{2 N}\Big] \leq C_{N, \varepsilon_0, \alpha}(t-s)^{2 N},\quad \forall\ s, t \geq \varepsilon_0 \text { and }\ a \in \mathbb{R}^d .
	\end{align*}
\end{lemma}

\begin{proof}
	Fix $N\geq1$, $\eps_0>0$ and $a\in\mathbb{R}^d$. Without loss of generality, we may assume $\varepsilon_0\le s\le t$. The case $\varepsilon_0\le t\le s$ follows by a similar argument. By H\"older's inequality, we obtain that
	\begin{align}\label{'bu.6.1}
		\mathbb{E}_\mu\Big[\big|X_0(P_t g_\alpha^a)-X_0(P_s g_\alpha^a)\big|^{2 N}\Big]
		&=\big|\langle\mu, P_t g_\alpha^a-P_s g_\alpha^a\rangle\big|^{2 N}\nn\\
		&\leq \mu(1)^{2 N-1} \int\big|P_t g_\alpha^a(x)-P_s g_\alpha^a(x)\big|^{2 N} \mu(d x) .
	\end{align}
	We apply Lemma \ref{*lem6.1} with $\tilde{\delta}=1$ to get
	\begin{align*}
\big|P_t g_\alpha^a(x)-P_s g_\alpha^a(x)\big| 
& =\Big|\int\big(p_t(x, y)-p_s(x, y)\big) d y \int_0^{\infty} e^{-\alpha r} p_r(a, y) d r \Big|\nn \\
& \leq \int_0^{\infty} e^{-\alpha r}\big|p_{t+r}(x, a)-p_{s+r}(x, a)\big| d r\nn\\
&\le 2(t-s)\int_0^{\infty} e^{-\alpha r}(s+r)^{-\frac{d}{2}-1}dr\nn\\
&\le 2(t-s)\Big[\int_0^1 (s+r)^{-\frac{d}{2}-1}dr+\int_1^{\infty}e^{-\alpha r} dr\Big]\le C_{\varepsilon_0,\alpha}(t-s),
\end{align*}
where $C_{\varepsilon_0,\alpha}$ is a constant depending on $d$. Then combining the above inequality with (\ref{'bu.6.1}) yields the desired result.
\end{proof}

We now consider the time difference of the second term on the right-hand side of (\ref{gfr}). Without loss of generality, we may assume $0\le s\le t$. The case $0\le t\le s$ follows by a similar argument. For any $N\geq1$ and $a\in\mathbb{R}^d$,
	\begin{align}\label{g.2}
		&\mathbb{E}_\mu\Big[\Big|\int_0^t\int P_{t-r}g^a_\alpha(x) dM(r,x)-\int_0^s\int P_{s-r}g^a_\alpha(x) dM(r,x)\Big|^{2 N}\Big]\nn\\
		\le&\ 2^{2N}\Big\{\mathbb{E}_\mu\Big[\Big|\int_0^s \int\big(P_{t-r} g_\alpha^a(x)-P_{s-r} g_\alpha^a(x)\big) d M(r, x)\Big|^{2 N}\Big]+\mathbb{E}_\mu\Big[\Big|\int_s^t \int P_{t-r} g_\alpha^a(x) d M(r, x)\Big|^{2 N}\Big]\Big\}\nn\\
		:=&\ 2^{2N}\big\{I_1^{N,\alpha,a}(s,t)+I_2^{N,\alpha,a}(s,t)\big\}.
	\end{align}
Then we only need to estimate $I_1^{N,\alpha,a}$ and $I_2^{N,\alpha,a}$. We present the proofs of these estimates in Appendix \ref{appi} and only state the results here as the following two lemmas.
 
 \begin{lemma}\label{bu.lem.b.2}
	Let $I_1^{N,\alpha,a}(s,t)$ be given by (\ref{g.2}), $N\ge 1$ and $T>0$.
	\begin{itemize}
		\item[\rm{(i)}] When $d=1$, for any $\tilde{\delta}\in (0,1)$, there exists a constant $C_{N,T,\alpha,\tilde{\delta}}>0$ such that  $I_1^{N,\alpha,a}(s,t)\le C_{N,T,\alpha,\tilde{\delta}}(t-s)^{2N\tilde{\delta}}$ for $0\le s\le t\le T$ and $a\in\mathbb{R}$.
		\item[\rm{(ii)}] When $d=2$, for any $\tilde{\delta}\in (0,1/2)$, there exists a constant $C_{N,T,\alpha,\tilde{\delta}}>0$ such that   $I_1^{N,\alpha,a}(s,t)\le C_{N,T,\alpha,\tilde{\delta}}(t-s)^{2N\tilde{\delta}}$ for $0\le s\le t\le T$ and $a\in\mathbb{R}^2$.
		\item[\rm{(iii)}] When $d=3$, suppose (\ref{de.1.5}) holds. Then for any $\tilde{\delta}\in (0, 1/6)$, there exists a constant $C_{N,T,\alpha,\tilde{\delta}}>0$ such that   $I_1^{N,\alpha,a}(s,t)\le C_{N,T,\alpha,\tilde{\delta}}(t-s)^{N\tilde{\delta}}$ for $0\le s\le t\le T$ and $a\in\mathbb{R}^3$ with $|a|\le T$.
	\end{itemize}
\end{lemma}

\begin{lemma}\label{bu.lem.b.3}
	Let $I_2^{N,\alpha,a}(s,t)$ be given by (\ref{g.2}), $N\ge 1$ and $T>0$.
	\begin{itemize}
		\item[\rm{(i)}] When $d=1$, there exists a constant $C_{N,T,\alpha}>0$ such that  $I_2^{N,\alpha,a}(s,t)\le C_{N,T,\alpha}(t-s)^{N}$ for any $0\le s\le t\le T$ and $a\in\mathbb{R}$.
		\item[\rm{(ii)}] When $d=2$, for any $\gamma_0 \in(0,1)$, there exists a constant $C_{N,T,\alpha,\gamma_0}>0$ such that   $I_2^{N,\alpha,a}(s,t)\le C_{N,T,\alpha,\gamma_0}(t-s)^{N\gamma_0/2}$ for any $0\le s\le t\le T$ and $a\in\mathbb{R}^2$.
		\item[\rm{(iii)}] When $d=3$, suppose (\ref{de.1.5}) holds. Then for any $\gamma_0 \in(0,1/2)$, there exists a constant $C_{N,T,\alpha,\gamma_0}>0$ such that   $I_2^{N,\alpha,a}(s,t)\le C_{N,T,\alpha,\gamma_0}(t-s)^{2N\gamma_0/3}$ for any $0\le s\le t\le T$ and $a\in\mathbb{R}^3$ with $|a|\le T$.
	\end{itemize}
\end{lemma}	
 
 Combining the Green function representation \eqref{gfr} with \eqref{g.2}, Lemmas \ref{bu.lem6.1}, \ref{bu.lem.b.2} and \ref{bu.lem.b.3} yields the following result.
 
\begin{proposition}\label{bu.prop.6.4}
	Let $N\ge 1$ and $0<\varepsilon_0<T$. Suppose (\ref{de.1.5}) holds when $d=2$ or $3$.
	\begin{itemize}
		\item[\rm{(i)}] When $d=1$, for every $\tilde{\delta}\in (0,1)$, there exists a constant $C_{N,\varepsilon_0,T,\alpha,\tilde{\delta}}>0$ such that for any $\varepsilon_0\le s,t\le T$ and $a\in\mathbb{R}$, 
		$$\mathbb{E}_\mu\Big[\big|X_t(g_\alpha^a)-X_s(g_\alpha^a)\big|^{2 N}\Big] \leq C_{N, \varepsilon_0, T, \alpha, \tilde{\delta}}|t-s|^{N \tilde{\delta}}.$$
		\item[\rm{(ii)}] When $d=2$, for every $\tilde{\delta}\in (0,1/2)$, there exists a constant $C_{N,\varepsilon_0,T,\alpha,\tilde{\delta}}>0$ such that for any $\varepsilon_0\le s,t\le T$ and $a\in\mathbb{R}^2$, 
		$$\mathbb{E}_\mu\Big[\big|X_t(g_\alpha^a)-X_s(g_\alpha^a)\big|^{2 N}\Big] \leq C_{N, \varepsilon_0, T, \alpha, \tilde{\delta}}|t-s|^{N \tilde{\delta}}.$$	
		\item[\rm{(iii)}] When $d=3$, for every $\tilde{\delta} \in(0, 1/6)$, there exists a constant $C_{N,\varepsilon_0,T,\alpha,\tilde{\delta}}>0$ such that for any $\varepsilon_0\le s,t\le T$ and $a\in\mathbb{R}^3$ with $|a|\le T$, 
		$$\mathbb{E}_\mu\Big[\big|X_t(g_\alpha^a)-X_s(g_\alpha^a)\big|^{2 N}\Big] \leq C_{N, \varepsilon_0, T, \alpha, \tilde{\delta}}|t-s|^{N \tilde{\delta}}.$$
	\end{itemize}
\end{proposition}

	Based on the estimates in Propositions \ref{prop.e.7.3}, \ref{prop.e.7.4} and \ref{bu.prop.6.4}, we derive the estimates on the time difference of \eqref{eb.7.5}.

\begin{theorem}\label{bu.th.6.5}
	Let $N\ge 1$ and $0<\varepsilon_0<T$. Suppose that (\ref{de.1.5}) holds when $d=2$ or $3$.
	\begin{itemize}
		\item[\rm{(i)}] When $d=1$, for every $\tilde{\delta}\in (0,1)$, there exists a constant $C_{N,\varepsilon_0,T,\alpha,\tilde{\delta}}>0$ such that for any $\varepsilon_0\le s,t\le T$ and $a\in\mathbb{R}$, 
		$$\mathbb{E}_\mu\Big[\big|Y(t,a)-Y(s,a)\big|^{2 N}\Big] \leq C_{N, \varepsilon_0, T, \alpha, \tilde{\delta}} |t-s|^{N \tilde{\delta}}.$$
		\item[\rm{(ii)}] When $d=2$, for every $\tilde{\delta}\in (0,1/2)$, there exists a constant $C_{N,\varepsilon_0,T,\alpha,\tilde{\delta}}>0$ such that for any $\varepsilon_0\le s,t\le T$ and $a\in\mathbb{R}^2$ with $|a|\le T$, 
		$$\mathbb{E}_\mu\Big[\big|Y(t,a)-Y(s,a)\big|^{2 N}\Big] \leq C_{N, \varepsilon_0, T, \alpha, \tilde{\delta}} |t-s|^{N \tilde{\delta}}.$$	
		\item[\rm{(iii)}] When $d=3$, for every $\tilde{\delta} \in(0, 1/6)$, there exists a constant $C_{N,\varepsilon_0,T,\alpha,\tilde{\delta}}>0$ such that for any $\varepsilon_0\le s,t\le T$ and $a\in\mathbb{R}^3$ with $|a|\le T$, 
		$$\mathbb{E}_\mu\Big[\big|Y(t,a)-Y(s,a)\big|^{2 N}\Big] \leq C_{N, \varepsilon_0, T, \alpha, \tilde{\delta}} |t-s|^{N \tilde{\delta}}.$$
	\end{itemize}
\end{theorem}

At the end of the section, we shall give the proof of the assertion (ii) in Theorem \ref{t1}.\\
	
\noindent\textbf{Proof of (ii) in Theorem \ref{t1}.}	Let $N\ge 1$ and $0<\varepsilon_0<T$. Suppose the assumptiom \eqref{de1.3} holds when $d=2$ or $3$ and let $\gamma\in (0,1\wedge(2-\frac{d}{2}))$ be given in this assumption. For any $s,t\in[\varepsilon_0,T]$ and $a,b\in \mathbb{R}^d$ with $|a|,|b|\le T$, it follows from Theorems \ref{th.e.7.2} and \ref{bu.th.6.5} that
\begin{align*}
	\mathbb{E}_\mu\Big[\big|Y(t,a)-Y(s,b)\big|^{2 N}\Big]
	&\le 2^{2 N}\Big\{\mathbb{E}_\mu\Big[\big|Y(t, a)-Y(t, b)\big|^{2 N}\Big]+\mathbb{E}_\mu\Big[\big|Y(t, b)-Y(s, b)\big|^{2 N}\Big]\Big\}\\
	&\le 
	\begin{cases}
	C_{N, T, \alpha, \gamma_1}|a-b|^{2 N \gamma_1}+C_{N, \varepsilon_0, T, \alpha, \gamma_1}|t-s|^{N \gamma_1}, & \text { if } d=1; \\ 
	C_{N, \varepsilon_0, T, \alpha, \gamma}|a-b|^{2 N \gamma}+C_{N, \varepsilon_0, T, \alpha, \tilde{\delta}}|t-s|^{N \tilde{\delta}}, & \text { if } d=2; \\ 
	C_{N, \varepsilon_0, T, \alpha, \gamma}|a-b|^{2 N \gamma}+C_{N, \varepsilon_0, T, \alpha, \hat{\delta}}|t-s|^{N \hat{\delta}}, & \text { if } d=3 ,
	\end{cases}
\end{align*}
where $\gamma_1\in(0,1)$, $\tilde{\delta}\in (0,1/2)$ and $\hat{\delta}\in (0,1/6)$. Hence by Kolmogorov's continuity criterion (see, e.g., Corollary 1.2 of Walsh \cite{W86}), we obtain that $Y(t,x)$ admits a joint continuous version for $x\in \mathbb{R}^d$ with $|x|\le T$ and $t\in[\varepsilon_0, T]$. Furthermore,
\begin{equation*}
			\left\{
			\begin{aligned}
				& \text{when}\ d=1,\ \text{the H\"older exponent of}\ Y(t,x)\ 
			\text{is arbitrarily close to}\ 1\
			\text{in}\ x\ \text{and}\ 1/2\ \text{in}\ t; \\
				& \text{when}\ d=2,\ \text{the H\"older exponent of}\ Y(t,x)\ 
			\text{is arbitrarily close to}\ \gamma\
			\text{in}\ x\ \text{and}\ 1/4\ \text{in}\ t;\\
			& \text{when}\ d=3,\ \text{the H\"older exponent of}\ Y(t,x)\ 
			\text{is arbitrarily close to}\ \gamma\
			\text{in}\ x\ \text{and}\ 1/12\ \text{in}\ t.
			\end{aligned}
			\right.
		\end{equation*}
Since $T, \varepsilon_0 > 0$ are arbitrary, the proof of Theorem \ref{t1} is complete.
$\hfill\blacksquare$

	\clearpage

	\bibliographystyle{plain}
	\def\cprime{$'$}

\appendix	
\section{The singularities in high dimensions}\label{appA}
In this section, we give the proofs of Theorems \ref{th.s.} and \ref{th.s..} for the singularities of $X$ and $Y$ with respect to Lebesgue measure in $d\ge 2$ and $d\ge 4$, respectively. To this end, we shall compare the Laplace transforms of $X$ and $Y$ with the super-Brownian motion and its occupation time process.

Fix $\mu\in M_F(\mathbb{R}^d)$. Let $\hat{X}=(\hat{X}_t, t\geq 0)$ be a super-Brownian motion starting from $\mu$. By slightly abusing the notation, we denote the law of $\hat{X}$ by $\mathbb{P}_{\mu}$ and its associated expectation by $\mathbb{E}_{\mu}$. According to the Laplace transform of the super-Brownian motion $\hat{X}$ (see, e.g., \cite[(1.1a) and (1.1b)]{DIP89}), it is easy to check that \cite[Theorem 2.18]{MX07} implies for any nonnegative bounded measurable functions $f$ and $\phi$,
\begin{align}\label{*a.1}
	\mathbb{E}_{\mu}\Big[ e^{-X_t(f)-Y_t(\phi) }\Big] \geq \mathbb{E}_{\mu}\Big[ e^{-\hat{X}_t(f)-\hat{Y}_t(\phi) }\Big],\quad\forall\ t\geq0,
\end{align}
where $\hat{Y}_t:=\int_0^t \hat{X}_s ds$ is the occupation time of the super-Brownian motion $\hat{X}$.
For any $x\in\mathbb{R}^d$ and $\varepsilon>0$, define
\begin{align*}
	B(x,\varepsilon):=\{y\in\mathbb{R}^d: |y-x|<\varepsilon\}.
\end{align*}

\noindent\textbf{Proof of Theorem \ref{th.s.}.} 
Let $d\ge 2$. For any $\lambda>0$, we set $f(y)=\lambda\mathbf{1}_{B(x,\varepsilon)}(y)$ and $\phi\equiv0$ in \eqref{*a.1} to get that for any $t\geq0$,
\begin{align*}
	\mathbb{E}_{\mu}\Big[ e^{-\lambda X_t(B(x,\varepsilon))}\Big] \geq \mathbb{E}_{\mu}\Big[ e^{-\lambda \hat{X}_t(B(x,\varepsilon))}\Big] .
\end{align*}
Letting $\lambda\rightarrow\infty$ in both sides of the above inequality implies
\begin{align*}
	\mathbb{P}_{\mu}\big( X_t(B(x,\varepsilon))>0 \big)	\leq \mathbb{P}_{\mu}\big( \hat{X}_t(B(x,\varepsilon))>0 \big).
\end{align*}
It then follows by the continuity of probabilities that
\begin{align*}
	\mathbb{P}_{\mu}\big( X_t(B(x,\varepsilon))>0,\ \forall\ \varepsilon>0\big)
	&= \lim_{\varepsilon\downarrow0}\mathbb{P}_{\mu}\big( X_t(B(x,\varepsilon))>0 \big) \nn\\
	&\leq \lim_{\varepsilon\downarrow0}\mathbb{P}_{\mu}\big( \hat{X}_t(B(x,\varepsilon))>0 \big) \nn\\
	&=\mathbb{P}_{\mu}\big( \hat{X}_t(B(x,\varepsilon))>0,\ \forall\ \varepsilon>0\big).
\end{align*}
Let $S(\nu)$ denote the closed support of a measure $\nu$. Note that for any measure $\nu$ and $x\in\mathbb{R}^d$,
\begin{align}\label{snv}
	\{ x\in S(\nu)\}=\{\nu(B(x,\varepsilon))>0,\ \forall\ \varepsilon>0\}.	
\end{align}
Therefore, 
\begin{align*}
	\mathbb{P}_{\mu}	\big(x\in S(X_t) \big) 
	&= \mathbb{P}_{\mu}\big( X_t(B(x,\varepsilon))>0,\ \forall\ \varepsilon>0\big) \nn\\
	&\leq \mathbb{P}_{\mu}\big( \hat{X}_t(B(x,\varepsilon))>0,\ \forall\ \varepsilon>0\big)= \mathbb{P}_{\mu}	\big(x\in S(\hat{X}_t) \big).
\end{align*}
Recall that when $d\geq 2$, with $\mathbb{P}_{\mu}$-probability one, the super-Brownian motion $\hat{X}_t$ is singular for any $t>0$ (see, e.g., \cite[Corollary D]{P88}). Thus we have $\mathbb{P}_{\mu}(x\in S(\hat{X}_t))=0$ for Lebesgue-almost every $x\in\mathbb{R}^d$. Hence,
\begin{align*}
	\mathbb{P}_{\mu} (x\in S({X}_t)) =0,\quad\text{for Lebesgue-a.e. $x\in\mathbb{R}^d$}.
\end{align*}
By Fubini's theorem, we get with $\mathbb{P}_{\mu}$-probability one, $x\notin S(X_t)$ for almost every $x\in\mathbb{R}^d$, which means $X_t$ is supported on a Lebesgue null set, i.e., $X_t$ is a singular measure. So we obtain 
\begin{align*}
	\mathbb{P}_{\mu} (X_t \text{ is singular})=1,\quad\forall\ t>0.	
\end{align*}
By another application of Fubini's theorem with $t>0$ we get
\begin{align*}
	\mathbb{P}_{\mu} (X_t \text{ is singular for Lebesgue-a.e. }t>0)=1.	
\end{align*}
The proof now is complete.
$\hfill\blacksquare$\\

\noindent\textbf{Proof of Theorem \ref{th.s..}.} 
The proof is similar to that of Theorem \ref{th.s.}. Let $d\ge 4$. For any $\lambda>0$, we set $\phi(y)=\lambda\mathbf{1}_{B(x,\varepsilon)}(y)$ and $f\equiv0$ in \eqref{*a.1} to get that for any $t\geq0$,
\begin{align*}
	\mathbb{E}_{\mu}\Big[ e^{-\lambda Y_t(B(x,\varepsilon))}\Big] \geq \mathbb{E}_{\mu}\Big[ e^{-\lambda \hat{Y}_t(B(x,\varepsilon))}\Big] .
\end{align*}
Letting $\lambda\rightarrow\infty$ and using the continuity of probabilities, we have
\begin{align*}
	\mathbb{P}_{\mu}\big( Y_t(B(x,\varepsilon))>0,\ \forall\ \varepsilon>0\big)
	\le \mathbb{P}_{\mu}\big( \hat{Y}_t(B(x,\varepsilon))>0,\ \forall\ \varepsilon>0\big).
\end{align*}
Recall that $S(\nu)$ is the closed support of a measure $\nu$. Then by \eqref{snv} we obtain that
\begin{align*}
	\mathbb{P}_{\mu}	\big(x\in S(Y_t) \big) 
	&= \mathbb{P}_{\mu}\big( Y_t(B(x,\varepsilon))>0,\ \forall\ \varepsilon>0\big) \nn\\
	&\leq \mathbb{P}_{\mu}\big( \hat{Y}_t(B(x,\varepsilon))>0,\ \forall\ \varepsilon>0\big)= \mathbb{P}_{\mu}	\big(x\in S(\hat{Y}_t) \big).
\end{align*}
Recall that when $d\geq 4$, with $\mathbb{P}_{\mu}$-probability one, the occupation time $\hat{Y}_t$ of the super-Brownian motion is singular for any $t\geq 0$ (see, e.g. \cite{DIP89} for $d>4$ and \cite{LG99} for $d=4$,). Thus we have $\mathbb{P}_{\mu}(x\in S(\hat{Y}_t))=0$ for Lebesgue-almost every $x\in\mathbb{R}^d$. Hence,
\begin{align*}
	\mathbb{P}_{\mu} (x\in S({Y}_t)) =0,\quad\text{for Lebesgue-a.e. $x\in\mathbb{R}^d$}.
\end{align*}
By Fubini's theorem, we get with $\mathbb{P}_{\mu}$-probability one, $x\notin S(Y_t)$ for almost every $x\in\mathbb{R}^d$, which means $Y_t$ is supported on a Lebesgue null set, i.e., $Y_t$ is a singular measure. So we obtain 
\begin{align*}
	\mathbb{P}_{\mu} (Y_t \text{ is singular})=1,\quad\forall\ t\geq0.	
\end{align*}
By another application of Fubini's theorem with $t\geq0$ we get
\begin{align*}
	\mathbb{P}_{\mu} (Y_t \text{ is singular for Lebesgue-a.e. }t\geq0)=1.	
\end{align*}
The proof now is complete.
$\hfill\blacksquare$

	\section{Moment bounds}\label{sec4}
	
\subsection{Bounds on $V^{\phi}_n(t,x)$}\label{appb1}
	
	Let $d\ge 1$. Throughout this subsection, we fix $\phi\in C_b^{2,+}(\R^d)$ and give some bounds on $V^{\phi}_n(t,x)$ defined as in \eqref{e2.1}.
	
	\begin{lemma}\label{lem5.1}
		Let $Q_t$ be given as in (\ref{q1}). For any fixed $T>0$ and $N\geq 1$, there is some constant $C_{N,T}>0$ {\rm(}independent of $\phi${\rm)} such that for any $0\leq t\leq T$ and $x_1, \cdots, x_N \in \R^d$, 
		\begin{align}\label{e.5.1}
			\E\big[V^{\phi}_1(t, x_1)\cdots V^{\phi}_1(t, x_N)\big]\leq C_{N,T} \prod_{i=1}^N Q_{t} \phi(x_i).
		\end{align}
	\end{lemma}
	
	\begin{proof}
		We shall give the proof by induction in $N\geq1$. For $N=1$, by \eqref{e10.5}, we see
		\begin{align*}
			\E\big[V^{\phi}_1(t, x_1)\big] = Q_{t} \phi(x_1),\quad x_1\in\R^d,
		\end{align*}
		which implies \eqref{e.5.1} holds. Now supposing \eqref{e.5.1} is satisfied when $N$ is replaced by $N-1$, we prove it also holds for $N$. Note that \eqref{e10.5} can be rewritten as
		\begin{align}\label{e.5.2}
			dV^{\phi}_1(t,x)= \frac{\Delta}{2} V^{\phi}_1(t,x) dt+V^{\phi}_1(t,x) {W}(dt,x)+\phi(x) dt.
		\end{align}
		For any $x_1,\cdots,x_N\in\R^d$, by It\^o's formula, it is easy to check that
		\begin{align}\label{e.5.3}
			d\bigg(\prod_{i=1}^N V^{\phi}_1(t, x_i)\bigg) = &\sum_{i=1}^{N}  \Big(  \prod_{\substack{1\leq j \leq N \\ j \neq i}} V^{\phi}_1(t, x_j) \Big) d V^{\phi}_1(t, x_i) \nonumber\\
			&+ \sum_{1\leq i< j\leq N}  \bigg(\prod_{\substack{1\leq l \leq N \\ l \neq i, j}} V^{\phi}_1 (t, x_l ) \bigg) d\langle V^{\phi}_1(\cdot, x_i), V^{\phi}_1(\cdot, x_j)\rangle_t .
		\end{align}
		In view of \eqref{e10.5}, a simple calculation shows that
		\begin{align*}
			d\langle V^{\phi}_1(\cdot, x_i), V^{\phi}_1(\cdot, x_j)\rangle_t = V^{\phi}_1(t, x_i) V^{\phi}_1(t, x_j) g(x_i, x_j) d t.
		\end{align*}
		This together with \eqref{e.5.2} and \eqref{e.5.3} yields  
		\begin{align*}
			d\bigg(\prod_{i=1}^N V^{\phi}_1(t, x_i)\bigg) 
			=&\ \sum_{i=1}^{N}  \Big(  \prod_{\substack{1 \leq j \leq N \\ j \neq i}} V^{\phi}_1(t, x_j) \Big) \Big[ \frac{\Delta}{2} V^{\phi}_1(t,x_i) dt+V^{\phi}_1(t,x_i) {W}(dt,x_i)+\phi(x_i) dt \Big] \nonumber\\
			&\ + \sum_{1\leq i< j\leq N}  \bigg(\prod_{\substack{1\leq l \leq N \\ l \neq i, j}} V^{\phi}_1 (t, x_l )  V^{\phi}_1(t, x_i) V^{\phi}_1(t, x_j) \bigg) g(x_i, x_j) d t \nonumber\\
			=&\ \frac{\Delta_{(x_1, \cdots, x_N)}}{2}\Big(\prod_{i=1}^N V^{\phi}_1(t, x_i)\Big) d t + \Big(\prod_{i=1}^N V^{\phi}_1(t, x_i) \Big)\Big( \sum_{j=1}^N W(d t, x_j) \Big)  \nonumber\\
			&\ + \sum_{i=1}^N \phi(x_i)  \Big(\prod_{\substack{1 \leq j \leq N \\ j \neq i}} V^{\phi}_1(t, x_j) \Big) d t +  \Big(\prod_{i=1}^N V^{\phi}_1(t, x_i) \Big) \Big(\sum_{1\leq i <j \leq N} g(x_i, x_j) \Big)d t.
		\end{align*}
		Therefore, 
		\begin{align*}
			d \E\bigg[\prod_{i=1}^N V^{\phi}_1(t, x_i)\bigg] = &\ \frac{\Delta_{(x_1, \cdots, x_N)}}{2} \E\bigg[\prod_{i=1}^N V^{\phi}_1(t, x_i)\bigg] d t + \sum_{i=1}^N \phi(x_i)  \E \bigg[\prod_{\substack{1 \leq j \leq N \\ j \neq i}} V^{\phi}_1(t, x_j) \bigg] d t \nonumber\\
			&+\E\bigg[\prod_{i=1}^N V^{\phi}_1(t, x_i)\bigg]  \Big(\sum_{1\leq i <j \leq N} g(x_i, x_j) \Big)d t.
		\end{align*}
		Since $\E[\prod_{i=1}^N V^{\phi}_1(0, x_i)]=0$, we use Feynman-Kac's formula to get
		\begin{align}\label{e.5.4}
			\E\bigg[&\prod_{i=1}^N V^{\phi}_1(t, x_i)\bigg]\\
			&= \Pi_{(x_1,\cdots,x_N)} \bigg\{ \int_{0}^t \sum_{i=1}^N \phi(B_s^i) \E\bigg[\prod_{\substack{1\leq j \leq N \\ j \neq i}} V^{\phi}_1(t-s, B_s^j) \bigg] e^{\int_0^s \Big(\sum_{1\leq i <j \leq N} g(B_r^i, B_r^j) \Big)d r} ds \bigg\} \nonumber,
		\end{align}
		where $B^1,\cdots,B^N$ are independent $d$-dimensional Brownian motions starting respectively from $x_1,\cdots,x_N\in\R^d$ under $\Pi_{(x_1,\cdots,x_N)}$. By the induction hypothesis, we see that for any $0\leq t\leq T$,
		\begin{align*}
			\E\bigg[\prod_{\substack{1\leq j \leq N \\ j\neq i}} V^{\phi}_1(t-s, x_j) \bigg] \leq  C_{N,T} \prod_{\substack{1\leq j \leq N \\ j\neq i}} Q_{t-s} \phi(x_j).
		\end{align*}
		For any $0\leq s \leq t \leq T$, we have
		\begin{align*}
			\int_0^s \Big(\sum_{1\leq i <j \leq N} g(B_r^i, B_r^j) \Big)d r \leq \frac{N(N-1)}{2} T \|g\|_{\infty}.
		\end{align*}
		Since $\phi\geq0$, we take the two preceding inequalities into \eqref{e.5.4} to obtain 
		\begin{align*}
			\E\bigg[\prod_{i=1}^N V^{\phi}_1(t, x_i)\bigg] 
			&\leq e^{ \frac{N(N-1)}{2} T \|g\|_{\infty} } C_{N,T}  \Pi_{(x_1,\cdots,x_N)} \bigg[ \int_{0}^t \sum_{i=1}^N \phi(B_s^i) \prod_{\substack{1\leq j \leq N \\ j\neq i}} Q_{t-s} \phi(B^j_s) ds \bigg] \nonumber\\
			&= e^{ \frac{N(N-1)}{2} T \|g\|_{\infty} } C_{N,T} \sum_{i=1}^N \int_{0}^t P_s\phi(x_i) ds \Big[\prod_{\substack{1\leq j \leq N \\ j\neq i}} \int_{s}^t  P_r \phi(x_j) dr \Big] \nonumber\\
			&\leq  C_{N,T} \prod_{i=1}^N Q_t\phi(x_i),
		\end{align*}
		where 
		the equality in the above comes from the independence of Brownian motions $B^i,i=1,\cdots,N$. The proof is now complete.
	\end{proof}

	\begin{lemma}\label{al1}
		Let $Q_t$ be given as in (\ref{q1}). Given $T>0$ and $N\geq 1$, for any $0\leq t_1, \cdots, t_N\leq T$ and $x_1, \cdots, x_N \in \R^d$, we have  
		\begin{align*}
			\E\big[V^{\phi}_1(t_1, x_1)\cdots V^{\phi}_1(t_N, x_N)\big]\leq {C}_{N,T} \prod_{i=1}^N {Q}_{t_i} \phi(x_i),
		\end{align*}
		where ${C}_{N,T}>0$ is a constant given as in Lemma \ref{lem5.1}.
	\end{lemma}
	
	\begin{proof}
		Fix $T>0$ and $N\geq 1$. By Lemma \ref{lem5.1}, we see that for $0\leq t \leq T$ and $x\in\R^d$,
		\begin{align*} 
			\E\big[V^{\phi}_1(t, x)^{N}\big] \leq C_{N,T} {Q}_{t} \phi(x)^{N}.
		\end{align*}
		Thus, for any $0\leq t_1, \cdots, t_N\leq T$ and $x_1$, $\cdots$, $x_N \in \R^d$, we apply a generalized H\" older's inequality 
		to see that 
		\begin{align*}
			\E\big[V^{\phi}_1(t_1, x_1)\cdots V^{\phi}_1(t_N, x_N)\big]
			\leq  \prod_{i=1}^N \Big(\E\Big[V^{\phi}_1(t_i, x_i)^{N}\Big]\Big)^{1/N} \leq {C}_{N,T} \prod_{i=1}^N  {Q}_{t_i} \phi(x_i).
		\end{align*}
		So the conclusion follows.
	\end{proof}

	For any $t\geq 0$ and $\phi\in C_b^{2,+}(\R^d)$, define 
	\begin{align}\label{g}
		H(\phi, t):=\sup_{y\in \R^d} \int_0^t P_s \phi(y) ds<\infty,
	\end{align}
	and let
	\begin{align}\label{g.1}
		G(\phi,t):=\max\{H(\phi,t),1\}.	
	\end{align}
	Note that
	\begin{align*}
		H(\phi,t)\leq G(\phi,t) \leq 1+ H(\phi,t) \leq 1+t\|\phi\|_{\infty}<\infty,	
	\end{align*}
	and $G(\phi,t)$ is non-decreasing on $t\geq 0$.
	Then the above lemma implies that
	\begin{align*}
		\E\big[V^{\phi}_1(t, x)^N\big]\leq C_{N,T} G(\phi, T)^N,\quad\forall\ 0\le t\le T,\ x\in\mathbb{R}^d.
	\end{align*}
	The following comparison lemma plays an important role in giving the bound for $\mathbb{E}[V^{\phi}_n(t,x)^{2N}]$.
	
	\begin{lemma}\label{l4.3}(\cite[Lemma 6.1]{FHX24})
		Let $d\geq 1$ and $T>0$. For any two continuous functions $F(t,x)$ and $G(t,x)$ defined on $[0, T] \times \R^d$, if there exist some function $\{\alpha(t,x): 0\leq t \leq T, x\in \R^d\}$ and constant $\beta_T>0$ such that for all $0\leq t\leq T$ and $x\in \R^d$,
		\begin{align*}
			G(t,x)\leq \alpha(t,x)+\beta_T\int_0^t ds \int p_{t-s}(x,y) G(s,y)dy,
		\end{align*}
		and
		\begin{align*}
			F(t,x)= \alpha(t,x)+\beta_T\int_0^t ds \int p_{t-s}(x,y) F(s,y)dy,
		\end{align*}
		then 
		\begin{align*}
			G(t,x)\leq F(t,x), \quad \forall\  0\leq t\leq T,\ x\in \R^d.
		\end{align*}
	\end{lemma}
	
	\begin{lemma}\label{lem5.3}
		Let $G(\phi,t)$ be defined by (\ref{g}). For every $n\geq 2$, $N\geq1$ and $T>0$. There is some constant $C_{n,N,T}>0$ {\rm(}independent of $\phi${\rm)} such that for any $0\leq t \leq T$ and $x \in \R^d$, 
		\begin{align}\label{e.5.5}
			\E\big[V^{\phi}_n(t, x)^{2N}\big]\leq C_{n,N,T} G(\phi, T)^{N 2^{n}}.
		\end{align}
	\end{lemma}
	
	\begin{proof}
		Fix $N\geq 1$, $T>0$, $0\leq t\leq T$ and $x\in\R^d$. We shall prove \eqref{e.5.5} holds by induction in $n\geq 2$. First, we deal with $n=2$. Recall from \eqref{e6.11} that
		\begin{align*}
			V^{\phi}_2(t, x)=  \int_0^t ds \int p_{t-s} (x,z)V^{\phi}_1(s, z)^2 dz +\int_0^t  \int p_{t-s}(x,z)  V^{\phi}_{2}(s,z) W(ds,z) dz.
		\end{align*}
		Using the above argument, we obtain
		\begin{align}\label{e5.12}
			\E\big[V^{\phi}_2(t,x)^{2N}\big] \leq &\ 2^{2N}  \E\Big[\Big(\int_0^t ds \int p_{t-s} (x,z)V^{\phi}_1(s,z)^2 dz\Big)^{2N} \Big]\nn\\
			&+2^{2N} \E\Big[ \Big(\int_0^t  \int p_{t-s}(x,z)  V^{\phi}_{2}(s,z) W(ds,z) dz\Big)^{2N}\Big].
		\end{align}
		For the first expectation on the right-hand side above, we get
		\begin{align}\label{e5.31}
			J_1:= \E\Big[&\Big(\int_0^t ds \int p_{t-s} (x,z)V^{\phi}_1(s,z)^2 dz\Big)^{2N} \Big]\nn\\
			=\E\Big[&\int_0^t  ds_1 \int  p_{t-s_1} (x,z_1) d z_1 \cdots \int_0^t  ds_{2N-1} \int  p_{t-s_{2N-1}} (x,z_{2N-1})   d z_{2N-1} \nn\\
			&\int_0^t  ds_{2N} \int  p_{t-s_{2N}} (x,z_{2N})  \prod_{i=1}^{2N}  V^{\phi}_1(s_i, z_i)^2  d z_{2N}  \Big].
		\end{align}
		By a generalized H\" older's inequality and \eqref{e.5.1}, we see
		\begin{align}
			\E\Big[ \prod_{i=1}^{2N}  V^{\phi}_1(s_i, z_i)^2 \Big] 
			&\leq 	\prod_{i=1}^{2N} \Big( \E\Big[   V^{\phi}_1(s_i, z_i)^{4N} \Big] \Big)^{1/(2N)} \leq C_{N,T} \prod_{i=1}^{2N} Q_{s_i}\phi(z_i)^2\nn.
		\end{align}
		By using $\phi\geq0$ and $0\leq s_i\leq t\leq T$, we see $Q_{s_i}\phi(z_i)\leq G(\phi, T)$. Hence, 
		\begin{align*}
			\E\Big[ \prod_{i=1}^{2N}  V^{\phi}_1(s_i, z_i)^2 \Big] \leq C_{N,T} G(\phi,T)^{4N}.
		\end{align*}
		We then apply the above inequality to \eqref{e5.31} to obtain
		\begin{align*}
			J_1&\leq C_{N,T} G(\phi,T)^{4N} \prod_{i=1}^{2N} \Big( \int_0^t  ds_i \int p_{t-s_i} (x,z_i)  d z_i  \Big)\\
			&=C_{N,T} G(\phi,T)^{4N} \Big(\int_0^t   ds  \int  p_{t-s } (x,z)  d z \Big)^{2N}.
		\end{align*}
		Note that
		\begin{align*}
			&\int_0^t   ds  \int  p_{t-s } (x,z)  d z \leq t \leq T.
		\end{align*}
		Then we conclude that
		\begin{align}\label{e5.32}
			J_1&\leq  C_{N,T} T^{2N} G(\phi, T)^{4N}.
		\end{align}
		Turning to the second expectation in \eqref{e5.12}, we apply Burkholder-Davis-Gundy's inequality to see that 
		\begin{align}\label{e5.21}
			J_2:=&\ \E\Big[ \Big(\int_0^t  \int p_{t-s}(x,z)  V^{\phi}_{2}(s,z) W(ds,z) dz\Big)^{2N}\Big]\nn\\
			\leq&\ C_N \E\Big[\Big(\int_0^t ds \int\!\!\int p_{t-s} (x,z)  p_{t-s} (x,w)V^{\phi}_{2}(s,z)  V^{\phi}_{2}(s,w)  g(z,w)dz dw\Big)^{N} \Big]\nn\\
			\leq&\  C_N \|g\|_{\infty}^{N} \E\Big[\Big(\int_0^t   \Big[\int p_{t-s} (x,z) V^{\phi}_{2}(s,z)   dz \Big]^2 ds\Big)^{N} \Big] \nn\\
			\leq&\ C_N \|g\|_{\infty}^{N} \E\Big[\Big(\int_0^t ds  \int p_{t-s} (x,z) V^{\phi}_{2}(s,z)^2  dz \Big)^{N} \Big],
		\end{align}
		where the constant $C_N>0$ and the last inequality follows from Cauchy-Schwarz's inequality. 		
		Applying H\"older's inequality with $p=N/(N-1)$ and $q=N$ yields 
		\begin{align}\label{e5.34}
			&\Big(\int_0^t ds  \int p_{t-s} (x,z) V^{\phi}_{2}(s,z)^2  dz\Big)^N\leq t^{N-1} \int_0^t ds  \int p_{t-s} (x,z) V^{\phi}_{2}(s,z)^{2N}  dz ,
		\end{align}
		which implies
		\begin{align}\label{e5.22}
			J_2  &\leq  C_N \|g\|_{\infty}^{N} T^{N-1} \E\Big[\int_0^t ds  \int p_{t-s} (x,z) V^{\phi}_{2}(s,z)^{2N}  dz \Big]\nn\\
			&= C_N \|g\|_{\infty}^{N} T^{N-1} \int_0^t ds  \int p_{t-s} (x,z) \E\big[V^{\phi}_{2}(s,z)^{2N}\big]  dz .
		\end{align}
		By combining the bounds for $J_1$ from \eqref{e5.32} and $J_2$ from \eqref{e5.22}, it is easily seen that \eqref{e5.12} becomes
		\begin{align}\label{e5.14}
			\E\big[V^{\phi}_2(t,x)^{2N}\big] \leq  C_{N,T}G(\phi, T)^{4N}   +C_{N,T} \int_0^t ds  \int p_{t-s} (x,z) \E\big[V^{\phi}_{2}(s,z)^{2N}\big]  dz,
		\end{align}
		Define $F_N(t,x)$ to be the solution of
		\begin{align}\label{e5.15}
			F_N(t,x) = C_{N,T}G(\phi, T)^{4N} +C_{N,T} \int_0^t ds  \int p_{t-s} (x,z) F_N(s,z)  dz.
		\end{align}
		In view of \eqref{e5.14} and \eqref{e5.15}, Lemma \ref{l4.3} implies
		\begin{align}\label{e.5.25}
			\E\big[V^{\phi}_2(t,x)^{2N}\big] \leq F_N(t,x), \quad \forall\  0\leq t\leq T,\ x\in \R^d.
		\end{align} 
		It suffices to find the bound for $F_N(t,x)$. By a simple calculation, we see that
		\begin{align*}
			\frac{\partial}{\partial t} F_N(t,x) =   \frac{\Delta}{2}  F_N(t,x)+C_{N,T} F_N(t,x).
		\end{align*}
		We then apply Feynman-Kac's formula with $F_N(0,x)=C_{N,T}G(\phi, T)^{4N}$ to get
		\begin{align*}
			F_N(t,x) = C_{N,T}G(\phi, T)^{4N} e^{C_{N,T} t} \leq C_{N,T} G(\phi, T)^{4N}.
		\end{align*}
		In view of \eqref{e.5.25}, we obtain
		\begin{align}\label{e.5.26}
			\E[V^{\phi}_2(t,x)^{2N}]\leq    C_{N,T}G(\phi, T)^{4N}, \quad \forall\ 0\leq t\leq T,\ x\in \R^d
		\end{align}
		as required. Assuming that \eqref{e.5.5} holds for all $2\leq k\leq n-1$ with some $n\geq 3$, we will prove \eqref{e.5.5} holds for $n$. To see this, by \eqref{e6.11}, note that 
		\begin{align}\label{e4.24}
			\E\big[V^{\phi}_n(t,x)^{2N}\big] \leq &\ C_{n,N} \E\Big[ \Big(\int_0^t  \int p_{t-s}(x,z)  V^{\phi}_n(s,z) W(ds,z) dz\Big)^{2N}\Big]\nn\\
			&+C_{n,N} \sum_{k=1}^{n-1} \E\Big[\Big(\int_0^t ds \int p_{t-s} (x,z)V^{\phi}_k(s,z) V^{\phi}_{n-k}(s,z) dz\Big)^{2N} \Big],
		\end{align}
		where $C_{n,N}>0$ is a constant depending on $n$ and $N$. For the first expectation above, by replacing $V^{\phi}_2(s,z)$ with $V^{\phi}_n(s,z)$ and repeating the arguments from \eqref{e5.21} to \eqref{e5.22}, we get
		\begin{align}\label{e4.52}
			\E\Big[ \Big(\int_0^t  \int p_{t-s}(x,z)  &V^{\phi}_n(s,z) W(ds,z) dz\Big)^{2N}\Big] \nn\\
			&\leq C_N \|g\|_{\infty}^N T^{N-1} \int_0^t ds  \int p_{t-s} (x,z)  \E\big[V^{\phi}_n(s,z)^{2N}\big] dz,
		\end{align}
		where the constant $C_N>0$. For the second expectation in \eqref{e4.24}, we use Cauchy-Schwarz's inequality twice to get that for each $1\leq k\leq n-1$,
		\begin{align}\label{e4.51}
			&J_2(k):=\E\Big[\Big(\int_0^t ds \int p_{t-s} (x,z)V^{\phi}_k(s,z) V^{\phi}_{n-k}(s,z) dz\Big)^{2N} \Big] \\
			\leq &\   \E\Big[ \Big(\int_0^t ds \int p_{t-s} (x,z)V^{\phi}_k(s,z)^2  dz\Big)^{N}\Big(\int_0^t ds \int p_{t-s} (x,z)V^{\phi}_{n-k}(s,z)^2  dz\Big)^{N} \Big]\nn\\
			\leq &\  \bigg\{  \E\Big[\Big(\int_0^t ds \int p_{t-s} (x,z)V^{\phi}_k(s,z)^2  dz\Big)^{2N} \Big] \bigg\}^{1/2} \bigg\{\E \Big[\Big(\int_0^t ds \int p_{t-s} (x,z)V^{\phi}_{n-k}(s,z)^2  dz\Big)^{2N}\Big] \bigg\}^{1/2}.\nn
		\end{align}
		If $k=1$, then in view of \eqref{e5.31} and \eqref{e5.32}, we get
		\begin{align*}
			\E\Big[\Big(\int_0^t ds \int p_{t-s} (x,z)V^{\phi}_1(s,z)^2  dz\Big)^{2N} \Big] \leq C_{N,T} T^{2N} G(\phi, T)^{4N}.
		\end{align*}
		If $2\leq k\leq n-1$, then by the induction hypothesis we have 
		\begin{align*} 
			\E\big[V^{\phi}_k(t,x)^{4N}\big]\leq C_{k,N,T} G(\phi, T)^{2^{k+1} N}, \quad \forall\ 0\leq t\leq T,\ x\in \R^d,
		\end{align*}
		where $C_{k,N,T}>0$. Similar to \eqref{e5.34}, we may apply H\"older's inequality with $p=2N/(2N-1)$ and $q=2N$ to get
		\begin{align*}
			&\int_0^t ds  \int p_{t-s} (x,z) V^{\phi}_k(s,z)^2  dz\leq t^{(2N-1)/(2N)}   \Big(\int_0^t ds  \int p_{t-s} (x,z) V^{\phi}_k(s,z)^{4N}  dz\Big)^{1/(2N)}.
		\end{align*}
		Thus, the two preceding inequalities imply
		\begin{align*}
			\E\Big[\Big(\int_0^t ds  \int p_{t-s} (x,z) V^{\phi}_k(s,z)^2  dz \Big)^{2N} \Big]
			&\leq  t^{2N-1} \int_0^t ds  \int p_{t-s} (x,z)  \E\big[V^{\phi}_k(s,z)^{4N}\big]  dz \\ 
			&\leq C_{k,N,T} T^{2N} G(\phi, T)^{2^{k+1} N}.
		\end{align*}
		This combined with \eqref{e4.51} yields that for any $1\leq k\leq n-1$, there exists a constant $C_{k,n,N,T}>0$ such that
		\begin{align*}
			J_2(k)  \leq   C_{k,n,N,T} G(\phi, T)^{2^{n-1}  N(2^{-(n-1-k)}+2^{-(k-1)})} \leq C_{k,n,N,T} G(\phi, T)^{2^{n}N},
		\end{align*}
		where the last inequality follows from $G(\phi,T)\geq 1$ and $1\leq k \leq n-1$. By taking the above inequality and \eqref{e4.52} into \eqref{e4.24}, we obtain
		\begin{align*} 
			\E\big[V^{\phi}_n(t,x)^{2N}\big] \leq&\ C_{n,N}  G(\phi, T)^{2^{n}  N}  \sum_{k=1}^{n-1} C_{k,n,N,T} \\
			&+ C_{n,N} C_N \|g\|_{\infty}^N T^{N-1} \int_0^t ds  \int p_{t-s} (x,z)  \E\big[V^{\phi}_n(s,z)^{2N}\big] dz\\
			\leq&\ C_{n,N,T} G(\phi, T)^{2^{n}  N} +C_{n,N,T}\int_0^t ds  \int p_{t-s} (x,z)  \E\big[V^{\phi}_n(s,z)^{2N}\big] dz,
		\end{align*}
		where the constant $C_{n,N,T}>0$. The above gives the same inequality as in \eqref{e5.14}, and hence by \eqref{e.5.26} we conclude that there exists a constant $C_{n,N,T}>0$ such that
		\begin{align*}
			\E[V^{\phi}_n(t,x)^{2N}]\leq   C_{n,N,T} G(\phi, T)^{2^{n}  N} 
		\end{align*}
		as required. The proof is complete by induction.
	\end{proof}
	
	\subsection{Bounds on $Y_t$}\label{appb2}
	
Recall that the sequence $\{\phi_{a,\gamma,r}\}\subset C_b^{2,+}(\mathbb{R}^d)$ is given by \eqref{phir} such that $\text{as}\ r\rightarrow\infty$,
\begin{align*}
	\phi_{a,\gamma,r}\uparrow \phi_{a,\gamma},\quad 0<\gamma<d,\ a\in \mathbb{R}^d,
\end{align*}
where $\phi_{a,\gamma}$ is defined by \eqref{phig}.
We denote $V^{\phi}_{n}(t,x)$ and $L^{\phi}_n(t)$  with respect to the text function $\phi=\phi_{a,\gamma,r}$ by $V_{n,r}(t,x)$ and $L_{n,r}(t)$, respectively. The dependence of $a$ and $\gamma$ will be suppressed for notational ease.  It follows from Lemma \ref{al1} that for any $0\leq t_1, \cdots, t_N\leq T$ and $x_1, \cdots, x_N \in \R^d$, 
		\begin{align}\label{'5.50}
			\E\big[V_{1,r}(t_1, x_1)\cdots V_{1,r}(t_N, x_N)\big]\leq {C}_{N,T} \prod_{i=1}^N {Q}_{t_i} \phi_{a,\gamma,r}(x_i)\leq {C}_{N,T} \prod_{i=1}^N {Q}_{t_i} \phi_{a,\gamma}(x_i).
		\end{align}
	Recall from (\ref{e6.11}) that $V_{n,r}(t,x)$ satisfies  
	\begin{align}\label{'5.51}
		V_{n,r}(t,x)=  &\sum_{k=1}^{n-1} \binom{n-1}{k}   \int_0^t ds \int p_{t-s} (x,z)V_{n-k,r}(s, z)V_{k,r}(s, z) dz\nn\\
		&+\int_0^t  \int p_{t-s}(x,z)  V_{n,r}(s,z) W(ds,z) dz,\quad \forall \ n\geq 2.
	\end{align}
	Based on (\ref{'5.50}) and (\ref{'5.51}), for $d=2,3$, we shall give upper bounds of $\mathbb{E}[V_{n,r}(t,x)^{2N}]$ and further of the moments of $Y_{t}(\phi_{a,\gamma,r})$
	, where the bounds are all independent of $r$.
	
	\subsubsection*{Case 1: d=2}
	
	\begin{lemma}\label{'lem5.13}
		Let $d=2$. For every $n \geq 2,\ N \geq 1$ and $T>0$, there exists a constant $C_{n, N, T, \gamma}>0$ {\rm(}independent of $r${\rm)} such that for any $0 \leq t \leq T,\ x, a \in \mathbb{R}^2$ and $r \geq 1$,
		\begin{align}\label{'5.52}
			\mathbb{E}\big[V_{n,r}(t, x)^{2 N}\big] \leq C_{n, N, T, \gamma}.
		\end{align}
	\end{lemma}
	\begin{proof}
		It follows from Lemma \ref{lem5.3} that
		\begin{align*}
			\mathbb{E}\big[V_{n,r}(t, x)^{2 N}\big] \leq C_{n, N, T}G(\phi_{a,\gamma,r},t)^{N2^n},
		\end{align*}
		where $G(\phi_{a,\gamma,r},t)$ is defined by (\ref{g.1}) with $\phi=\phi_{a,\gamma,r}$. Since
		$0\leq \phi_{a,\gamma,r}\uparrow \phi_{a,\gamma}$, we get that for any $0\leq t\leq T$,
		\begin{align*}
			H(\phi_{a,\gamma, r},t)\leq \sup_{x\in\mathbb{R}^2} \int_0^T P_s\phi_{a,\gamma}(x)ds:= H(\phi_{a,\gamma},T).	
		\end{align*}
		By Lemma \ref{lem5.1'}, we see
		\begin{align*}
			H(\phi_{a,\gamma},T) = \sup_{x\in\mathbb{R}^2} \int_0^T ds \int \phi_{a,\gamma}(y) p_s(x,y) dy \leq C_{\gamma} T^{1-\frac{\gamma}{2}}.
		\end{align*}
		Hence,
		\begin{align}\label{'5.53}
			G(\phi_{a,\gamma},T):=\max\{H(\phi_{a,\gamma},T),1\}\le 1+ C_{\gamma}T^{1-\frac{\gamma}{2}},
		\end{align}
		which implies the desired result.
	\end{proof}
	
	\begin{lemma}\label{'lem5.14}
		Let $d=2$. For every $n \geq 1$ and $T>0$, there exists some constant $C_{n, T, \gamma}>0$ {\rm(}independent of $r${\rm)} such that for any $0 \leq t \leq T,\ a \in \mathbb{R}^2$ and $r \geq 1$,
		\begin{align*}
			\mathbb{E}_{\mu}\big[Y_{t}(\phi_{a,\gamma,r})^{n}\big] \leq C_{n, T, \gamma}.
		\end{align*}
	\end{lemma}
	\begin{proof}
		Recall from Proposition \ref{l2.1} with $\phi=\phi_{a,\gamma,r}$ that
		\begin{align}\label{''y1}
\mathbb{E}[Y_t(\phi_{a,\gamma,r})^n]=\mathbb{E}[L_{n,r}(t)].
\end{align}
Fix $T>0$ and let $0\le t\le T$, for any $n\ge 1,\ a \in \mathbb{R}^2$ and $r \geq 1$, we shall prove that
\begin{align}\label{'5.55}
\mathbb{E}_\mu[L_{n,r}(t)^N] \leq C_{n,N, T,\gamma},\quad \forall\  N\ge 1,
\end{align}
for some constant $C_{n,N, T,\gamma}>0$ independent of $r$. The conclusion follows immediately from 
(\ref{''y1}) and the above with $N=1$. 
When $n=1$, for any $N\ge 1$ we have
\begin{align}\label{'5.56}
\mathbb{E}[L_{1,r}(t)^N] =\mathbb{E}[\langle \mu,V_{1,r}(t,\cdot)\rangle^N]&=\int \mu(d x_1) \cdots \int \mathbb{E}\bigg[\prod_{i=1}^N V_{1,r}(t, x_i)\bigg] \mu(d x_n)\nn\\
	&\le C_{N, T}\bigg[\int Q_t \phi_{a,\gamma}(x) \mu(d x)\bigg]^N\nn\\
	&\le C_{N, T}[\mu(1)G(\phi_{a,\gamma},T)]^N \nn\\
	&\le C_{N, T}\big[\mu(1)(1+C_{\gamma}T^{1-\frac{\gamma}{2}} )\big]^N,
\end{align}
where the first and the last inequality follow by (\ref{'5.50}) and (\ref{'5.53}), respectively. Hence (\ref{'5.55}) holds for the case $n=1$. Assume that (\ref{'5.55}) holds for all $1\le k\le n-1$ with some $n\ge 2$. Then for the case $n$, using Lemma \ref{'lem5.13}, we repeat the proof in Lemma \ref{lem5.5} with $\phi$ replaced by $\phi_{a,\gamma,r}$ to get the desired result.
	\end{proof}
	
	\subsubsection*{Case 2: d=3}
	
	\begin{lemma}\label{'lem5.16}
		Let $d=3$. For every $n \geq 2,\ N \geq 1$ and $T>0$, there exists a constant $C_{n, N, T, \gamma}>0$ \rm{(}independent of $r$\rm{)} such that for any $0 \leq t \leq T,\ x, a \in \mathbb{R}^3$ and $r \geq 1$,
		\begin{align}\label{'5.59*}
			\mathbb{E}\big[V_{n,r}(t, x)^{2 N}\big] \leq C_{n, N, T, \gamma}.
		\end{align}
	\end{lemma}
	\begin{proof}
		This is similar to the proof of Lemma \ref{lem5.3}. Indeed, 		\begin{align}\label{'5.59}
			\E\big[V_{2,r}(t,x)^{2N}\big] \leq &\ 2^{2N}  \E\Big[\Big(\int_0^t ds \int p_{t-s} (x,z)V_{1,r}(s,z)^2 dz\Big)^{2N} \Big]\nn\\
			&+2^{2N} \E\Big[ \Big(\int_0^t  \int p_{t-s}(x,z)  V_{2,r}(s,z) W(ds,z) dz\Big)^{2N}\Big].
		\end{align}
		For the first expectation on the right-hand side above, we get
		\begin{align}\label{'5.60}
			J_{1,r}:= \E\Big[&\Big(\int_0^t ds \int p_{t-s} (x,z)V_{1,r}(s,z)^2 dz\Big)^{2N} \Big]\nn\\
			=\E\Big[&\int_0^t  ds_1 \int  p_{t-s_1} (x,z_1) d z_1 \cdots \int_0^t  ds_{2N-1} \int  p_{t-s_{2N-1}} (x,z_{2N-1})   d z_{2N-1} \nn\\
			&\int_0^t  ds_{2N} \int  p_{t-s_{2N}} (x,z_{2N})  \prod_{i=1}^{2N}  V_{1,r}(s_i, z_i)^2  d z_{2N}  \Big].
		\end{align}
		By a generalized H\" older's inequality and \eqref{'5.50}, we see
		\begin{align}
			\E\Big[ \prod_{i=1}^{2N}  V_{1,r}(s_i, z_i)^2 \Big] 
			&\leq 	\prod_{i=1}^{2N} \Big( \E\Big[   V_{1,r}(s_i, z_i)^{4N} \Big] \Big)^{1/(2N)} \leq C_{N,T} \prod_{i=1}^{2N} Q_{s_i}\phi_{a,\gamma}(z_i)^2\nn.
		\end{align}
		We then apply the above inequality to \eqref{'5.60} to obtain
		\begin{align}\label{'5.61}
			J_{1,r}&\leq C_{N,T} \Big(\int_0^t   ds  \int  p_{t-s } (x,z) \big(Q_{s}\phi_{a,\gamma}(z)\big)^2 d z \Big)^{2N}\nn\\
			&\leq C_{N,T} \Big(\int_0^t   ds  \int  p_{t-s } (x,z) \Big(\int_0^s du\int\frac{1}{|a-y|^\gamma}p_u(z,y)dy\Big)^2 d z \Big)^{2N}.
		\end{align}
		Now for any $0\le t\le T$, we proceed to estimate
		\begin{align*}
			I(\gamma,t):=\int_0^t   ds  \int  p_{t-s } (x,z) \Big(\int_0^s du\int\frac{1}{|a-y|^\gamma}p_u(z,y)dy\Big)^2 d z.
		\end{align*}
		When $0<\gamma<2$, it follows from Lemma \ref{lem5.1'} that
		\begin{align}\label{'bu5.60}
			\int_0^s du\int\frac{1}{|a-y|^\gamma}p_u(z,y)dy\le C_{\gamma} s^{1-\frac{\gamma}{2}}
		\end{align}
		with the constant $C_{\gamma}>0$. Then we have
		\begin{align*}
			I(\gamma,t)\le C_{\gamma} \int_0^t s^{2-\gamma}ds\int p_{t-s}(x,z)dz\le C_{\gamma,T}.
		\end{align*}
		When $\gamma=2$, using Lemma 5.1 of \cite{H18} gets that 
		\begin{align}\label{'bu5.62}
			\int_0^s du\int\frac{1}{|a-y|^2}p_u(z,y)dy
			\le 2\log^+\frac{1}{|a-z|}+C_T. 
		\end{align}
		Note that for any $x, y\in\mathbb{R}^d$,
		\begin{align}\label{'bu5.63}
			\Big( \log^+\frac{1}{|x-y|} \Big)^2 \leq \frac{1}{|x-y|}\quad \text{and} \quad \log^+\frac{1}{|x-y|} \leq \frac{1}{|x-y|}
		\end{align}
		hold for $d\geq1$. Then by Lemma 3.2 of \cite{H18} we have
		\begin{align*}
			I(\gamma,t)&\le \int_0^t   ds  \int  p_{t-s } (x,z)\Big[4\Big(\log^+\frac{1}{|a-z|} \Big)^2+4 C_T\log^+\frac{1}{|a-z|}\Big] dz+ TC_T^2\nn\\
			&\le C_T \int_0^t   ds  \int  p_{t-s } (x,z) \frac{1}{|a-z|}dz+C_T\le C_T.
		\end{align*}
		When $2<\gamma<3$, by the estimate in \cite[p.24]{H18} and Lemma \ref{lem5.1'} we obtain
		\begin{align*}
			I(\gamma,t)\le C_\gamma \int_0^t   ds  \int  p_{t-s } (x,z) \frac{1}{|a-z|^{2(\gamma-2)}}dz\le C_{\gamma,T}.
		\end{align*}
		Combining the above gives that for any $0<\gamma<3$,
		\begin{align}\label{'5.62}
			I(\gamma,t)\le C_{\gamma,T},\quad \forall\ 0\le t\le T, 
		\end{align}
		with the constant $C_{\gamma,T}>0$. 
		Then we conclude that 
		$J_{1,r}\le C_{N,T,\gamma}$. On the other hand, we repeat the proof in Lemma \ref{lem5.3} with $V^{\phi}_2(s,z)$ replaced by $V_{2,r}(s,z)$ to get
		\begin{align*}
		J_{2,r}&:=\E\Big[ \Big(\int_0^t  \int p_{t-s}(x,z)  V_{2,r}(s,z) W(ds,z) dz\Big)^{2N}\Big]\nn\\
		&\le C_{N,T} \int_0^t ds \int p_{t-s}(x,z)  \E\big[V_{2,r}(s,z)^{2N}\big] dz. 
		\end{align*}
		By a similar argument used in Lemma \ref{lem5.3}, we get the desired result for $n=2$ and further for any $n\ge 2$.  
	\end{proof}
	
	\begin{lemma}\label{'lem5.17}
		Let $d=3$. Suppose \eqref{de.1.5} holds for $d=3$. For every $n \geq 1$ and $T>0$,  there exists some constant $C_{n, T, \gamma}>0$ {\rm(}independent of $r${\rm)} such that for any $r \geq 1$, $0 \leq t \leq T$ and $a\in\mathbb{R}^3$ with $|a|\leq T$,
		\begin{align}\label{'5.65*}
			\mathbb{E}_{\mu}\big[Y_{t}(\phi_{a,\gamma,r})^{n}\big] \leq C_{n, T, \gamma}.
		\end{align}
	\end{lemma}
	\begin{proof}
		We shall omit some parts of the proof since the arguments are similar to those in Lemma \ref{'lem5.14}. Based on this, it is sufficient to give a new estimate of $\mathbb{E}[L_{n,r}(t)^N]$ for $n=1$. Indeed, when $n=1$, for any $N\ge 1$ we have
		\begin{align*}
\mathbb{E}[L_{1,r}(t)^N] =\mathbb{E}[\langle \mu,V_{1,r}(t,\cdot)\rangle^N]&=\int \mu(d x_1) \cdots \int \mathbb{E}\bigg[\prod_{i=1}^N V_{1,r}(t, x_i)\bigg] \mu(d x_n)\nn\\
	&\le C_{N, T}\bigg[\int Q_t \phi_{a,\gamma}(x) \mu(d x)\bigg]^N,
\end{align*}
where the inequality follows by (\ref{'5.50}). For any $0\le t\le T$, by \eqref{'bu5.60}, \eqref{'bu5.62} and \eqref{'bu5.63} we gets
\begin{align*}
	\int Q_t \phi_{a,\gamma}(x) \mu(d x)
	&=\int \mu(dx)\int_0^t ds\int \frac{1}{|a-y|^\gamma}p_s(x,y)dy\nn\\
	&\leq 
	\begin{cases}
		 C_{\gamma}  t^{1-\frac{\gamma}{2}}\mu(1) \leq C_{\gamma}  T^{1-\frac{\gamma}{2}}\mu(1), &\text{if } 0<\gamma<2; \\
		2 \int \frac{1}{|a-x|} \mu(dx) +C_T \mu(1), &\text{if } \gamma=2; \\
		C_{\gamma} \int \frac{1}{|a-x|^{\gamma-2}} \mu(dx), &\text{if } 2<\gamma<3,
	\end{cases}
\end{align*}
where the estimate in $2<\gamma<3$ comes from that in \cite[p.24]{H18}. Note that for any $\gamma\in(2,3)$,
\begin{align*}
\int\frac{1}{|a-x|^{\gamma-2}}\mu(dx) 
&\leq \int_{\{|a-x|<1\}} \frac{1}{|a-x|^{\gamma-2}}\mu(dx) + \mu(1)  \nn\\
&\leq \int \frac{1}{|a-x|}\mu(dx) + \mu(1).
\end{align*}
Since we assume \eqref{de.1.5} holds for $d=3$, the two preceding inequalities imply that for any $0\leq t\leq T$ and $a\in\mathbb{R}^3$ with $|a|\leq T$,
\begin{align*}
	\int Q_t \phi_{a,\gamma}(x) \mu(d x) \leq C_{\gamma, T}.
\end{align*}
Here the constant $C_{\gamma, T}$ depends on $\mu$ and we omit it to ease notation. Therefore, the desired result for $n=1$ is proved, and thus for any $n\ge 1$ follows by Lemma \ref{'lem5.16}.
	\end{proof}

	\subsection{Bounds on $X_t$}
	
	
	
	Recall that the sequence $\{\phi_{a,\gamma,r}\}\subset C_b^{2,+}(\mathbb{R}^d)$ is given by \eqref{xphir} such that $\text{as}\ r\rightarrow\infty$,
\begin{align*}
	\phi_{a,\gamma,r}\uparrow \phi_{a,\gamma},\quad a\in \mathbb{R}^d,
\end{align*}
where $\phi_{a,\gamma}$ is defined by \eqref{xphig}, $0<\gamma<1$ when $d=2$ and $0<\gamma<5/2$ when $d=3$. We denote $\widetilde{V}^{\phi}_{n}(t,x)$ and $\widetilde{L}^{\phi}_n(t)$  with respect to the text function $\phi=\phi_{a,\gamma,r}$ by $\widetilde{V}_{n,r}(t,x)$ and $\widetilde{L}_{n,r}(t)$, respectively. The dependence of $a$ and $\gamma$ will also be suppressed for notational ease. Fix $N\geq 1$ and $T>0$. It follows from \cite[Lemma 2.7]{MX07} ($\phi_t$ in the proof there is $\widetilde{V}_{1,r}$) and a generalized H\"older's inequality that for any $0\leq t_1, \cdots, t_N\leq T$ and $x_1, \cdots, x_N \in \R^d$, we have  
		\begin{align}\label{'b5.26}
			\E\big[\widetilde{V}_{1,r}(t_1, x_1)\cdots \widetilde{V}_{1,r}(t_N, x_N)\big]\leq {C}_{N,T} \prod_{i=1}^N {P}_{t_i} \phi_{a,\gamma,r}(x_i)\leq {C}_{N,T} \prod_{i=1}^N {P}_{t_i} \phi_{a,\gamma}(x_i),
		\end{align}
		where the constant ${C}_{N,T}>0$ is independent of $\phi_{a,\gamma,r}$.
	Recall from (\ref{ve10.5'}) that $\widetilde{V}_{n,r}(t,x)$ satisfies  
	\begin{align}\label{'b5.26*}
		\widetilde{V}_{n,r}(t,x)=  &\sum_{k=1}^{n-1} \binom{n-1}{k}   \int_0^t ds \int p_{t-s} (x,z)\widetilde{V}_{n-k,r}(s, z)\widetilde{V}_{k,r}(s, z) dz\nn\\
		&+\int_0^t  \int p_{t-s}(x,z)  \widetilde{V}_{n,r}(s,z) W(ds,z) dz,\quad\forall\ n\geq2.
	\end{align}
	Based on (\ref{'b5.26}) and (\ref{'b5.26*}), for $d=2,3$, we shall give upper bounds of $\mathbb{E}[\widetilde{V}_{n,r}(t,x)^{2N}]$ and further of the moments of $X_{t}(\phi_{a,\gamma,r})$ , where the bounds are all independent of $r$.
	
	\subsubsection*{Case 1: d=2 and $\mathbf{0<\gamma<1}$}
	
	\begin{lemma}\label{'lemb5.6}
		Let $d=2$ and $0<\gamma<1$. For every $n \geq 2,\ N \geq 1$ and $T>0$, there exists a constant $C_{n, N, T, \gamma}>0$ {\rm(}independent of $r${\rm)} such that for any $0 \leq t \leq T,\ x, a \in \mathbb{R}^2$ and $r \geq 1$,
		\begin{align}\label{'b5.27}
			\mathbb{E}\big[\widetilde{V}_{n,r}(t, x)^{2 N}\big] \leq C_{n, N, T, \gamma}.
		\end{align}
	\end{lemma}

	\begin{proof}
		This is similar to the proof of Lemma \ref{lem5.3}. Indeed, 		
		\begin{align}\label{'b5.28}
			\E\big[\widetilde{V}_{2,r}(t,x)^{2N}\big] \leq &\ 2^{2N}  \E\Big[\Big(\int_0^t ds \int p_{t-s} (x,z)\widetilde{V}_{1,r}(s,z)^2 dz\Big)^{2N} \Big]\nn\\
			&+2^{2N} \E\Big[ \Big(\int_0^t  \int p_{t-s}(x,z)  \widetilde{V}_{2,r}(s,z) W(ds,z) dz\Big)^{2N}\Big].
		\end{align}
		For the first expectation on the right-hand side above, we get
		\begin{align}\label{'b5.29}
			\widetilde{J}_{1,r}:= \E\Big[&\Big(\int_0^t ds \int p_{t-s} (x,z)\widetilde{V}_{1,r}(s,z)^2 dz\Big)^{2N} \Big]\nn\\
			=\E\Big[&\int_0^t  ds_1 \int  p_{t-s_1} (x,z_1) d z_1 \cdots \int_0^t  ds_{2N-1} \int  p_{t-s_{2N-1}} (x,z_{2N-1})   d z_{2N-1} \nn\\
			&\int_0^t  ds_{2N} \int  p_{t-s_{2N}} (x,z_{2N})  \prod_{i=1}^{2N}  \widetilde{V}_{1,r}(s_i, z_i)^2  d z_{2N}  \Big].
		\end{align}
		By a generalized H\"older's inequality and \eqref{'b5.26}, we see that for any $0\leq s_i\leq t \leq T$,
		\begin{align}
			\E\Big[ \prod_{i=1}^{2N}  \widetilde{V}_{1,r}(s_i, z_i)^2 \Big] 
			\leq 	\prod_{i=1}^{2N} \Big( \E\Big[   \widetilde{V}_{1,r}(s_i, z_i)^{4N} \Big] \Big)^{1/(2N)} \leq C_{N,T} \prod_{i=1}^{2N} P_{s_i}\phi_{a,\gamma}(z_i)^2\nn.
		\end{align}
		It follows from Lemma 3.1 of \cite{H18} that
		\begin{align}
			P_{s_i}\phi_{a,\gamma}(z_i)=\int\frac{1}{|a-y|^\gamma}p_{s_i}(z_i,y)dy\le C_{\gamma}\frac{1}{|a-z_i|^\gamma}\nn.
		\end{align}
		For any $0\le t\le T$, we then apply the above inequality to \eqref{'b5.29} and use Lemma \ref{lem5.1'} under $0<\gamma<1$ to obtain
		\begin{align*}
			\widetilde{J}_{1,r}&\leq C_{N,T,\gamma} \Big(\int_0^t   ds  \int  p_{t-s } (x,z) \frac{1}{|a-z|^{2\gamma}} d z \Big)^{2N}
			\leq C_{N,T,\gamma}\cdot T^{2N(1-\gamma)}=C_{N,T,\gamma}.
		\end{align*}
		 On the other hand, we repeat the proof in Lemma \ref{lem5.3} with $V^{\phi}_2(s,z)$ replaced by $\widetilde{V}_{2,r}(s,z)$ to get
		\begin{align*}
		\widetilde{J}_{2,r}&:=\E\Big[ \Big(\int_0^t  \int p_{t-s}(x,z)  \widetilde{V}_{2,r}(s,z) W(ds,z) dz\Big)^{2N}\Big]\nn\\
		&\le C_{N,T} \int_0^t ds \int p_{t-s}(x,z)  \E\big[\widetilde{V}_{2,r}(s,z)^{2N}\big] dz. 
		\end{align*}
		By a similar argument used in Lemma \ref{lem5.3}, we get the desired result for $n=2$ and further for any $n\ge 2$.  
	\end{proof}

	\begin{lemma}\label{'lemb5.7}
		Let $d=2$ and $0<\gamma<1$. For every $n \geq 1$ and $0<\eps_0<T$, there exists some constant $C_{n, \eps_0,T, \gamma}>0$ {\rm(}independent of $r${\rm)} such that for any $\eps_0 \leq t \leq T$, $a\in\mathbb{R}^2$ and $r \geq 1$,
		\begin{align}\label{'b5.40}
			\mathbb{E}_{\mu}\big[X_{t}(\phi_{a,\gamma,r})^{n}\big] \leq C_{n, \eps_0, T, \gamma}.
		\end{align}
	\end{lemma}
	\begin{proof}
		Recall from Proposition \ref{l2.2} with $\phi=\phi_{a,\gamma,r}$ that
		\begin{align}\label{''by1}
\mathbb{E}_\mu[X_t(\phi_{a,\gamma,r})^n]=\mathbb{E}[\widetilde{L}_{n,r}(t)] ,
\end{align}
Fix $0<\eps_0<T$ and let $\eps_0\le t\le T$, for any $n\ge 1,\ a \in \mathbb{R}^2$ and $r \geq 1$, we shall prove that
\begin{align}\label{'b5.41}
\mathbb{E}[\widetilde{L}_{n, r}(t)^N] \leq C_{n, N, \eps_0, T,\gamma},\quad \forall\ N\ge 1.
\end{align} 
The conclusion follows immediately from 
(\ref{''by1}) and the above with $N=1$. 
When $n=1$, for any $N\ge 1$ we have
\begin{align*}
\mathbb{E}[\widetilde{L}_{1,r}(t)^N] =\mathbb{E}[\langle \mu,\widetilde{V}_{1,r}(t,\cdot)\rangle^N]&=\int \mu(d x_1) \cdots \int \mathbb{E}\bigg[\prod_{i=1}^N \widetilde{V}_{1,r}(t, x_i)\bigg] \mu(d x_n)\nn\\
	&\le C_{N,T} \Big[ \int P_t \phi_{a,\gamma}(x) \mu(dx)\Big]^N \nn\\
	&= C_{N,T} \Big[ \int \mu(dx) \int \frac{1}{|a-y|^{\gamma}}p_t(x, y) dy\Big]^N \nn\\
	&\le C_{N,T,\gamma} t^{-\gamma N/2}\mu(1)^N \nn\\
	&\leq C_{N,T,\gamma}\eps_0^{-\gamma N/2}\mu(1)^N :=C_{N,\eps_0,T,\gamma},\quad \forall\ \eps_0\leq t\leq T,
\end{align*}
where the second inequality follows by $d-\gamma-1=1-\gamma>0$ and \eqref{8.1}, and we omit the initial measure $\mu$ in the notation of the above constant.
Hence (\ref{'b5.40}) holds for the case $n=1$. Assume that (\ref{'b5.40}) holds for all $1\le k\le n-1$ with some $n\ge 2$. Then for the case $n$, using Lemma \ref{'lemb5.6} and a similar argument used in Lemma \ref{lem5.5} gives the desired result.
	\end{proof}
	
	\subsubsection*{Case 2: d=3 and $\mathbf{0<\gamma<5/2}$}
	
	\begin{lemma}\label{'lemb5.9}
		Given $d=3$ and $0<\gamma<5/2$. For every $n \geq 2,\ N \geq 1$ and $0<\eps_0<T$, there exists a constant $C_{n, N, \eps_0, T, \gamma}>0$ {\rm(}independent of $r${\rm)} such that for any $\eps_0 \leq t \leq T,\ x, a \in \mathbb{R}^3$ and $r \geq 1$,
		\begin{align}\label{'b5.59*}
			\mathbb{E}\big[\widetilde{V}_{n, r}(t, x)^{2 N}\big] \leq C_{n, N, \eps_0, T, \gamma}.
		\end{align}
	\end{lemma}
	
	\begin{proof}
		The proof is similar to that of Lemma \ref{lem5.3}, so we only describe the difference. It suffices to prove that \eqref{'b5.59*} holds for $n=2$. Note that
		\begin{align}\label{be.5.93}
			\E\big[\widetilde{V}_{2,r}(t, x)^{2N}\big] \leq &\ 2^{2N}  \E\Big[\Big(\int_0^t ds \int p_{t-s} (x, z)\widetilde{V}_{1,r}(s, z)^2 dz\Big)^{2N} \Big]\nn\\
			&+2^{2N} \E\Big[ \Big(\int_0^t  \int p_{t-s}(x, z)  \widetilde{V}_{2,r}(s, z) W(ds, z) dz\Big)^{2N}\Big].\end{align}
			For the first term on the right-hand side in the above, it is obvious that
			\begin{align}\label{be.5.94}
			& \E\Big[\Big(\int_0^t ds \int p_{t-s} (x, z)\widetilde{V}_{1,r}(s, z)^2 dz\Big)^{2N} \Big] \nn\\
					 =&\ \int_0^t ds_1\int p_{t-s_1}(x, z_1)dz_1 \cdots \int_0^t ds_{2N-1}\int p_{t-s_{2N-1}}(x, z_{2N-1})dz_{2N-1} \nn\\
					 &\ \int_0^t ds_{2N}\int p_{t-s_{2N}}(x, z_{2N}) \E\Big[ \prod_{i=1}^{2N}\widetilde{V}_{1,r}(s_i, z_i)^2\Big] dz_{2N}.	
			\end{align}
			By a generalized H\"older's inequality and \eqref{'b5.26} we get that for any $0\leq s_i\leq t \leq T$,
			\begin{align}\label{be.5.95}
			\E\Big[ \prod_{i=1}^{2N}  \widetilde{V}_{1,r}(s_i, z_i)^2 \Big] 
			\leq 	\prod_{i=1}^{2N} \Big( \E\Big[   \widetilde{V}_{1,r}(s_i, z_i)^{4N} \Big] \Big)^{1/(2N)} \leq C_{N,T} \prod_{i=1}^{2N} P_{s_i}\phi_{a,\gamma}(z_i)^2,\quad \forall\ r\geq1.
			\end{align}
			Note that the interval $(0\vee(1-2/\gamma),1\wedge(3/\gamma-1))\neq\varnothing$ when $0<\gamma<5/2$, then we take
			\begin{align*}
				\kappa:= \kappa(\gamma)=  \frac{1}{2\gamma}	\wedge 1 \in \Big(0\vee(1-\frac{2}{\gamma}), 1\wedge(\frac{3}{\gamma}-1)\Big),
			\end{align*}
			which implies 
			\begin{align*}
				\gamma(1-\kappa)<2	\quad\text{and}\quad \gamma(1+\kappa)<3.
			\end{align*}
			Hence,
			\begin{align}\label{be.5.96}
				P_{s_i}\phi_{a,\gamma}(z_i)^2 =  P_{s_i}\phi_{a,\gamma}(z_i)^{1+\kappa} \times P_{s_i}\phi_{a,\gamma}(z_i)^{1-\kappa}.
			\end{align}
			By Jensen's inequality we get
			\begin{align*}
				P_{s_i}\phi_{a,\gamma}(z_i)^{1+\kappa} = \Big[ \int \phi_{a,\gamma}(y) p_{s_i}(z_i, y) dy \Big]^{1+\kappa} \leq \int \phi_{a,\gamma}(y)^{1+\kappa} p_{s_i}(z_i, y) dy.
			\end{align*}
			By \eqref{8.1}, we have
			\begin{align*}
				P_{s_i}\phi_{a,\gamma}(z_i) = \int \frac{1}{|a-y|^{\gamma}} p_{s_i}(z_i, y) dy	\leq C_{\gamma} s_i^{-\gamma/2}, 	
				\end{align*}
				here the constant $C_{\gamma}$ is finite since $d-\gamma-1=2-\gamma>-1/2$. By the two preceding inequalities, \eqref{be.5.96} becomes
				\begin{align*}
					P_{s_i}\phi_{a,\gamma}(z_i)^2\leq C_{\gamma}^{1-\kappa} s_i^{-\gamma(1-\kappa)/2} \int \phi_{a,\gamma}(y)^{1+\kappa} p_{s_i}(z_i, y) dy .
				\end{align*}
				This together with \eqref{be.5.94} and \eqref{be.5.95} implies that for any $\eps_0\leq t \leq T$,
				\begin{align*}
					& \E\Big[\Big(\int_0^t ds \int p_{t-s} (x, z)\widetilde{V}_{1,r}(s, z)^2 dz\Big)^{2N} \Big] \nn\\
					 \leq &\ C_{N,T} C_{\gamma}^{2N(1-\kappa)} \prod_{i=1}^{2N} \Big[ \int_0^t s_i^{-\gamma(1-\kappa)/2}ds_i \int p_{t-s_i} (x, z_i) dz_i \int \phi_{a,\gamma}(y)^{1+\kappa} p_{s_i}(z_i, y) dy \Big] \nn\\
					 = &\ C_{N,T,\gamma} \Big[ \int \phi_{a,\gamma}(y)^{1+\kappa} p_t(x, y)dy\Big]^{2N} \times\Big[ \int_0^{t} s^{-\gamma(1-\kappa)/2}ds \Big]^{2N}.
				\end{align*}
				Since $\gamma(1-\kappa)<2$, we have
				\begin{align*}
					\int_{0}^t  s^{-\gamma(1-\kappa)/2}ds = \Big(1-\frac{\gamma(1-\kappa)}{2}\Big)^{-1} t^{1-\frac{\gamma(1-\kappa)}{2}}	.		
					\end{align*}
				A similar calculation based on \eqref{8.1} shows that
				\begin{align*}
					\int \phi_{a,\gamma}(y)^{1+\kappa} p_t(x, y)dy =\int \frac{1}{|a-y|^{\gamma(1+\kappa)}} p_{t}(x, y) dy	\leq C_{\gamma} t^{-\gamma(1+\kappa)/2} ,
				\end{align*}
				where $C_{\gamma}$ is finite since $2-\gamma(1+\kappa)>-1$.
				Combining the three preceding inequalities gives that for any $\eps_0\leq t\leq T$,
				\begin{align*}
					\E\Big[\Big(\int_0^t ds \int p_{t-s} (x, z)\widetilde{V}_{1,r}(s, z)^2 dz\Big)^{2N} \Big] &\leq C_{N,T,\gamma} t^{2N(1-\gamma)} \nn\\
					&\leq C_{N,T,\gamma} (\eps_0^{2N(1-\gamma)}+T^{2N(1-\gamma)})\nn\\
					&= C_{N,\eps_0,T,\gamma},
				\end{align*}
				where the second inequality comes from $\gamma\in(0,5/2)$.
				By repeating the remaining proof of Lemma \ref{lem5.3} for $n=2$ and replacing $V^{\phi}_2(s, r)$ with $\widetilde{V}_{2,r}(s, r)$ we obtain \eqref{'b5.59*} holds for $n=2$. Hence, we get the desired result.
	\end{proof}
	
	\begin{lemma}\label{'lemb5.10}
		Given $d=3$ and $0<\gamma<5/2$. For every $n \geq 1$ and $0<\eps_0<T$, there exists some constant $C_{n, \eps_0, T, \gamma}>0$ {\rm(}independent of $r${\rm)} such that for any $\eps_0 \leq t \leq T$, $a\in\mathbb{R}^3$ and $r \geq 1$,
		\begin{align}\label{'be5.40}
			\mathbb{E}_{\mu}\big[X_{t}(\phi_{a,\gamma,r})^{n}\big] \leq C_{n, \eps_0, T, \gamma}.
		\end{align}
	\end{lemma}
	
	\begin{proof}
		We shall omit some parts of the proof since the arguments are similar to those in Lemma \ref{'lemb5.7}. Based on this and Lemma \ref{'lemb5.9}, it is sufficient to give a new estimate of $\mathbb{E}[\widetilde{L}_{n,r}(t)^N]$ for $n=1$.
	When $n=1$, for any $N\geq1$ and $0\le t\le T$, by \eqref{'b5.26} we have
	\begin{align*}
		\mathbb{E}[\widetilde{L}_{1,r}(t)^N] =\mathbb{E}[\langle \mu,\widetilde{V}_{1,r}(t,\cdot)\rangle^N]&=\int \mu(d x_1) \cdots \int \mathbb{E}\bigg[\prod_{i=1}^N \widetilde{V}_{1,r}(t, x_i)\bigg] \mu(d x_n)\nn\\
	&\le C_{N,T} \Big[ \int P_t \phi_{a,\gamma}(x) \mu(dx)\Big]^N.	
	\end{align*}
	Since $2-\gamma>-1/2$, it follows from Lemma \ref{lem5.1'} that for any $\eps_0\le t\le T$,
	\begin{align*}
		\int P_t \phi_{a,\gamma}(x) \mu(dx)
		&=\int \mu(dx)\int \frac{1}{|a-y|^\gamma}p_t(x,y)dy\nn\\
		&\le \mu(1) C_{\gamma} t^{-\gamma/2}
		\le \mu(1) C_{\gamma}\eps_0^{-\gamma/2}.
	\end{align*}
	Then for any $\eps_0\le t\le T$ and $r\ge 1$,
	\begin{align*}
		\mathbb{E}[\widetilde{L}_{1,r}(t)^N]\le C_{ N, \eps_0, T, \gamma}.	
	\end{align*}
	Therefore, the desired result for $n=1$ is proved.
	\end{proof}
	
\section{Proofs of \eqref{e.4.15} and \eqref{e4.28}}\label{appc}

\noindent\textbf{Proof of \eqref{e4.28}.} 	Recall from \eqref{e4.3} to see that for any $h,h^{\prime}>0$ and $x\in\mathbb{R}^d$,
		\begin{align*}
			V_t^{p_h^x, p_{h^{\prime}}^x}(u, v) &\leq e^{t\|g\|_{\infty}} \bigg[ \int_0^t p_{s+h}^x(u) d s \int_s^{t}  p_{r+h^{\prime}}^x(v) dr + \int_0^t p_{s+h^{\prime}}^x(v) d s \int_s^{t}  p_{r+h}^x(u) dr \bigg],\quad \forall\  u,v\in\mathbb{R}^d.
		\end{align*}
	 
		(i) When $d=1$, it is easy to check that for any $h,h^{\prime}>0$,
			\begin{align}\label{e4.25}
				V_t^{p_h^x, p_{h^{\prime}}^x}(u, v) \leq 2e^{t\|g\|_{\infty}} \int_{0}^t s^{-1/2} ds \int_{s}^t r^{-1/2} dr \leq  8t e^{t\|g\|_{\infty}}:=\bar{V}_{t,1},
			\end{align}
			where the last inequality comes from $\int_{s}^t r^{-1/2} dr \leq 2 t^{1/2}$. Since $\mu$ is a finite measure, we obtain
			\begin{align}
				&\int\!\!\int \bar{V}_{t,1} \mu(du) \mu(dv)= 8\mu(1)^2  te^{t\|g\|_{\infty}} <\infty,\label{e4.26}\\
				&\int_0^t ds \int \mu(du) \int p_{t-s}(u,v) \hk{\bar{V}_{s,1}} dv = 8\mu(1)t^2 e^{t\|g\|_{\infty}} <\infty.\label{e4.27}
			\end{align}
		
			(ii) When $d=2$, by a simple calculation, we get that for any $0<h,h^{\prime}<1$, 
			\begin{align}\label{e.4.16}
				 \int_0^t p_{s+h}^x(u) d s \int_s^{t}  p_{r+h^{\prime}}^x(v) dr \leq \int_0^{t+1} p_{s}^x(u) d s \int_0^{t+1}  p_{r}^x(v) dr.
			\end{align}
			Thus, we get
			\begin{align}\label{e4.8}
				V_t^{p_h^x, p_{h^{\prime}}^x}(u, v) &\leq e^{t\|g\|_{\infty}} \bigg( \int_0^{t+1} p_{s}^x(u) d s \int_0^{t+1}  p_{r}^x(v) dr  + \int_0^{t+1} p_{s}^x(v) d s \int_0^{t+1}  p_{r}^x(u) dr \bigg) \nonumber\\
				&\leq 2e^{t\|g\|_{\infty}} \int_0^{t+1} p_{s}^x(u) d s \int_0^{t+1}  p_{r}^x(v) dr :=\bar{V}_{t,2,1}^x(u,v).
			\end{align}
			We next consider the integrability of $\bar{V}_{t,2,1}^x(u,v)$. \hk{It follows from \cite[Appendix C(i)]{H18} with $\alpha=1$ that	 for any $x,u\in\mathbb{R}^2$,		
			\begin{align*}
				\int_{0}^{\infty} e^{-s} p_s^x(u) ds \leq C+  \log^{+} \frac{1}{|x-u|},
			\end{align*} 
			where the constant $C>0$. Thus, 
			\begin{align}\label{e.4.12}
				\int_0^t p_s^x(u) ds \leq e^{t} \int_{0}^{t} e^{-s} p_s^x(u)  ds \leq e^{t} \int_{0}^{\infty} e^{-s} p_s^x(u)  ds  
				\leq C e^{t}+e^t\log^+\frac{1}{|x-u|}.
			\end{align}}
			Therefore,
			\begin{align}\label{e4.2.2}
				\int\!\!\int \bar{V}_{t,2,1}^x(u,v) \mu(du)\mu(dv) &= 2e^{t\|g\|_{\infty}} \bigg( \int \mu(du)  \int_0^{t+1} p_{s}^x(u) d s \bigg)^2 \nonumber\\
				&\hk{\leq 2e^{t\|g\|_{\infty}} \bigg( Ce^{t+1}\mu(1)+e^{t+1}\int \log^{+} \frac{1}{|x-u|} \mu(du) \bigg)^2  }<\infty,
			\end{align}
			where the finiteness comes from \hk{the assumption \eqref{de.1.5} holds} for $d=2$. In particular, if $u=v$, recall from (3.43) of Sugitani \cite{S89} that 
			\begin{align}\label{e4.9}
				p_t^x(v)p_s^x(v) \leq (st)^{-d/4} p_{ts/(t+s)}^x(v),\quad s,t>0,\ x,v\in\mathbb{R}^d.
			\end{align}
			It follows from \eqref{e4.8} that for any $x,v\in\mathbb{R}^2$,
			\begin{align}\label{e4.2.3}
				V_t^{p_h^x, p_{h^{\prime}}^x}(v, v) \leq 2e^{t\|g\|_{\infty}} \int_{0}^{t+1} ds \int_{0}^{t+1} (sr)^{-1/2} p_{sr/(s+r)}^x(v) dr := \bar{V}_{t,2,2}^x(v,v).
			\end{align}
		We shall prove $\bar{V}_{t,2,2}^x(v,v)$ is integrable. Assume \eqref{de.1.5} holds for $d=2$, then 
		\begin{align}\label{e4.11}
			&\int_0^t ds \int \mu(du) \int p_{t-s}(u,v) \bar{V}_{s,2,2}^x(v,v) dv \nonumber\\
			=&\ 2 \int _0^t e^{s\|g\|_{\infty}}  ds \int \mu(du) \int_{0}^{s+1} ds_1 \int_{0}^{s+1} (s_1r)^{-1/2}  dr \int p_{t-s}(u,v) p_{s_1r/(s_1+r)}^x(v) dv \nonumber\\
			=&\ 2 \int _0^t e^{s\|g\|_{\infty}}  ds \int \mu(du) \int_{0}^{s+1} ds_1 \int_{0}^{s+1} (s_1r)^{-1/2}p_{t-s+\frac{s_1r}{s_1+r}}^x(u)   dr  \nonumber\\
			\leq&\ 2  e^{t\|g\|_{\infty}}   \int \mu(du) \int_{0}^{t+1} ds_1 \int_{0}^{t+1} (s_1r)^{-1/2} dr \int _0^t p_{t-s+\frac{s_1r}{s_1+r}}^x(u)  ds \nonumber\\
			\leq&\ 2  e^{t\|g\|_{\infty}}   \int \mu(du) \int_{0}^{t+1} ds_1 \int_{0}^{t+1} (s_1r)^{-1/2} dr \int _0^{2t+1} p_{s}^x(u)  ds \nonumber\\
			\leq&\ \hk{2  e^{t\|g\|_{\infty}}    \int_{0}^{t+1} ds_1 \int_{0}^{t+1} (s_1r)^{-1/2} dr  \int \Big(C e^{2t+1}+e^{2t+1} \log^+\frac{1}{|x-u|} \Big) \mu (du)}  \nonumber\\
			=&\ 8 (t+1) e^{t\|g\|_{\infty}}  \hk{\Big(C e^{2t+1} \mu(1)+ e^{2t+1} \int\log^+\frac{1}{|x-u|}  \mu (du) \Big)} <\infty,
		\end{align}
		where the last inequality comes from \eqref{e.4.12}.

		(iii) When $d=3$, \hk{\eqref{e4.8} still holds, that is}, for any $0<h,h^{\prime}<1$,
		\begin{align}\label{e4.10}
			V_t^{p_h^x, p_{h^{\prime}}^x}(u, v) \leq 2e^{t\|g\|_{\infty}} \int_0^{t+1} p_{s}^x(u) d s \int_0^{t+1}  p_{r}^x(v) dr :=\bar{V}_{t,3,1}^x(u,v).
		\end{align}
		(We point out that $\bar{V}_{t,2,1}^x(u,v)\neq \bar{V}_{t,3,1}^x(u,v)$ since the dimensions are different.) Note that
		\begin{align}\label{e.4.13}
			\int_{0}^{\infty} p_s^x(u) ds = \frac{1}{2\pi} \frac{1}{|x-u|},\quad x,u\in\mathbb{R}^3.
		\end{align}
		We shall consider the integrability. Suppose that \eqref{de.1.5} holds for $d=3$, then
		\begin{align}\label{e4.3.2}
			\int\!\!\int \bar{V}_{t,3,1}^x(u,v) \mu(du)\mu(dv) &\leq 2e^{t\|g\|_{\infty}} \bigg( \int \mu(du) \int_0^{\infty} p_{s}^x(u) d s \bigg)^2 \nonumber\\
			&=\frac{1}{2\pi^2} e^{t\|g\|_{\infty}} \bigg( \int \frac{1}{|x-u|} \mu(du) \bigg)^2<\infty.
		\end{align}
		On the other hand, we apply \eqref{e4.9} to \eqref{e4.10} for $u=v$ to get that for any $x,v\in\mathbb{R}^3$, 
		\begin{align}\label{e4.3.3}
			V_t^{p_h^x, p_{h^{\prime}}^x}(v, v) \leq 2e^{t\|g\|_{\infty}} \int_{0}^{t+1} ds \int_{0}^{t+1} (sr)^{-3/4} p_{sr/(s+r)}^x(v) dr := \bar{V}_{t,3,2}^x(v,v).
		\end{align}
		Based on \eqref{e.4.13}, by repeating the arguments for deriving \eqref{e4.11} we  obtain
		\begin{align}\label{e4.3.4}
			&\int_0^t ds \int \mu(du) \int p_{t-s}(u,v) \bar{V}_{s,3,2}^x(v,v) dv \nonumber\\
			&\leq\ 2  e^{t\|g\|_{\infty}}   \int \mu(du) \int_{0}^{t+1} ds_1 \int_{0}^{t+1} (s_1r)^{-3/4} dr \int _0^t p_{t-s+\frac{s_1r}{s_1+r}}^x(u)  ds \nonumber\\
			&\leq\ \hk{2  e^{t\|g\|_{\infty}}   \int_{0}^{t+1} ds_1 \int_{0}^{t+1} (s_1r)^{-3/4} dr \bigg( \frac{1}{2\pi} \int \frac{1}{|x-u|} \mu (du) \bigg)} \nonumber\\
			&= \frac{16}{\pi} (t+1)^{1/2} e^{t\|g\|_{\infty}} \int \frac{1}{|x-u|} \mu (du) <\infty,
		\end{align}
		when \eqref{de.1.5} holds for $d=3$.
		
	Based on the above, once we establish \hk{the above pointwise convergences of $V_t^{p_h^x, p_{h^{\prime}}^x}(u, v)$ and $V_s^{p_h^x, p_{h^{\prime}}^x}(v, v)$} (i.e., \eqref{e.4.15} and \eqref{bue.4.15}), together with \hk{the integral upper bounds $\bar{V}_{t,1}$,$\bar{V}_{t,2,1}^x(u,v)$, $\bar{V}_{s,2,2}^x(v,v),\bar{V}_{t,3,1}^x(u,v)$ and  $\bar{V}_{s,3,2}^x(v,v)$}, we apply dominated convergence to \hk{the right-hand side of \eqref{e.4.14}} to get \eqref{e4.28}.
	$\hfill\blacksquare$\\
	
	\noindent\textbf{Proof of \eqref{e.4.15}.}
	Recall from \eqref{e10.23} to see that
		\begin{align*}
			 V_t^{p_h^x, p_{h^{\prime}}^x}(u, v)=&\Pi_{(u, v)}\bigg\{\int_0^t \Big[p_h^x(B_s) Q_{t-s} p_{h^{\prime}}^x(\tilde{B}_s) + p_{h^{\prime}}^x(\tilde{B}_s) Q_{t-s} p_{h}^x(B_s) \Big] e^{\int_0^s g(B_r, \tilde{B}_r) d r} d s\bigg\}.
		\end{align*}
		Note that
		\begin{align*}
			e^{\int_0^s g(B_r, \tilde{B}_r) d r}-1=\int_0^s g(B_r, \tilde{B}_r) e^{\int_0^r g(B_a, \tilde{B}_a) d a} d r.
		\end{align*}
		It then follows that
		\begin{align*}
			V_t^{p_h^x,  p_{h^{\prime}}^x}(u, v) =&\ \Pi_{(u, v)}\bigg\{\int_0^t \Big[p_h^x(B_s) Q_{t-s} p_{h^{\prime}}^x(\tilde{B}_s) + p_{h^{\prime}}^x(\tilde{B}_s) Q_{t-s} p_{h}^x(B_s) \Big] ds \bigg\}\nn \\
			& + \Pi_{(u, v)}\bigg\{\int_0^t p_h^x(B_s) Q_{t-s} p_{h^{\prime}}^x(\tilde{B}_s)  ds \int_0^s g(B_r, \tilde{B}_r) e^{\int_0^r g(B_a, \tilde{B}_a) da} dr  \bigg\} \nonumber\\
			& + \Pi_{(u, v)}\bigg\{\int_0^t p_{h^{\prime}}^x(\tilde{B}_s) Q_{t-s} p_{h}^x(B_s)  ds \int_0^s g(B_r, \tilde{B}_r) e^{\int_0^r g(B_a, \tilde{B}_a) da} dr  \bigg\} \nonumber\\
			=&\ \int_0^t p_{s+h}^x(u) ds \int_{s}^t p_{r+h^{\prime}}^x(v) dr + \int_0^t p_{s+h^{\prime}}^x(v) ds \int_{s}^t p_{r+h}^x(u) dr\nonumber\\
			& +\int_0^t ds \int_0^s \Pi_{(u, v)}\Big\{ p_h^x(B_s) Q_{t-s} p_{h^{\prime}}^x(\tilde{B}_s) \cdot  g(B_r, \tilde{B}_r) e^{\int_0^r g(B_a, \tilde{B}_a) da}  \Big\} dr\nn\\
			& +\int_0^t ds \int_0^s \Pi_{(u, v)}\bigg\{ p_{h^{\prime}}^x(\tilde{B}_s) Q_{t-s} p_{h}^x(B_s) \cdot  g(B_r, \tilde{B}_r) e^{\int_0^r g(B_a, \tilde{B}_a) da}  \Big\} dr.
		\end{align*}
		By the Markov property \hk{we have
		\begin{align*}
			&\ \Pi_{(u, v)}\bigg\{ \Big[p_{h}^x(B_s) Q_{t-s} p_{h^{\prime}}^x(\tilde{B}_s)+p_{h^{\prime}}^x(\tilde{B}_s) Q_{t-s} p_{h}^x(B_s) \Big]  g(B_r, \tilde{B}_r) e^{\int_0^r g(B_a, \tilde{B}_a) da} \bigg\} \nonumber\\
			=&\ \Pi_{(u, v)}\bigg\{ \Big[p_{s-r+h}^x(B_r) Q_{t-s} p_{s-r+h^{\prime}}^x(\tilde{B}_r) + p_{s-r+h^{\prime}}^x(\tilde{B}_r) Q_{t-s} p_{s-r+h}^x(B_r) \Big]   g(B_r, \tilde{B}_r) e^{\int_0^r g(B_a, \tilde{B}_a) da}  \bigg\}. 
		\end{align*}}
		Hence, 
		\begin{align*}
			V_t^{p_h^x,  p_{h^{\prime}}^x}(u, v)= &\int_0^t p_{s+h}^x(u) ds \int_{s}^t p_{r+h^{\prime}}^x(v) dr + \int_0^t p_{s+h^{\prime}}^x(v) ds \int_{s}^t p_{r+h}^x(u) dr \\
			&+ \int_0^t ds \int_0^s \Pi_{(u, v)}\Big\{ p_{s-r+h}^x(B_r) Q_{t-s} p_{s-r+h^{\prime}}^x(\tilde{B}_r) \cdot g(B_r, \tilde{B}_r) e^{\int_0^r g(B_a, \tilde{B}_a) da}  \Big\} dr\nn\\
			&+ \int_0^t ds \int_0^s \Pi_{(u, v)}\Big\{ p_{s-r+h^{\prime}}^x(\tilde{B}_r) Q_{t-s} p_{s-r+h}^x(B_r) \cdot g(B_r, \tilde{B}_r) e^{\int_0^r g(B_a, \tilde{B}_a) da}  \Big\} dr.\nn
		\end{align*}
	\hk{Since 
	\begin{align*}
		0\leq \mathbf{1}_{\{h\leq s\leq t+h\}} p_s^x(u) \mathbf{1}_{\{s-h+h^{\prime}\leq r\leq t+h^{\prime}\}} p_r^x(v) \leq \mathbf{1}_{\{0\leq s\leq t+1\}} p_s^x(u) \mathbf{1}_{\{0\leq r\leq t+1\}} p_r^x(v) 
	\end{align*}
	for $0<h,h^{\prime}<1$. By \eqref{e.4.12} and \eqref{e.4.13}, we use dominated convergence to get that for $\mu\times\mu$-almost every $(u,v)\in\mathbb{R}^{2d}$,
	\begin{align}\label{e.4.17}
		\lim_{h, h^{\prime}\downarrow0} \int_0^t p_{s+h}^x(u) ds \int_{s}^t p_{r+h^{\prime}}^x(v) dr &=\lim_{h, h^{\prime}\downarrow0} \int_h^{t+h} p_{s}^x(u) ds \int_{s-h+h^{\prime}}^{t+h^{\prime}} p_{r}^x(v) dr\nonumber\\
		&= \int_0^t p_{s}^x(u) ds \int_{s}^t p_{r}^x(v) dr,
	\end{align}
	and similarly,
	\begin{align*}
		\lim_{h, h^{\prime}\downarrow0} \int_0^t p_{s+h^{\prime}}^x(v) ds \int_{s}^t p_{r+h}^x(u) dr &= \int_0^t p_{s}^x(v) ds \int_{s}^t p_{r}^x(u) dr,	
	\end{align*}
when \eqref{de.1.5} holds for $d=2,3$.} It suffices to show that
	\begin{align}\label{e.4.20}
		\int_0^t ds &\int_{0}^s \bigg|  \Pi_{(u, v)}\Big\{p_{s-r+h}^x(B_r) Q_{t-s} p_{s-r+h^{\prime}}^x(\tilde{B}_r) \cdot g(B_r, \tilde{B}_r) e^{\int_0^r g(B_a, \tilde{B}_a) da} \Big\} \nonumber\\ 
		&\qquad+\Pi_{(u, v)}\Big\{p_{s-r+h^{\prime}}^x(\tilde{B}_r) Q_{t-s} p_{s-r+h}^x(B_r) \cdot  g(B_r, \tilde{B}_r) e^{\int_0^r g(B_a, \tilde{B}_a) da} \Big\} \nonumber\\ 
		&\qquad- \Pi_{(u, v)}\Big\{p_{s-r}^x(B_r) Q_{t-s} p_{s-r}^x(\tilde{B}_r) \cdot   g(B_r, \tilde{B}_r) e^{\int_0^r g(B_a, \tilde{B}_a) da}  \Big\} \nonumber\\
		&\qquad- \Pi_{(u, v)}\Big\{p_{s-r}^x(\tilde{B}_r) Q_{t-s} p_{s-r}^x(B_r) \cdot   g(B_r, \tilde{B}_r) e^{\int_0^r g(B_a, \tilde{B}_a) da}  \Big\} \bigg| dr \xrightarrow[h,h^{\prime}\downarrow0]{}0,
	\end{align}
	$\text{for }\mu\times\mu\text{-almost every }(u,v)\in\mathbb{R}^{2d}.$ Note that the integrand function in the above inequality is bounded by the sum of
	\begin{align*}
		A_{1,h,h^{\prime}}(s,r):= \bigg|  \Pi_{(u, v)}\bigg\{ \Big[&p_{s-r+h}^x(B_r) Q_{t-s} p_{s-r+h^{\prime}}^x(\tilde{B}_r) \nonumber\\
		&- p_{s-r}^x(B_r) Q_{t-s} p_{s-r}^x(\tilde{B}_r)  \Big] g(B_r, \tilde{B}_r) e^{\int_0^r g(B_a, \tilde{B}_a) da} \bigg\} \bigg|
	\end{align*}
	and
	\begin{align*}
		A_{2,h,h^{\prime}}(s,r):=\bigg|  \Pi_{(u, v)}\bigg\{ \Big[&p_{s-r+h^{\prime}}^x(\tilde{B}_r) Q_{t-s} p_{s-r+h}^x({B}_r) \nonumber\\
		&- p_{s-r}^x(\tilde{B}_r) Q_{t-s} p_{s-r}^x(B_r)  \Big] g(B_r, \tilde{B}_r) e^{\int_0^r g(B_a, \tilde{B}_a) da} \bigg\} \bigg|.
	\end{align*}
	Thus, we first prove $A_{1,h,h^{\prime}}(s,r)$  \hk{tends} to zero as $h,h^{\prime}\downarrow0$. \hk{Note that}
	\begin{align}\label{e.4.18}
		A_{1,h,h^{\prime}}(s,r) \leq  &\ \bigg|  \Pi_{(u, v)}\bigg\{ \Big[p_{s-r+h}^x(B_r) Q_{t-s} p_{s-r+h^{\prime}}^x(\tilde{B}_r) \nonumber\\
		&\qquad\qquad\quad- p_{s-r}^x(B_r) Q_{t-s} p_{s-r+h^{\prime}}^x(\tilde{B}_r)  \Big] g(B_r, \tilde{B}_r) e^{\int_0^r g(B_a, \tilde{B}_a) da} \bigg\} \bigg| \nonumber\\
		&\hk{+}\  \bigg|  \Pi_{(u, v)}\bigg\{ \Big[p_{s-r}^x(B_r) Q_{t-s} p_{s-r+h^{\prime}}^x(\tilde{B}_r) \nonumber\\
		&\qquad\qquad\quad- p_{s-r}^x(B_r) Q_{t-s} p_{s-r}^x(\tilde{B}_r)  \Big] g(B_r, \tilde{B}_r) e^{\int_0^r g(B_a, \tilde{B}_a) da} \bigg\} \bigg| \nonumber\\
		\leq&\ \|g\|_{\infty}  e^{r\|g\|_{\infty}} \Pi_{(u, v)}\bigg\{ Q_{t-s} p_{s-r+h^{\prime}}^x(\tilde{B}_r) \Big|p_{s-r+h}^x(B_r)-p_{s-r}^x(B_r) \Big| \bigg\} \nonumber\\
		&+ \|g\|_{\infty}  e^{r\|g\|_{\infty}} \Pi_{(u, v)}\bigg\{ p_{s-r}^x(B_r) Q_{t-s}\Big|p_{s-r+h^{\prime}}^x-p_{s-r}^x\Big|(\tilde{B}_r) \bigg\} \nonumber\\
		\leq&\ \|g\|_{\infty}  e^{r\|g\|_{\infty}} Q_{t-s} p_{s+h^{\prime}}^x(v)  \Pi_{(u, v)}\bigg\{\Big|p_{s-r+h}^x(B_r)-p_{s-r}^x(B_r)\Big|\bigg\} \nonumber\\
		&+ \|g\|_{\infty}  e^{r\|g\|_{\infty}}  p_s^x(u)  \Pi_{(u, v)}\bigg\{Q_{t-s}\Big|p_{s-r+h^{\prime}}^x-p_{s-r}^x\Big|(\tilde{B}_r)\bigg\},
	\end{align}
	where the last inequality follows by the independence between $B_r$ and $\tilde{B}_r$. 
	By \hk{choosing the similar integrable bounded function of \eqref{e.4.17} and then using dominated convergence,} we have 
	\begin{align}\label{e.4.19}
		\lim_{h^{\prime}\downarrow0} Q_{t-s} p_{s+h^{\prime}}^x(v) =Q_{t-s} p_{s}^x(v)\hk{,\quad\text{for $\mu$-almost every } v\in\mathbb{R}^d}.
	\end{align}It is easy to check that
	\begin{align}\label{e.4.21}
		\Big|p_{s-r+h}^x(B_r)-p_{s-r}^x(B_r)\Big| \leq p_{s-r+h}^x(B_r)+p_{s-r}^x(B_r)
	\end{align}
	and for $\mu\times\mu$-almost every $(u,v)\in\mathbb{R}^{2d}$,
	\begin{align*}
		& \lim_{h\downarrow0} \Pi_{(u,v)} \Big\{ p_{s-r+h}^x(B_r)+p_{s-r}^x(B_r) \Big\} = 2 p_s^x(u) = \Pi_{(u,v)} \Big\{ \lim_{h\downarrow0} \big( p_{s-r+h}^x(B_r)+p_{s-r}^x(B_r) \big) \Big\},\nonumber\\
		& \lim_{h^{\prime}\downarrow0} \Pi_{(u,v)} \Big\{ Q_{t-s} ( p_{s-r+h^{\prime}}^x+p_{s-r}^x)(\tilde{B}_r) \Big\}=2Q_{t-s} p_s^{x}(v) = \Pi_{(u,v)} \Big\{ \lim_{h^{\prime}\downarrow0} Q_{t-s} ( p_{s-r+h^{\prime}}^x+p_{s-r}^x)(\tilde{B}_r) \Big\}.
	\end{align*}
	Based on the above two relations, we use the generalized dominated convergence (see, e.g., Exercise 20 of Chp.2 of \cite{F99}) to obtain that for $\mu\times\mu$-almost every $(u,v)\in\mathbb{R}^{2d}$,
	\begin{align*}
		&\lim _{h \downarrow 0} \Pi_{(u, v)}\Big\{\big|p_{s-r+h}^x(B_r)-p_{s-r}^x(B_r)\big|\Big\}=\Pi_{(u, v)}\Big\{\lim _{h \downarrow 0}\big|p_{s-r+h}^x(B_r)-p_{s-r}^x(B_r)\big|\Big\}=0,\\
		&\lim _{h^{\prime} \downarrow 0} \Pi_{(u, v)}\Big\{Q_{t-s} \big| p_{s-r+h^{\prime}}^x-p_{s-r}^x\big|(\tilde{B}_r)\Big\}=\Pi_{(u, v)}\Big\{\lim _{h^{\prime} \downarrow 0}Q_{t-s} \big| p_{s-r+h^{\prime}}^x-p_{s-r}^x\big|(\tilde{B}_r)\Big\}=0.
	\end{align*} 
	Combining the above two convergences with \eqref{e.4.18} and \eqref{e.4.19} yields 
	\begin{align*}
		\lim_{h,h^{\prime}\downarrow0} A_{1,h,h^{\prime}}(s,r)=0,\quad{\mu\times\mu}\text{-a.e.}.
	\end{align*}
	In view of \hk{\eqref{e.4.18}-\eqref{e.4.21}}, we apply the generalized dominated convergence to obtain
	\begin{align*}
		\lim_{h,h^{\prime}\downarrow0} \int_0^t ds \int_{0}^s A_{1,h,h^{\prime}}(s,r) dr=0,\quad\mu\times\mu\text{-a.e.}.
	\end{align*}
	By repeating the proof of the above convergence, we also have 
	\begin{align*}
		\lim_{h,h^{\prime}\downarrow0} \int_0^t ds \int_{0}^s A_{2,h,h^{\prime}}(s,r) dr=0,\quad\mu\times\mu\text{-a.e.}.
	\end{align*} 
	Hence, these two convergences imply \eqref{e.4.20} and then \eqref{e.4.15} hold. 
	 $\hfill\blacksquare$

\section{Proof of \eqref{hc1} and \eqref{hc2}}\label{AA}

	When $d=2$ and $\alpha>0$, we consider the following function
	\begin{align*}
		f_{\alpha}(x)=\int_0^{\infty} e^{-\alpha s} p_s(0,x) ds - \frac{1}{\pi} \log^+\frac{1}{|x|},\quad x\in\mathbb{R}^2.
	\end{align*}
	Let $\gamma\in (0,1)$ be given in the assumption \eqref{de1.3}, we shall prove \eqref{hc1}, i.e., $\mu f_{\alpha}(a)$ is locally $\gamma$-H\"older continuous with respect to $a$.\\ 

\noindent\textbf{Proof of \eqref{hc1}.}	
	Note that
 	\begin{align*}
 		f_{\alpha}(x)
 		&=\Big[ \int_0^{1} \frac{1}{2\pi s} e^{-\frac{|x|^2}{2s}} ds - \frac{1}{\pi} \log^+\frac{1}{|x|} \Big] + \int_0^1 \frac{1}{2\pi s}(e^{-\alpha s}-1) e^{-\frac{|x|^2}{2s}} ds + \int_1^{\infty} \frac{1}{2\pi s}e^{-\alpha s} e^{-\frac{|x|^2}{2s}}ds  \nonumber\\
 		&:= f_{\alpha,1}(|x|)+f_{\alpha,2}(|x|)+f_{\alpha,3}(|x|).
 	\end{align*}
 	Hence, for any $a, b\in\mathbb{R}^2$,
 	\begin{align}\label{bua.3}
 		\big|\mu f_{\alpha}(a)-\mu f_{\alpha}(b)\big| &\leq \int \big| f_{\alpha}(a-x)-f_{\alpha}(b-x) \big| \mu(dx) \nn\\
 		&\leq \sum_{l=1}^3 \int \big|f_{\alpha,l}(|a-x|)- f_{\alpha,l}(|b-x|)\big| \mu(dx).
 	\end{align}
 	Note that $f_{\alpha,l}\ (l=1,2,3)$ are actually functions on $[0,\infty)$. We claim that for every $l=1,2,3$,
 	\begin{align}\label{bua.4}
 		\big|f_{\alpha,l}(u)- f_{\alpha,l}(v)\big| \leq C_{\alpha,l} |u-v|,\quad \forall\ u,v\in[0,\infty).
 	\end{align}
	It then follows that
	\begin{align*}
		\big|f_{\alpha,l}(|a-x|)-f_{\alpha,l}(|b-x|)\big| 
		\leq C_{\alpha,l} \big||a-x|- |b-x|\big| \leq C_{\alpha,l} |a-b|.
	\end{align*}
	This together with \eqref{bua.3} gives the local $\gamma$-H\"older continuity. So it suffices to prove \eqref{bua.4} holds.
 	
 	(i) By the definition of $f_{\alpha,3}$, we see that for any $u, v\geq0$,
 	\begin{align}\label{aaa.4}
 		\big|f_{\alpha,3}(u)-f_{\alpha,3}(v)\big| 
 		=  \Big| \int_1^{\infty} \frac{1}{2\pi s} e^{-\alpha s} \Big( e^{-\frac{u^2}{2s}} -e^{-\frac{v^2}{2s}} \Big) ds \Big|  
 		\leq \int_1^{\infty} s^{-1} \Big| e^{-\frac{u^2}{2s}} -e^{-\frac{v^2}{2s}} \Big| ds.
 	\end{align}
 	By the mean value theorem, we have
 	\begin{align}\label{bua.7}
 		\Big|e^{-\frac{u^2}{2s}} -e^{-\frac{v^2}{2s}} \Big|=|u-v| \frac{\xi}{s}e^{-\frac{\xi^2}{2s}}=|u-v|\sqrt{\frac{\xi^2}{s^2}e^{-\frac{\xi^2}{s}}},	
 	\end{align}
	as $\xi\geq0$ lies between $u$ and $v$. 
	Since $\sup_{r\geq 0} r e^{-r}<\infty$, we have
	\begin{align*}
		\big|f_{\alpha,3}(u)-f_{\alpha,3}(v)\big| 
 		\leq |u-v| \int_1^{\infty} s^{-3/2}	\sqrt{\frac{\xi^2}{s} e^{-\frac{\xi^2}{s}}} ds 
 		\leq C|u-v|.
	\end{align*}
	
	(ii) Note that $\sup_{0\leq s\leq 1} s^{-1} (1-e^{-\alpha s})=C_{\alpha}<\infty$, then together with \eqref{bua.7} we get that for any $u, v\geq0$,
	\begin{align*}
		\big| f_{\alpha,2}(u)- f_{\alpha,2}(v)\big| 
 		&\leq \int_0^1 s^{-1} (1-e^{-\alpha s}) \Big| e^{-\frac{u^2}{2s}} -e^{-\frac{v^2}{2s}} \Big| ds \nonumber\\
 		&\leq C_{\alpha} |u-v| \int_0^1  \sqrt{\frac{\xi^2}{s^2} e^{-\frac{\xi^2}{s}}} ds.
 	\end{align*}
	By Jensen's inequality we get
	\begin{align*}
		\int_0^1  \sqrt{\frac{\xi^2}{s^2} e^{-\frac{\xi^2}{s}}} ds \leq \Big( \xi^2 \int_0^1  s^{-2} e^{-\frac{\xi^2}{s}} ds \Big)^{1/2}	= e^{-\frac{\xi^2}{2}}\leq 1.
	\end{align*}
	It then follows that
	\begin{align*}
		\big|f_{\alpha,2}(u)-f_{\alpha,2}(v)\big| 
 		\leq  C_{\alpha} |u-v|.	
	\end{align*}
	
	(iii) We shall estimate the difference of $f_{\alpha,1}(u)$ with respect to $u\in[0,\infty)$ in the following three cases.  
	(a) For $u, v\geq1$, by \eqref{bua.7} we obtain
	\begin{align}\label{aaa.12}
		\big| f_{\alpha,1}(u)- f_{\alpha,1}(v) \big|
		&= \Big| \int_0^1 \frac{1}{2\pi s} e^{-\frac{u^2}{2s}}ds - \int_0^1 \frac{1}{2\pi s} e^{-\frac{v^2}{2s}}ds\Big| \nonumber\\	
		&\leq \int_0^1 s^{-1} \Big| e^{-\frac{u^2}{2s}}-e^{-\frac{v^2}{2s}}\Big| ds \nonumber\\
		&\leq|u-v| \int_0^1 \xi s^{-2}  e^{-\frac{\xi^2}{2s}} ds \nonumber\\
		&\leq \frac{2|u-v|}{\xi}\leq 2|u-v|,
	\end{align}
	as we note $\xi\geq \min\{u, v\}\geq1$. 
	
	(b) For $0\leq u, v\leq1$, a simple calculation shows that if $0\leq u\leq1$,
	\begin{align*}
		f_{\alpha,1}(u)	&= \int_0^1 \frac{1}{2\pi s} e^{-\frac{u^2}{2s}}ds - \frac{1}{\pi} \log \frac{1}{u}	\nonumber\\
		&=\Big[\int_{u^2/2}^1 \frac{1}{2\pi s}ds - \frac{1}{\pi} \log\frac{1}{u} \Big]+\int_{u^2/2}^1 \frac{1}{2\pi s}(e^{-s}-1)ds +\int_1^{\infty} \frac{1}{2\pi s} e^{-s} ds,
	\end{align*}
	and
	\begin{align*}
		\int_{u^2/2}^1 \frac{1}{2\pi s}ds - \frac{1}{\pi} \log\frac{1}{u} = \frac{1}{2\pi} \log 2,\qquad \int_1^{\infty} \frac{1}{2\pi s} e^{-s} ds<\infty.
	\end{align*}
	Note that $\sup_{r\geq0}r^{-1}(1-e^{-r^2/2})<\infty$, then we use the mean value theorem to obtain
	\begin{align*}
		\big| f_{\alpha,1}(u)-f_{\alpha,1}(v) \big|	 
		&= \Big| \int_{u^2/2}^1 \frac{1}{2\pi s}(e^{-s}-1)ds - \int_{v^2/2}^1 \frac{1}{2\pi s}(e^{-s}-1)ds\Big| \nn\\
		&= |u-v| \frac{1}{\pi\hat{\xi}} \Big(1-e^{-\frac{\hat{\xi}^2}{2}}\Big) \leq C|u-v|,	
	\end{align*}
	where $\hat{\xi}\geq0$ lies between $u$ and $v$. 
		
	(c) For $u\geq1, 0\leq v\leq1$, by the assertions (a) and (b), we have
		\begin{align*}
			\big| f_{\alpha,1}(u)-f_{\alpha,1}(v) \big| 
			&\leq \big| f_{\alpha,1}(u)-f_{\alpha,1}(1) \big| + \big| f_{\alpha,1}(1)-f_{\alpha,1}(v) \big|\nonumber\\
			&\leq C \big(\big| u-1 \big| + \big| 1-v \big|\big) \nonumber\\
			&= C (u-v)= C |u-v|.
		\end{align*}
	Similarly, the above inequality holds for $0\leq u\leq1, v\geq1$. Summing up, 
	\begin{align*}
		\big|f_{\alpha,1}(u)- f_{\alpha,1}(v)\big| \leq C |u-v|.
	\end{align*}
	Therefore, \eqref{bua.4} follows by the above assertions (i)-(iii).
	$\hfill\blacksquare$\\
	
	When $d=3$ and $\alpha\geq0$, we consider the different function 
	\begin{align*}
			\widetilde{f}_{\alpha}(x)=\int_0^{\infty} e^{-\alpha s} p_s(0,x) ds - \frac{1}{2\pi|x|},\quad x\in\mathbb{R}^3.
	\end{align*}
	Recall that 
 	\begin{align*}
		\mu\widetilde{f}_{\alpha}(a)= \int \widetilde{f}_{\alpha}(a-x) \mu(dx).
	\end{align*}
	In the following, we shall prove \eqref{hc2}, i.e.,  $\mu\widetilde{f}_{\alpha}(a)$ is locally $\gamma$-H\"older continuous with respect to $a$, where $\gamma\in (0,1/2)$ is given in the assumption \eqref{de1.3}. \\
	
	\noindent\textbf{Proof of \eqref{hc2}.}
	If $\alpha=0$, then $\widetilde{f}_{\alpha}(x)=0$, so the local $\gamma$-H\"older continuity for $\mu\widetilde{f}_{\alpha}$ obviously holds. If $\alpha>0$, note that
	\begin{align*}
		\widetilde{f}_{\alpha}(x)=	\int_0^{\infty} e^{-\alpha s} p_s(0,x) ds - \int_0^{\infty} p_s(0,x) ds = \int \big( \frac{1}{2\pi s}\big)^{3/2} ( e^{-\alpha s} -1) e^{-\frac{|x|^2}{2s}} ds.
	\end{align*}
	For any $a=(a_1,a_2,a_3), b=(b_1,b_2,b_3)\in\mathbb{R}^3$,
	\begin{align}\label{aa3.1}
		\big| \mu\widetilde{f}_{\alpha}(a)-\mu\widetilde{f}_{\alpha}(b)	 \big| 
		&\leq \int |\widetilde{f}_{\alpha}(a-x) - \widetilde{f}_{\alpha}(b-x) \big| \mu(dx) \nonumber \\
		&\leq \int \mu(dx) \int_0^{\infty} s^{-3/2} (1-e^{-\alpha s}) \Big| e^{-\frac{|a-x|^2}{2s}} -e^{-\frac{|b-x|^2}{2s}} \Big| ds.
	\end{align}
	By the mean value theorem, there exists $\theta\in(0,1)$ such that
	\begin{align*}
		e^{-\frac{|a-x|^2}{2s}} -e^{-\frac{|b-x|^2}{2s}} = \sum_{i=1}^{3}-\frac{\xi_i(x)}{s} (a_i-b_i) e^{-\frac{|\xi(x)|^2}{2s}},
	\end{align*}
	where $\xi(x)=(\xi_1(x),\xi_2(x),\xi_3(x))$ and $\xi_i(x)=b_i-x_i+\theta(a_i-b_i)$ for $i=1,2,3$. It then follows that
	\begin{align*}
		\Big| e^{-\frac{|a-x|^2}{2s}} -e^{-\frac{|b-x|^2}{2s}} \Big| 
		&\leq \sqrt{3 \sum_{i=1}^3 \frac{\xi_i(x)^2}{s^2} (a_i-b_i)^2 e^{-\frac{|\xi(x)|^2}{s}} } \nonumber\\
		&\leq |a-b| \sqrt{3\frac{|\xi(x)|^2}{s^2}e^{-\frac{|\xi(x)|^2}{s}}}.
	\end{align*}
	This combined with \eqref{aa3.1} yields that
	\begin{align}\label{aa3.2}
		  \big| \mu\widetilde{f}_{\alpha}(a)-\mu\widetilde{f}_{\alpha}(b)	 \big| 
		\leq  &\int \mu(dx) \int_0^1 s^{-3/2} (1-e^{-\alpha s}) \Big| e^{-\frac{|a-x|^2}{2s}} -e^{-\frac{|b-x|^2}{2s}} \Big|^\gamma \Big| e^{-\frac{|a-x|^2}{2s}} -e^{-\frac{|b-x|^2}{2s}} \Big|^{1-\gamma} ds \nonumber\\
		&+ \int \mu(dx) \int_1^{\infty} s^{-3/2} (1-e^{-\alpha s}) \Big| e^{-\frac{|a-x|^2}{2s}} -e^{-\frac{|b-x|^2}{2s}} \Big| ds \nonumber\\
		\leq & 2^{1-\gamma} 3^{\gamma/2} |a-b|^{\gamma} \int_0^1 s^{-(3+\gamma)/2} (1-e^{-\alpha s}) \Big( \frac{|\xi(x)|^2}{s} e^{-\frac{|\xi(x)|^2}{s}} \Big)^{\gamma/2} ds \nonumber\\
		&+ \sqrt{3} |a-b| \int \mu(dx) \int_{1}^{\infty} s^{-2} \sqrt{\frac{|\xi(x)|^2}{s} e^{-\frac{|\xi(x)|^2}{s}}} ds \nonumber\\
		\leq & C_{\gamma} \mu(1) |a-b|^{\gamma} \int_0^1 s^{-(1+\gamma)/2}s^{-1}(1-e^{-\alpha s}) ds + C \mu(1)|a-b|\int_1^{\infty} s^{-2} ds \nonumber\\
		\leq & C_{\alpha,\gamma} \mu(1) |a-b|^{\gamma} \int_0^1 s^{-(1+\gamma)/2} ds + C \mu(1)|a-b| \nonumber\\
		\leq & C_{\alpha,\gamma} \mu(1) |a-b|^{\gamma} + C \mu(1)|a-b|,
	\end{align}
	where the third and the fourth inequalities follow by the two facts $\sup_{r\geq 0} r e^{-r}<\infty$ and $\sup_{0\leq s\leq 1} s^{-1}(1-e^{-\alpha s})=C_{\alpha}<\infty$ for $\alpha>0$, respectively, the last inequality comes from $\gamma<1/2<1$. Hence, the local $\gamma$-H\"older continuity holds by using the inequality \eqref{aa3.2}.
	$\hfill\blacksquare$

\section{Proof of the Green function representation}\label{appg}
	In this section, we shall give the proof of Theorem \ref{cor3.4*}, which provides the Green function representation of $X$. Recall that the superprocess $X$ satisfies the martingale problem \eqref{eMP}. 

\begin{lemma}\label{lem3.3*}
	For	any $t\geq0$ and $\phi\in C_b^2(\mathbb{R}^d)$, we have 
	\begin{align}\label{bu.3.5*}
		X_t(\phi) = X_0(P_t \phi) + \int_0^t\int P_{t-s}\phi(x) dM(s,x),\quad \mathbb{P}_{\mu}\text{-a.s.},	
	\end{align}
	where $P_t$ is the transition probability given by \eqref{q} and $M$ is the martingale measure associated with $X$.
\end{lemma}
	

\begin{proof}
	The proof is quite standard. For example, it follows from the proof of Proposition II.5.7 in Perkins \cite{P02} that for the ``nice enough'' function $\phi_s:[0, t] \rightarrow C_b^{2}(\mathbb{R})$, we have
	\begin{align*}
X_t(\phi_t)-X_0(\phi_0)=\int_0^t X_s\Big(\frac{\Delta}{2} \phi_s+\dot{\phi}_s\Big) d s & +\int_0^t \int\phi_s(x) dM(s,x),\quad\mathbb{P}_{\mu}\text{-a.s.},
\end{align*}
where $\dot{\phi}_s:=\frac{\partial \phi_s}{\partial s}$.
Take $\phi_s=P_{t-s} \phi$ for $\phi \in C_b^{2}(\mathbb{R})$ to obtain (\ref{bu.3.5*}) as required.
\end{proof}
	
\begin{lemma}\label{lem3.3**}
	Let $d\leq 3$. Suppose that \eqref{de.1.5} holds for $d=2,3$. Recall that for any $a\in\mathbb{R}^d$, $g_\alpha^a(x)=\int_0^{\infty} e^{-\alpha t} p_t^a(x) d t$ is a function on $\mathbb{R}^d$. Then for any $t>0$, $$\int_0^u \int P_{t-s}g_\alpha^a(x) d M(s,x),\ 0 \le u \le t$$ is a continuous $L^2$-bounded martingale with
\begin{align*}
\Big\langle\int_0^{\cdot} \int P_{t-s} g_\alpha^a(x) d M(s, x) \Big\rangle_{u}=&\ \int_0^u X_s\big((P_{t-s} g_\alpha^a)^2\big) d s\\&+\int_0^u d s \int\!\!\!\int g(z, w) P_{t-s} g_\alpha^a(z) P_{t-s} g_\alpha^a(w) X_s(d z) X_s(d w) .
\end{align*}
\end{lemma}

\begin{proof}
	We only need to show that for any given $u\in[0, t]$,
	\begin{align}\label{bu**1}
		\mathbb{E}_\mu\Big[\int_0^u X_s\big((P_{t-s} g_\alpha^a)^2\big) d s+\int_0^u d s \int\!\!\!\int g(z, w) P_{t-s} g_\alpha^a(z) P_{t-s} g_\alpha^a(w) X_s(d z) X_s(d w)\Big]<\infty.
	\end{align}
	For the first term, by (\ref{e.7.0}) and \eqref{e7.24} we have
	\begin{align*}
		\mathbb{E}_\mu\Big[\int_0^u X_s\big((P_{t-s} g_\alpha^a)^2\big) d s\Big]&\le C_\alpha \mathbb{E}_\mu\Big[ \int_0^u X_s\big((P_{t-s} g_0^a)^2\big) ds\Big]\\
		&\leq C_\alpha \mathbb{E}_\mu\Big[\int_0^u X_s\big((g_0^a)^2\big)ds\Big]\\
		&=C_\alpha \mathbb{E}_\mu\big[Y_u\big((g_0^a)^2\big)\big]<\infty,
		\end{align*}
		where the finiteness comes from \eqref{ebc.1} and the finiteness of $I_3(t)$. 
	Turning to the second term, it is easy to check that	
	\begin{align*}
		\mathbb{E}_\mu&\Big[\int_0^u d s \int\!\!\!\int g(z, w) P_{t-s} g_\alpha^a(z) P_{t-s} g_\alpha^a(w) X_s(d z) X_s(d w)\Big]\\
		&\le\|g\|_{\infty}\mathbb{E}_\mu\Big[\int_0^u \big(X_s(P_{t-s} g_\alpha^a)\big)^2 d s\Big]\le C_{\alpha} \|g\|_{\infty}\mathbb{E}_\mu\Big[\int_0^u \big(X_s(P_{t-s} g_0^a)\big)^2 d s\Big]\\
		&\le C_{\alpha} \|g\|_{\infty} \int_0^u \mathbb{E}_\mu\big[ \big(X_s(g_0^a)\big)^2 \big] d s.
	\end{align*}
	Since the finiteness of the above expectation follows from the proof of \eqref{e.7.28}, we have \eqref{bu**1} holds.
\end{proof}

Based on Lemmas \ref{lem3.3*} and \ref{lem3.3**}, the Green function representation of $X_t(g^a_\alpha)$ (Theorem \ref{cor3.4*}) now follows by monotone convergence and dominated convergence.

\section{Proof of Lemmas \ref{bu.lem.b.2} and \ref{bu.lem.b.3}}\label{appi}

\noindent\textbf{Proof of Lemma \ref{bu.lem.b.2}.}
	For any $0\leq s\leq t$ and $a\in\mathbb{R}^d$, a similar argument used in Lemma \ref{lem3.3**} shows that
	\begin{align*}
\int_0^u \int(P_{t-r}  g_\alpha^a(x)-P_{s-r} g_\alpha^a(x)) d M(r, x), \quad 0 \leq u \leq s
\end{align*}
	is a continuous $L^2$-bounded martingale with
	\begin{align*}
\langle\int_0^\cdot &\int(P_{t-r}  g_\alpha^a(x)-P_{s-r} g_\alpha^a(x)) d M(r, x)\rangle_u=\int_0^u X_r((P_{t-r}g_\alpha^a-P_{s-r}g_\alpha^a)^2)dr\\
&+\int_0^u d r \int\!\!\int g(z, w) \big(P_{t-r} g_\alpha^a(z)-P_{s-r} g_\alpha^a(z)\big) \big(P_{t-r} g_\alpha^a(w)-P_{s-r} g_\alpha^a(w)\big) X_r(d z) X_r(d w) .
\end{align*}
Then by Burkholder-Davis-Gundy's inequality we have
\begin{align*}
	I_1^{N,\alpha,a} (s,t)
	\le&\ 2^N C_N\Big\{\mathbb{E}_\mu\Big[\Big(\int_0^s X_r\big((P_{t-r} g_\alpha^a-P_{s-r} g_\alpha^a)^2\big) d r\Big)^N\Big]\nn\\
	&\ +\|g\|_{\infty}^{N}\ \mathbb{E}_\mu\Big[\Big(\int_0^s \big(X_r(|P_{t-r} g_\alpha^a-P_{s-r} g_\alpha^a|)\big)^2 d r\Big)^N\Big]\Big\}.
\end{align*}
By Cauchy-Schwarz's inequality and Lemma \ref{lemx5.1} we get that for any $0\leq s\leq t\leq T$,
\begin{align}\label{d.0}
	 I_1^{N,\alpha,a} (s,t) 
	\leq &\ \Big\{ 2^N C_N +\|g\|_{\infty}^N \Big( \mathbb{E}_{\mu}\big[\big(\sup_{0\leq r\leq T} X_r(1) \big)^{2N} \big] \Big)^{1/2} \Big\} \nn\\
	&\qquad\qquad\times\Big\{ \mathbb{E}_{\mu}\Big[ \Big( \int_0^s X_r\big( (P_{t-r}g_{\alpha}^a- P_{s-r}g_{\alpha}^a)^2 \big) dr\Big)^{2N}\Big] \Big\}^{1/2} \nn\\
	\leq &\ C_{N,T} \Big\{ \mathbb{E}_{\mu}\Big[ \Big( \int_0^s X_r\big( (P_{t-r}g_{\alpha}^a- P_{s-r}g_{\alpha}^a)^2 \big) dr\Big)^{2N}\Big] \Big\}^{1/2}.
\end{align}
Notice that
\begin{align*}
	\big|P_{t-r}  g_\alpha^a(x)-P_{s-r} g_\alpha^a(x)\big|
	&=\Big|\int\big(p_{t-r}(x, y)-p_{s-r}(x, y)\big) d y \int_0^{\infty} e^{-\alpha r_1} p_{r_1}(a, y) d r_1\Big|\\
	&\le \int_0^{\infty} e^{-\alpha r_1}\big|p_{t-r+r_1}(x,a)-p_{s-r+r_1}(x,a)\big|dr_1.
\end{align*}
Since $0\le r\le s\le t$, for any $\tilde{\delta}\in(0,1)$, by Lemma \ref{*lem6.1} we have
\begin{align*}
	\big|p_{t-r+r_1}(x,a)-p_{s-r+r_1}(x,a)\big|
	&\le \Big[(t-s)(s-r+r_1)^{-\frac{d}{2}-1}\Big]^{\tilde{\delta}}\Big[p_{t-r+r_1}(x,a)^{1-\tilde{\delta}}+p_{s-r+r_1}(x,a)^{1-\tilde{\delta}}\Big]\\
	&\le 2(t-s)^{\tilde{\delta}}(s-r+r_1)^{-\tilde{\delta}-\frac{d}{2}}.
\end{align*}
It follows that
\begin{align}\label{d.1}
	\big|P_{t-r}  g_\alpha^a(x)-P_{s-r} g_\alpha^a(x)\big|
	&\le 2(t-s)^{\tilde{\delta}}\int_0^{\infty}e^{-\alpha r_1}(s-r+r_1)^{-\tilde{\delta}-\frac{d}{2}}dr_1\nn\\
	&\le 2(t-s)^{\tilde{\delta}}\Big[\int_0^1 (s-r+r_1)^{-\tilde{\delta}-\frac{d}{2}}dr_1 +\frac{1}{\alpha}\Big].
\end{align}
Moreover, it is easy to obtain that for any $\tilde{\delta}\in (0,1)$,
		\begin{equation}\label{d.2}
			\int_0^1 (s-r+r_1)^{-\tilde{\delta}-\frac{d}{2}}dr_1
			\le\left\{
			\begin{aligned}
				& \Big|\frac{1}{2}-\tilde{\delta}\Big|^{-1}\Big[(s-r+1)^{\frac{1}{2}-\tilde{\delta}}+(s-r)^{\frac{1}{2}-\tilde{\delta}}\Big], & \text{if } d=1; \\
				& \tilde{\delta}^{-1}(s-r)^{-\tilde{\delta}}, & \text{if } d=2;\\
				&\Big(\frac{1}{2}+\tilde{\delta}\Big)^{-1}(s-r)^{-(\frac{1}{2}+\tilde{\delta})}, & \text{if } d=3.
			\end{aligned}
			\right.
		\end{equation}

(i) When $d=1$, by (\ref{d.1}) and (\ref{d.2}) we get for any $0\le s\le t$ and $\tilde{\delta}\in(0,1)$,
	\begin{align*}
\big|P_{t-r} g_\alpha^a(x)-P_{s-r} g_\alpha^a(x)\big| & \leq 2(t-s)^{\tilde{\delta}}\Big\{C_{\tilde{\delta}}\Big[(s-r+1)^{\frac{1}{2}-\tilde{\delta}}+(s-r)^{\frac{1}{2}-\tilde{\delta}}\Big]+\frac{1}{\alpha}\Big\} \\
& =(t-s)^{\tilde{\delta}}\Big[C_{\tilde{\delta}}(s-r+1)^{\frac{1}{2}-\tilde{\delta}}+C_{\tilde{\delta}}(s-r)^{\frac{1}{2}-\tilde{\delta}}+C_\alpha\Big] .
\end{align*}
Then for any $0\le s\le T$ and $\tilde{\delta}\in (0,1)$, by Lemmas \ref{lem5.5} and \ref{lemx5.1} we have
\begin{align}\label{d.3}
&\ \mathbb{E}_\mu\Big[\Big(\int_0^s X_r\big((P_{t-r} g_\alpha^a-P_{s-r} g_\alpha^a)^2\big) d r\Big)^{2N}\Big]\nn \\
\leq &\ \mathbb{E}_\mu\Big[\Big((t-s)^{2 \tilde{\delta}} \int_0^s\Big[C_{\tilde{\delta}}(s-r+1)^{\frac{1}{2}-\tilde{\delta}}+C_{\tilde{\delta}}(s-r)^{\frac{1}{2}-\tilde{\delta}}+C_\alpha\Big]^2 X_r(1) d r\Big)^{2N}\Big]\nn\\
\le &\ (t-s)^{4 \tilde{\delta} N}\Big\{C_{N, \tilde{\delta}}\ \mathbb{E}_\mu\Big[\Big(\int_0^s(s-r+1)^{1-2 \tilde{\delta}} X_{r}(1) dr\Big)^{2N}+\Big(\int_0^s(s-r)^{1-2 \tilde{\delta}} X_{r}(1) d r\Big)^{2N}\Big]\nn\\
&\qquad\qquad\qquad +C_{N, \alpha}\ \mathbb{E}_{\mu}\Big[\Big(\int_0^s X_r(1) dr\Big)^{2N}\Big]\Big\}\nn\\
\le &\ (t-s)^{4 \tilde{\delta} N}\Big\{C_{N, \tilde{\delta}}\ \mathbb{E}_\mu\Big[\Big(\sup _{0 \leq r \leq T} X_r(1)\Big)^{2N}\Big]\Big[\Big(\int_0^s(s-r+1)^{1-2 \tilde{\delta}} d r\Big)^{2N}+\Big(\int_0^s(s-r)^{1-2 \tilde{\delta}} d r\Big)^{2N}\Big]\nn\\
&\qquad\qquad\qquad +C_{N, \alpha}\ \mathbb{E}_{\mu}\big[Y_s(1)^{2N}\big]\Big\}\nn\\
\le &\ (t-s)^{4 \tilde{\delta} N}\Big\{C_{N, T,\tilde{\delta}}\ \mathbb{E}_\mu\Big[\big(\sup _{0 \leq r \leq T} X_r(1)\big)^{2N}\Big]+C_{N, \alpha}\ \mathbb{E}_{\mu}\big[Y_s(1)^{2N}\big]\Big\}\nn\\
\le &\ C_{N, T,\alpha, \tilde{\delta}}\ (t-s)^{4 \tilde{\delta} N}.
\end{align}
Combining the above inequality with \eqref{d.0} yields Lemma \ref{bu.lem.b.2}(i).

(ii) When $d=2$, by (\ref{d.1}) and (\ref{d.2}) we get for any $0\le s\le t$ and $\tilde{\delta}\in(0,1)$,
	\begin{align*}
\big|P_{t-r} g_\alpha^a(x)-P_{s-r} g_\alpha^a(x)\big| & \leq 2(t-s)^{\tilde{\delta}}\Big[\tilde{\delta}^{-1}(s-r)^{-\tilde{\delta}}+\frac{1}{\alpha}\Big] \\
& =(t-s)^{\tilde{\delta}}\Big[C_{\tilde{\delta}}(s-r)^{-\tilde{\delta}}+C_\alpha\Big] .
\end{align*}
Notice that for any $0\le s\le T$ and $\tilde{\delta}\in(0,\frac{1}{2})$, there exists a constant $C_{T,\tilde{\delta}}$ such that
\begin{align*}
\int_0^s(s-r)^{-2 \tilde{\delta}} d r\le C_{T,\tilde{\delta}}<\infty.
\end{align*}
Then by similar arguments used in \eqref{d.3} we have for any $0\leq s\leq t\leq T$,
\begin{align*}
&\ \mathbb{E}_\mu\Big[\Big(\int_0^s X_r\big((P_{t-r} g_\alpha^a-P_{s-r} g_\alpha^a)^2\big) d r\Big)^{2N}\Big]\nn \\
\le &\ (t-s)^{4 \tilde{\delta} N}\Big\{C_{N, \tilde{\delta}}\ \mathbb{E}_\mu\Big[\big(\sup _{0 \leq r \leq T} X_r(1)\big)^{2N}\Big]\Big(\int_0^s(s-r)^{-2 \tilde{\delta}} d r\Big)^{2N}+C_{N, \alpha}\ \mathbb{E}_{\mu}\big[Y_s(1)^{2N}\big]\Big\}\nn\\
\le &\ C_{N, T,\alpha, \tilde{\delta}}\ (t-s)^{4 \tilde{\delta} N}.
\end{align*}
This together with \eqref{d.0} gives Lemma \ref{bu.lem.b.2}(ii).

(iii) When $d=3$, by (\ref{d.1}) and (\ref{d.2}) we get for any $0\le s\le t$ and $\tilde{\delta}\in(0,1)$,
	\begin{align}\label{d.1.5}
\big|P_{t-r} g_\alpha^a(x)-P_{s-r} g_\alpha^a(x)\big| & \leq 2(t-s)^{\tilde{\delta}}\Big[\Big(\frac{1}{2}+\tilde{\delta}\Big)^{-1}(s-r)^{-(\frac{1}{2}+\tilde{\delta})}+\frac{1}{\alpha}\Big]\nn \\
& =(t-s)^{\tilde{\delta}}\Big[C_{\tilde{\delta}}(s-r)^{-(\frac{1}{2}+\tilde{\delta})}+C_\alpha\Big] .
\end{align}
It follows from \eqref{e.7.0} and \eqref{e7.24} that
\begin{align}\label{d.1.6}
\big|P_{t-r} g_\alpha^a(x)-P_{s-r} g_\alpha^a(x)\big| 
\leq C_\alpha\big[P_{t-r} g_0^a(x)+P_{s-r} g_0^a(x)\big] \leq C_\alpha g_0^a(x) .
\end{align}
Then we get
\begin{align*}
X_r\big((P_{t-r} g_\alpha^a-P_{s-r} g_\alpha^a)^2\big)
&\le (t-s)^{\tilde{\delta}} C_\alpha\Big[C_{\tilde{\delta}}(s-r)^{-(\frac{1}{2}+\tilde{\delta})}+C_\alpha\Big] X_r(g_0^a)\\
&\le (t-s)^{\tilde{\delta}} C_{\alpha,\tilde{\delta}}\Big[(s-r)^{-(\frac{1}{2}+\tilde{\delta})}+1\Big] X_r(g_0^a).
\end{align*}
For any $\tilde{\delta}\in(0,1/6)$, the interval $(3/2, (\tilde{\delta}+\frac{1}{2})^{-1})\neq\varnothing$. Thus we take
\begin{align}\label{xiu.e.59}
	\gamma_2:=\gamma_2(\tilde{\delta})=\frac{3}{4}+\frac{1}{1+2\tilde{\delta}}.	
\end{align}
It is easily seen that
\begin{align}\label{xiu.e.60}
	\gamma_2\in(3/2, 7/4)\quad\text{and}\quad \frac{\gamma_2}{\gamma_2 -1} \in(7/3, 3).	
\end{align}
Using H\"older's inequality with $p=\gamma_2$ and $q=\gamma_2/(\gamma_2-1)$ yields
\begin{align*}
&\int_0^s X_r\big((P_{t-r} g_\alpha^a-P_{s-r} g_\alpha^a)^2\big) d r\\
\le&\ C_{\alpha, \tilde{\delta}}(t-s)^{\tilde{\delta}} \int_0^s\Big[(s-r)^{-\left(\frac{1}{2}+\tilde{\delta}\right)}+1\Big] X_r\left(g_0^a\right) d r\\
\le&\ C_{\alpha, \tilde{\delta}}(t-s)^{\tilde{\delta}}\Big\{\int_0^s\Big[(s-r)^{-(\frac{1}{2}+\tilde{\delta})}+1\Big]^{\gamma_2} d r\Big\}^{\frac{1}{\gamma_2}}\Big\{\int_0^s \big(X_r(g_0^a)\big)^{\frac{\gamma_2}{\gamma_2-1}} d r\Big\}^{\frac{\gamma_2-1}{\gamma_2}}\\
\le&\ C_{\alpha, \tilde{\delta}}(t-s)^{\tilde{\delta}}\Big\{\int_0^s\Big[(s-r)^{-\gamma_2(\frac{1}{2}+\tilde{\delta})}+1\Big] d r\Big\}^{\frac{1}{\gamma_2}}\Big\{\int_0^s \big(X_r(g_0^a)\big)^{\frac{\gamma_2}{\gamma_2-1}} d r\Big\}^{\frac{\gamma_2-1}{\gamma_2}}.
\end{align*}
Then we have
\begin{align}\label{d.5}
&\ \mathbb{E}_\mu\Big[\Big(\int_0^s X_r\big((P_{t-r} g_\alpha^a-P_{s-r} g_\alpha^a)^2\big) d r\Big)^{2N}\Big]\nn \\
\le &\ C_{N, \alpha, \tilde{\delta}}(t-s)^{2N \tilde{\delta}}\Big(\int_0^s\Big[(s-r)^{-\gamma_2(\frac{1}{2}+\tilde{\delta})}+1\Big] d r\Big)^{\frac{2N}{\gamma_2}}\mathbb{E}_\mu\Big[\Big(\int_0^s \big(X_r(g_0^a)\big)^{\frac{\gamma_2}{\gamma_2-1}} d r\Big)^{\frac{2N(\gamma_2-1)}{\gamma_2}}\Big].
\end{align}
Similarly, we again apply H\"older's inequality with $p=\gamma_2 $ and $q=\gamma_2/(\gamma_2-1)$ to obtain that
\begin{align*}
	\big(X_r(g_0^a)\big)^{\frac{\gamma_2}{\gamma_2-1}}\leq X_r(1)^{\frac{1}{\gamma_2-1}} X_r\big(( g_0^a)^{\frac{\gamma_2}{\gamma_2-1}}\big) .
\end{align*}
Then by Cauchy-Schwarz's inequality, Jensen's inequality, and Lemma \ref{lemx5.1}, for any $0\le s\le T$,
\begin{align}\label{xiu.e.58}
\mathbb{E}_\mu\Big[\Big(\int_0^s \big(X_r(g_0^a)\big)^{\frac{\gamma_2}{\gamma_2-1}} d r\Big)^{\frac{2N(\gamma_2-1)}{\gamma_2}}\Big]
\le&\ \mathbb{E}_\mu\Big\{\Big[\int_0^s X_r(1)^{\frac{1}{\gamma_2-1}} X_r\big(( g_0^a)^{\frac{\gamma_2}{\gamma_2-1}}\big) d r\Big]^{\frac{2N(\gamma_2-1)}{\gamma_2}}\Big\}\nn\\
\le&\ \mathbb{E}_\mu\Big[\Big(\sup _{0 \leq r \leq T} X_r(1)\Big)^{\frac{2N}{\gamma_2}} \Big(Y_T\big((g_0^a)^{\frac{\gamma_2}{\gamma_2-1}}\big)\Big)^{\frac{2N(\gamma_2-1)}{\gamma_2}}\Big]\nn\\
\le&\ \bigg\{ \mathbb{E}_{\mu}\Big[\Big(\sup _{0 \leq r \leq T} X_r(1)\Big)^{\frac{4N}{\gamma_2}} \Big] \mathbb{E}_{\mu}\Big[ \Big(Y_T\big((g_0^a)^{\frac{\gamma_2}{\gamma_2-1}}\big)\Big)^{\frac{4N(\gamma_2-1)}{\gamma_2}}\Big] \bigg\}^{1/2} \nn \\
\le&\ \bigg\{\mathbb{E}_\mu\Big[\Big(\sup _{0 \leq r \leq T} X_r(1)\Big)^{4N} \Big]\bigg\}^{\frac{1}{2\gamma_2}}\bigg\{\mathbb{E}_\mu\Big[\Big(Y_T\big((g_0^a)^{\frac{\gamma_2}{\gamma_2-1}}\big)\Big)^{4N}\Big]\bigg\}^{\frac{\gamma_2-1}{2\gamma_2}} \nn\\
\le&\ C_{N,T,\tilde{\delta}}\Big\{\mathbb{E}_\mu\Big[\Big(Y_T\big((g_0^a)^{\frac{\gamma_2}{\gamma_2-1}}\big)\Big)^{4N}\Big]\Big\}^{\frac{\gamma_2-1}{2\gamma_2}}.
\end{align}
Notice that $g_0^a(x)=\frac{1}{2 \pi|x-a|} \leq \frac{1}{|x-a|}$ and \eqref{de.1.5} holds. It follows by \eqref{xiu.e.60} and Corollary \ref{'cor5.18} that for any $a\in \mathbb{R}^3$ with $|a|\le T$,
\begin{align}\label{d.6}
	\mathbb{E}_\mu\Big[\Big(Y_T\big((g_0^a)^{\frac{\gamma_2}{\gamma_2-1}}\big)\Big)^{4N}\Big] \leq \mathbb{E}_\mu\Big[\Big(\int \frac{1}{|x-a|^{\gamma_2/(\gamma_2-1)}} Y_T(d x)\Big)^{4N}\Big] \leq C_{N,T,\tilde{\delta}}<\infty .
\end{align}
In view of \eqref{xiu.e.59}, for any $\tilde{\delta}\in(0,1/6)$ we have
\begin{align*}
	\gamma_2 \Big(\frac{1}{2}+\tilde{\delta}\Big) = \frac{3}{4}\Big(\frac{1}{2}+\tilde{\delta}\Big) + \frac{1}{2} \in (7/8,1).
\end{align*}
Therefore, for any $0\leq s\leq T$,
\begin{align}\label{d.7}
	\int_0^s\Big[(s-r)^{-\gamma_2(\frac{1}{2}+\tilde{\delta})}+1\Big] d r\le \int_0^T(s-r)^{-\gamma_2(\frac{1}{2}+\tilde{\delta})} d r +T\le C_{T,\tilde{\delta}}<\infty.
\end{align}
Combining \eqref{d.5}, \eqref{xiu.e.58}, \eqref{d.6} and \eqref{d.7} yields for any $\tilde{\delta} \in(0, 1/6)$,
\begin{align*}
\mathbb{E}_\mu\Big[\Big(\int_0^s X_r\big((P_{t-r} g_\alpha^a-P_{s-r} g_\alpha^a)^2\big) d r\Big)^{2N}\Big]\le C_{N, T, \alpha, \tilde{\delta}}(t-s)^{2N \tilde{\delta}}.
\end{align*}
Therefore, Lemma \ref{bu.lem.b.2}(iii) follows from \eqref{d.0} and the above inequality.
$\hfill\blacksquare$\\

\noindent\textbf{Proof of Lemma \ref{bu.lem.b.3}.}
	Let $0\le s\le t\le T$. By Lemma \ref{lem3.3**} and Burkholder-Davis-Gundy's inequality we have
	\begin{align*}
		I&_2^{N,\alpha,a}(s,t)=\mathbb{E}_\mu\Big[\Big|\int_s^t \int P_{t-r} g_\alpha^a(x) d M(r, x)\Big|^{2 N}\Big]\\
		&\le C_N \mathbb{E}_\mu\Big[\Big(\int_s^t X_r\big((P_{t-r} g_\alpha^a)^2\big) d r+\int_s^t d r \int\!\!\!\int g(z, w) P_{t-r} g_\alpha^a(z) P_{t-r} g_\alpha^a(w) X_r(d z) X_r(d w)\Big)^N\Big]\\
		&\le 2^N C_N\Big\{\mathbb{E}_\mu\Big[\Big(\int_s^t X_r\big((P_{t-r} g_\alpha^a)^2\big) dr\Big)^N\Big]+\|g\|_{\infty}^N\mathbb{E}_\mu\Big[\Big(\int_s^t \big(X_r(P_{t-r} g_\alpha^a)\big)^2 dr\Big)^N\Big]\Big\}\\
		&\le C_{N,\alpha}\Big\{\mathbb{E}_\mu\Big[\Big(\int_s^t X_r\big(( g_0^a)^2\big) dr\Big)^N\Big]+\|g\|_{\infty}^N\mathbb{E}_\mu\Big[\Big(\int_s^t \big(X_r(g_0^a)\big)^2 dr\Big)^N\Big]\Big\},
		\end{align*}
where the last inequality follows by \eqref{e.7.0} and \eqref{e7.24}. For any $0\le t\le T$, by arguments similar to those used in \eqref{xiu.e.57} with $g_{\alpha}^a$ replaced by $g_0^a$, we get
\begin{align*}
	\mathbb{E}_\mu\Big[\Big(\int_s^t \big(X_r(g_0^a)\big)^2 dr\Big)^N\Big]
	\le C_{N,T}\Big\{\mathbb{E}_\mu\Big[\Big(\int_s^t X_r\big((g_0^a)^2\big) dr\Big)^{2N}\Big]\Big\}^{1/2}.
\end{align*}
Combining the above two inequalities gives that 
\begin{align*}
	I_2^{N,\alpha,a}(s,t)
	&\le C_{N,\alpha}\Big\{\mathbb{E}_\mu\Big[\Big(\int_s^t X_r\big(( g_0^a)^2\big) dr\Big)^N\Big]+\|g\|_{\infty}^NC_{N,T}\Big\{\mathbb{E}_\mu\Big[\Big(\int_s^t X_r\big((g_0^a)^2\big) dr\Big)^{2N}\Big]\Big\}^{1/2}\Big\}\\
	&\le C_{N,T,\alpha}\Big\{\mathbb{E}_\mu\Big[\Big(\int_s^t X_r\big((g_0^a)^2\big) dr\Big)^{2N}\Big]\Big\}^{1/2},
\end{align*}
where $C_{N, T, \alpha} \in(0, \infty)$ depends on $\|g\|_{\infty}$, and the last inequality follows by applying Caudy-Schwarz's inequality for the first expectation. Then the desired results followed by \eqref{6.53}, \eqref{6.57} and \eqref{6.59}.
$\hfill\blacksquare$

\end{document}